\documentclass[3p,11pt,authoryear]{elsarticle}
\pdfoutput=1
\usepackage[pagewise]{lineno}
\usepackage[hidelinks]{hyperref}

\usepackage{appendix}
\usepackage{bbm}
\usepackage{xcolor}
\usepackage{enumitem}
\usepackage[ruled]{algorithm2e}
\usepackage{adjustbox}
\usepackage{graphicx}
\usepackage{amsfonts}
\usepackage{amssymb}
\usepackage{amsmath}
\usepackage{multirow}
\usepackage{makecell}
\usepackage{rotating}
\usepackage{pdflscape}
\usepackage{arydshln}
\usepackage{etoolbox}
\usepackage{booktabs}
\usepackage{mathtools}
\usepackage{mathrsfs}
\usepackage{subcaption}
\usepackage{comment}
\usepackage{bm}
\usepackage{chngpage}
\usepackage{setspace}
\usepackage{url}
\usepackage{ltablex}
\usepackage{tabularx}
\usepackage[normalem]{ulem}
\allowdisplaybreaks

\makeatletter

\def\ps@pprintTitle{%
	\let\@oddhead\@empty
	\let\@evenhead\@empty
	\def\@oddfoot{\reset@font\hfil\thepage\hfil}
	\let\@evenfoot\@oddfoot
}
\makeatother

\makeatletter
\def\@author#1{\g@addto@macro\elsauthors{\normalsize%
		\def\baselinestretch{1}%
		\upshape\authorsep#1\unskip\textsuperscript{%
			\ifx\@fnmark\@empty\else\unskip\sep\@fnmark\let\sep=,\fi
			\ifx\@corref\@empty\else\unskip\sep\@corref\let\sep=,\fi
		}%
		\def\authorsep{\unskip,\space}%
		\global\let\@fnmark\@empty
		\global\let\@corref\@empty 
		\global\let\sep\@empty}%
	\@eadauthor={#1}
}
\makeatother

\begin{document}
\begin{frontmatter}


\title{Planning ride-pooling services with detour restrictions for spatially heterogeneous demand: A multi-zone queuing network approach}


\author{Yining Liu}
\author{Yanfeng Ouyang\corref{correspondingauthor}}
\address{Department of Civil and Environmental Engineering, University of Illinois at Urbana-Champaign, Urbana, IL 61801, USA}
\cortext[correspondingauthor]{Corresponding author. Tel: +1 217 333 9858}
\ead{yfouyang@illinois.edu}

\begin{abstract}
This study presents a multi-zone queuing network model for steady-state ride-pooling operations that serve heterogeneous demand, and then builds upon this model to optimize the design of ride-pooling services. Spatial heterogeneity is addressed by partitioning the study region into a set of relatively homogeneous zones, and a set of criteria are imposed to avoid significant detours among matched passengers. A generalized multi-zone queuing network model is then developed to describe how vehicles' states transition within each zone and across neighboring zones, and how passengers are served by idle or partially occupied vehicles. A large system of equations is constructed based on the queuing network model to analytically evaluate steady-state system performance. Then, we formulate a constrained nonlinear program to optimize the design of ride-pooling services, such as zone-level vehicle deployment, vehicle routing paths, and vehicle rebalancing operations. A customized solution approach is also proposed to decompose and solve the optimization problem. The proposed model and solution approach are applied to a hypothetical case and a real-world Chicago case study, so as to demonstrate their applicability and to draw insights. Agent-based simulations are also used to corroborate results from the proposed analytical model. These numerical examples not only reveal interesting insights on how ride-pooling services serve heterogeneous demand, but also highlight the importance of addressing demand heterogeneity when designing ride-pooling services.
\end{abstract}


\end{frontmatter}
\thispagestyle{empty}
\nolinenumbers

\section{Introduction}
The advancement of wireless communication technologies has fueled the booming growth of a spectrum of demand-responsive services offered by transportation network companies (TNCs), e.g., the driver may have a travel need \citep{wang2018stable,wang2019ridesourcing} or provide dedicated chauffeur service, while the passenger may share the ride with the driver only \citep{chen2019optimal,ozkan2020dynamic} or with one or more other passengers \citep{Daganzo2019Paper,bei2018algorithms}. 

In this paper, we focus on one of the popular services in the current market, chauffeured ride-pooling 
(also referred to as ride-splitting or shared taxi), e.g., UberPool, Lyft Line, and Didi ExpressPool. It serves multiple passengers with similar spatiotemporal characteristics in one ride. Individual participants may detour from their shortest paths, but they can potentially reduce their own costs by sharing the expenses with matched participants. At the societal level, properly operated ride-pooling services are expected to increase vehicle occupancy, alleviate congestion, and reduce emissions \citep{caulfield2009estimating,chan2012ridesharing,furuhata2013ridesharing,li2016empirical,yu2017environmental,alisoltani2021can}. 


Ride-pooling services face many similar challenges as the non-shared services (such as conventional taxi and crowd-sourcing services); e.g., the travel demand in cities usually exhibits high levels of spatial heterogeneity, and the imbalance of trip origins and destinations requires a large number of service vehicles to be systematically and regularly relocated. Failure to properly plan and rebalance service vehicles often leads to inefficient operating modes (e.g., when most of the vehicles are trapped in long-distance deadheading trips), also known as the wild goose chase (WGC) phenomenon \citep{daganzo2010public,castillo2017surge,Daganzo2019Paper}.
These operational challenges have been extensively studied in the literature, mostly with a focus on non-shared services. For example, pricing strategies have been thoroughly studied as an instrument to manage demand and supply: on the demand side, the platform uses OD-based or even path-based prices to induce or suppress trip demand in certain areas \citep{castillo2017surge, lei2019path}; on the supply side, drivers are incentivized to relocate to desired areas \citep{bimpikis2019spatial,zhu2021mean,chen2021spatial,li2021spatial,afifah2022spatial}. More recently, \citet{xu2021generalized} develops a fluid model to capture the system dynamics of ride-sourcing services, and uses it as the basis to optimize TNC pricing, matching, and vehicle dispatching decisions. 
In addition, various strategies are developed to reduce unproductive vehicle deadheading by matching travel requests and drivers more effectively, e.g., by limiting geometric matching within an adaptive radius \citep{xu2020supply}, using block matching \citep{feng2022approximating}, or swapping passenger-vehicle assignments dynamically \citep{shen4264716dynamic, ouyang2023measurement}.

However, ride-pooling services also face new challenges, which come exactly from the sharing of rides -- a new passenger may be picked up by a vehicle already serving another passenger. These passengers have different origins, destinations, or desired departure times. Pooling them into a shared vehicle most likely will create extra waiting and/or detours, and these inconveniences in turn directly affect these passengers' willingness to use such a service \citep{krueger2016preferences,chen2017understanding,storch2021incentive}. Hence, it is critical to design effective passenger matching algorithms or strategies to avoid significant waiting or detours. In so doing, a large body of studies have focused on optimizing ride-pooling services, mainly at the operational level, e.g., dynamic matching algorithms \citep{santi2014quantifying,alonso2017demand,di2013optimization,bei2018algorithms,zeng2020exploring}, designing meeting points \citep{stiglic2015benefits,chen2019ride}, and demand management and pricing \citep{ke2020pricing,jacob2021ride,zhu2020analysis}. Readers are referred to \citet{tafreshian2020frontiers} for a comprehensive review of recent efforts. 
  
In contrast, there are relatively few models that address the planning of ride-pooling systems. Recently, \citet{Daganzo2019Paper} proposed an overarching modeling framework that can cover a range of demand-responsive services (e.g., taxi, ride-pooling, and dial-a-ride), and analytically connect (in closed-form formulas) TNCs' fleet deployment decisions to passengers' expected steady-state travel experience. Later, this modeling framework is generalized in a series of studies to address ride-pooling-related challenges, e.g., restricting the detour and waiting time of matched passengers \citep{Daganzo2020,ouyang2021performance}, and integrating the ride-pooling service into transit network design \citep{Liu2021}. 

These above-mentioned studies on strategic planning of ride-pooling services have generally assumed spatially homogeneous demand, probably because it is far more complex to predict how passengers from heterogeneous neighborhoods share rides: now, partially-occupied vehicles passing a certain neighborhood can join the local idle vehicles in serving arriving passengers, and hence how these vehicles are routed across neighborhoods, and how they are probabilistically matched with passengers under waiting/detour limits, would make a significant difference. The closest literature we can find address operational level decisions in this regard. For example, path-based pricing schemes are proposed, based on dynamic programming, to optimize TNCs' pricing and vehicle path decisions for ride-pooling \citep{lei2019path} and dial-a-ride services \citep{shen2021path}. Yet, to the authors' best knowledge, the problem of strategically planning ride-pooling services to serve spatially heterogeneous demand has not been covered in the existing literature.

In light of the challenges mentioned above, this paper aims to build a steady-state ride-pooling service model for spatially heterogeneous demand, and use it to optimize the design of the ride-pooling service. 
The TNC needs to determine the spatial deployment of vehicles, their approximate 
routing plan during service, the matching among vehicles and passengers, and the strategic rebalancing of idle vehicles. We first represent the heterogeneous service 
region by a set of relatively homogeneous zones, and propose a set of detour/waiting criteria for the ride-pooling service to match different types of passenger trips within and across zones. Then, we extend the steady-state aspatial queuing network model in \citet{Daganzo2019Paper} and \citet{Daganzo2020} into a multi-zone version, which contains a set of vehicle states (characterized by a vehicle's workload and location) and a set of vehicle state transitions inside a zone (driven by passenger pickup, drop-off, and assignment events) and across zones (driven by vehicle crossing zone boundaries). A large system of nonlinear equations is formulated to connect the vehicle count in each state and the state transition flow rates, based on vehicle flow conservation within and across zones, vehicle-passenger matching criteria, and the expected duration of each vehicle state. The solution to the system of equations can be used to estimate expected system performance (such as a passenger's expected travel time). Hence, these equations are further embedded into a constrained nonlinear program so that TNCs can make optimal vehicle deployment, routing, and rebalancing decisions to minimize the average system-wide cost. A customized solution approach is also proposed to decompose and effectively solve the optimization problem. 

The proposed model and solution approach are implemented in a series of hypothetical cases with varying levels of demand heterogeneity and zone sizes, to demonstrate their applicability and to cast insights into the optimal system characteristics. Idle vehicles and partially-occupied vehicles are found to complement each other in providing services, and for this exact reason, the relative locations of high- or low-demand zones tend to influence the volume of passing vehicles as well as idle vehicle deployments. In general, when demand is highly heterogeneous, TNCs can deploy fewer vehicles because it is easier to match riders under more concentrated demand, but these vehicles also must be engaged in more rebalancing activities. The impacts of zone size choices on system performance are also discussed. These model results are found to closely match the observations from agent-based simulations. 
 In addition, we also apply the proposed model to design ride-pooling services in downtown and suburban Chicago, where trip rates across zones vary by orders of magnitude, to demonstrate its applicability to real-world applications. The optimal service design from our model is found to be significantly better than its benchmark counterpart (obtained under homogeneous demand assumption) when the demand is actually heterogeneous, because: (i) it achieves a significant reduction in system-wide cost (especially when the passengers' value of time is low), and (ii) it avoids the WGC (which arises when the benchmark design is used). This observation highlights the importance of properly addressing demand heterogeneity when planning ride-pooling services.

The remainder of this paper is organized as follows. Section \ref{sec:methodology} describes the multi-zone version of the queuing network model, including the matching criteria, the system of equations, the performance metrics, the constrained non-linear program formulation, and the suggested solution approach. Section \ref{sec:numerical} presents the series of numerical experiments. 
Section \ref{sec:conclusion} concludes the paper and discusses possible future extensions.

\section{Methodology}\label{sec:methodology}
\subsection{Model Setting}
\begin{figure}[t]
     \centering
     \begin{subfigure}[b]{\textwidth}
        \centering
        \includegraphics[width=0.95\textwidth]{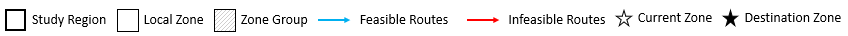}
     \end{subfigure}
     \begin{subfigure}[b]{0.47\textwidth}
         \centering
         \includegraphics[width=\textwidth]{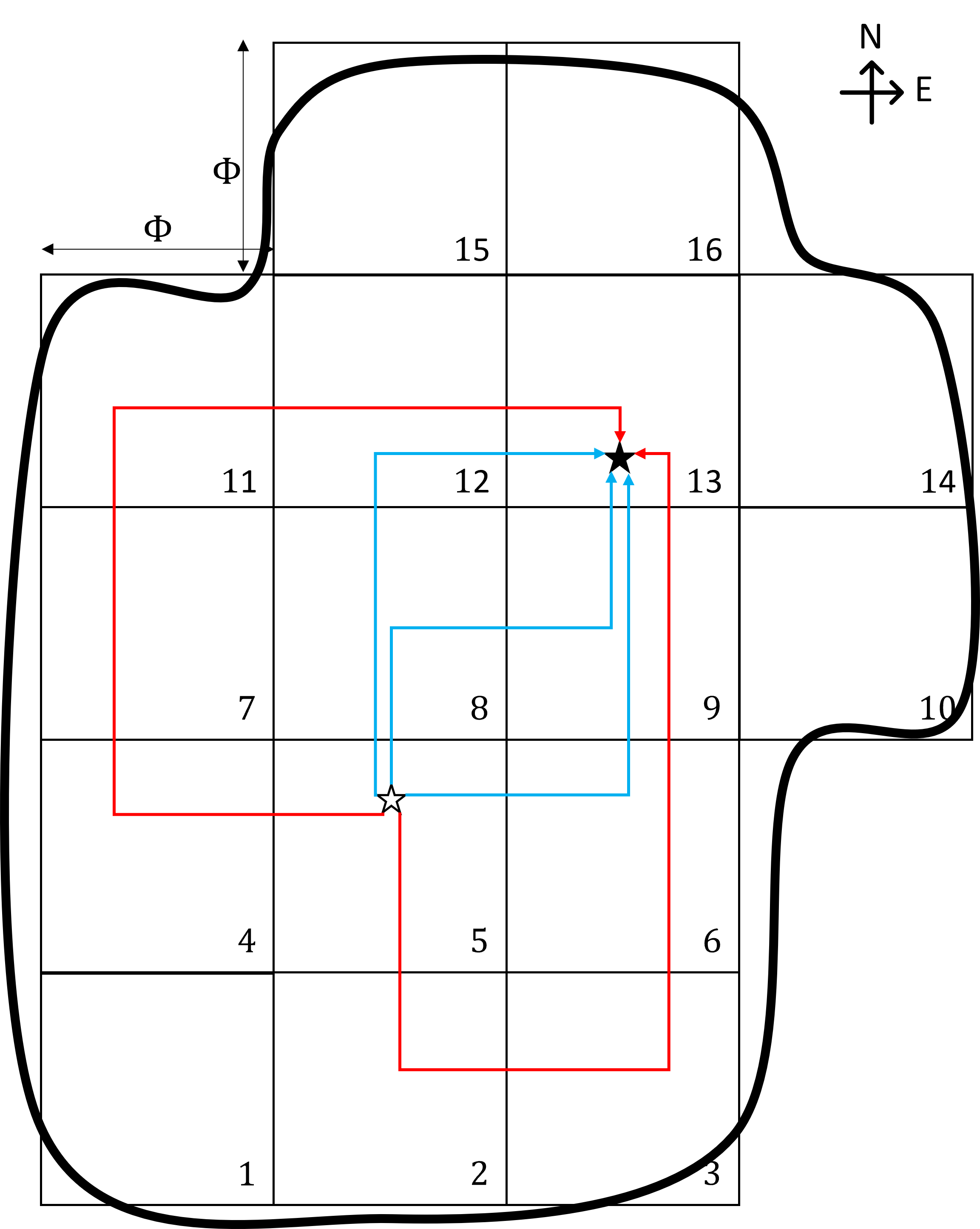}
         \caption{Study region and local zones, with route examples.}
         \label{fig:study region 1}
     \end{subfigure}
     \begin{subfigure}[b]{0.47\textwidth}
         \centering
         \includegraphics[width=\textwidth]{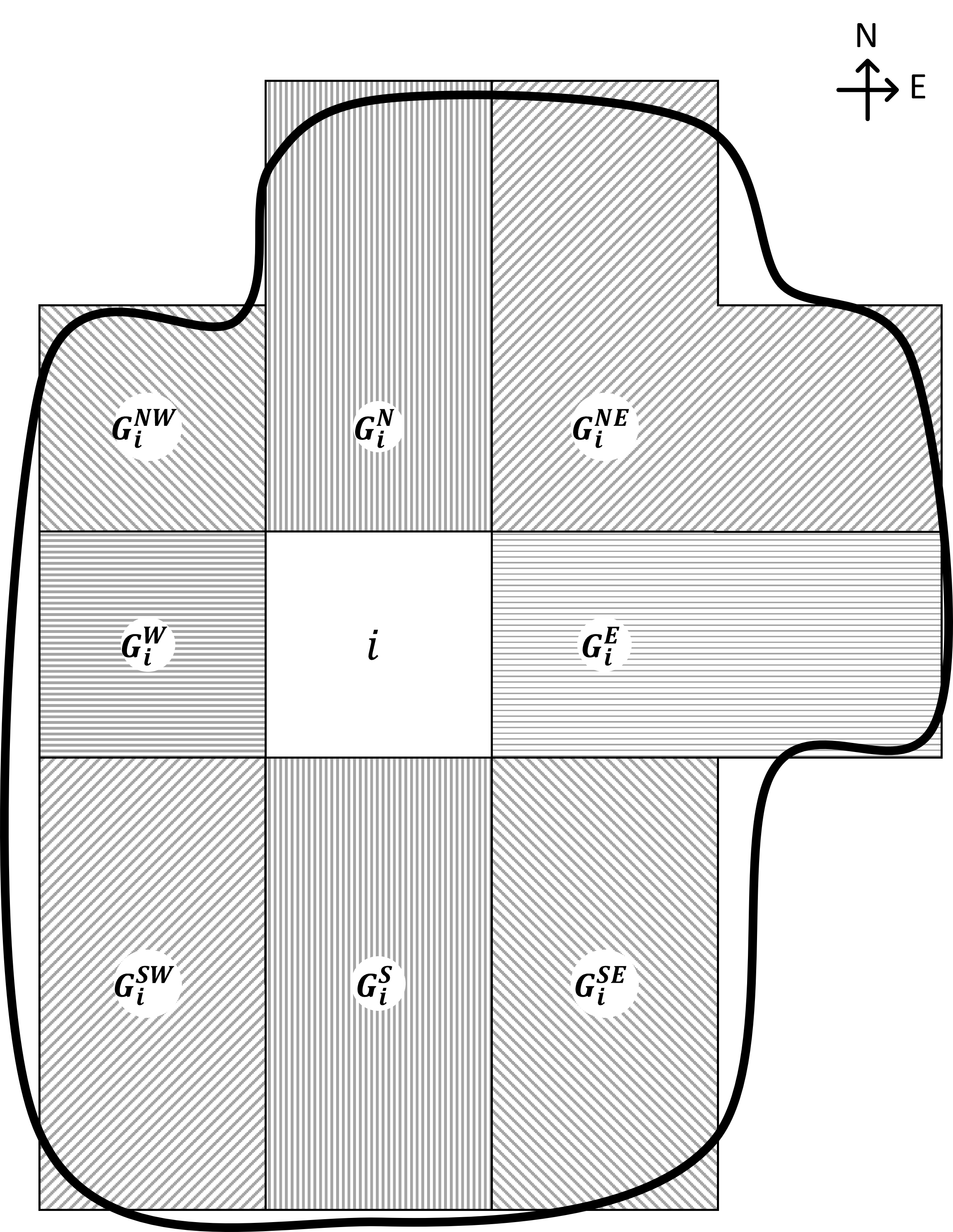}
         \caption{Subgroup of zones relative to zone $i$.}
         \label{fig:study region 2}
     \end{subfigure}
    \caption{Illustration of study region, local zones, and zone groups.}
    \label{fig:study region}
\end{figure}

We assume that a transportation network company (TNC) needs to provide ride-pooling services to inelastic and spatially heterogeneous demand in a study region. A dense grid of streets covers the region, and service vehicles can travel on the street network along the E-W or N-S direction. Here, we let E, NE, N, NW, W, SW, S, and SE stand for East, Northeast, North, Northwest, West, Southwest, South, and Southeast directions, respectively. The region is simply connected\footnote{If the study region is not simply connected (i.e., containing ``holes''), the general modeling framework is still applicable, but the definition of vehicle travel directions and the estimate of travel distance should be carefully updated. This topic will be left for future research.} and can be approximated by a set of relatively homogeneous square zones,\footnote{In this study, we assume that the partition of the region is given. Methods for partitioning spatially heterogeneous areas have been proposed in various contexts; e.g., for spatial quantization \citep{du1999centroidal}, for developing macroscopic fundamental diagrams \citep{ji2012spatial}, for school districting \citep{ferland1990decision}, for shape control \citep{ouyang2006discretization,ouyang2007design}, and for network districting \citep{xie2016railroad}. Additionally, the square zone shape is chosen in this paper for simplicity; optimal zone partition shapes under a few cost metrics are discussed in \citet{xie2015optimal}.} each with side length $\Phi$ [km]. Figure \ref{fig:study region 1} shows an example, in which the study region is approximated by a grid of 16 zones. Let $\mathcal{K}$ denote the set of zones, i.e., $\mathcal{K}=\{1,2,\cdots,|\mathcal{K}|\}$. The shortest distance between the centroids of zones $i$ and $j$ via the street network is denoted by $L_{ij}$ [km]; e.g., $L_{5, 13} = L_{13, 5} = 3\Phi$ in Figure \ref{fig:study region 1}. We define $A_i^r$ to be the index of the adjacent zone on the $r$ side of zone $i$, where $r\in\{\text{E}, \text{N}, \text{W}, \text{S}\}$, and $A_i^r = 0$ (i.e., index for a dummy zone, added for modeling convenience) if no such a zone exists. In Figure \ref{fig:study region 1}, for zone $6$, $A_6^\text{E} = 0$, $A_6^\text{N} = 9$, $A_6^\text{W} = 5$, and $A_6^\text{S} = 3$. Furthermore, zones are grouped based on their relative positions to each zone $i\in\mathcal{K}$. We let $G_i^r$ denote the indices of zones in the $r$ direction of zone $i$, where $r\in \{\text{E}, \text{NE}, \text{N}, \text{NW}, \text{W}, \text{SW}, \text{S}, \text{SE}\}$. For convenience, we also define the set of zones that are aligned along an E-W or N-S line of zone $i$, i.e., $G^{+}_i = G_i^\text{N} \cup G_i^\text{W} \cup G_i^\text{S} \cup G_i^\text{E}$; and all the ``diagonal" ones, 
i.e., $G^{\times}_i = G_i^\text{NE} \cup G_i^\text{NW} \cup G_i^\text{SW} \cup G_i^\text{SE} = \mathcal{K} \setminus (G^{+}_i\cup\{i\})$
, respectively. Figure \ref{fig:study region 2} shows an illustration, where for example $G_8^{\text{E}} = \{9, 10\}$, $G_8^{\text{NE}} = \{13, 14, 16\}$, $G^{+}_8 = \{2, 5, 7, 9, 10, 12, 15\}$, and $G^{\times}_8 = \{1, 3, 4, 6, 11, 13, 14, 16\}$. A table of key notation is summarized in \ref{sec:appendix notation}.

The trips from zone $i \in \mathcal{K}$ to zone $j \in \mathcal{K}$ are generated from a homogeneous Poisson process at rate $\lambda_{ij}$ [trip/hr]. Within each zone, the origins and destinations (OD) of all its associated trips are distributed uniformly over space and time. All travelers have the same value of time $\beta$ [\$/trip-hr]. The TNC uses a fleet of service vehicles driven by fully compliant drivers (e.g., salaried drivers or robo-taxis). Each vehicle travels at a constant speed $v$ [km/h], and the time for passenger boarding and alighting from a vehicle is negligible. The total hourly vehicle cost (including driver salary, prorated vehicle acquisition, and operation costs) is denoted as $\gamma$ [\$/veh-hr]. 

The TNC needs to make a range of planning decisions. First, it must determine the steady-state distribution of vehicles (especially the idle vehicles) across zones. Second, it must plan the zone-level paths for vehicles traveling between every OD zone pair\footnote{The zone-level path decisions can be easily enforced by the platform in case of robo-taxis or when drivers are fully compliant (e.g., see \citet{lei2019path}). In case of more independent human drivers, the platform could suggest these routes on the app’s navigation system and encourage the drivers to follow them (possibly with incentives).} --- since 
the travel demand is heterogeneous across zones, a vehicle's path across zones may affect its probability to encounter a passenger en-route. Here, we assume that the TNC platform only considers vehicle zone-level paths that do 
not have zone-level detours \citep{lei2019path}; i.e., no vehicle service path can simultaneously contain zone-level movements in both north and south (or east and west) directions. For example, in Figure \ref{fig:study region 1}, the blue zone-level paths are feasible, but the red ones are not allowed. Third, it must have a plan to rebalance idle vehicles between zones because the incoming and outgoing trip generation rates of a zone could differ.

In this study, we assume that the TNC platform adopts the following operating strategies.\footnote{These strategies aim at reducing passengers' out-of-vehicle waiting time and limiting service detours. The vehicle capacity is also limited to two purely for simplicity of graphical illustration and explanation of the multi-zone queuing network model. Larger vehicle capacities and other operation strategies, as those in \citet{Daganzo2019Paper} and \citet{Daganzo2020}, can be accommodated similarly without affecting the structure and validity of the overarching modeling framework.} A vehicle can serve at most two passengers at any time (either for pickup or drop-off); in addition, no new passenger will be assigned to a vehicle already en-route to pick up another passenger (otherwise, the passengers may experience an excessive wait). As a result, the vehicles suitable for a new passenger with origin in zone $i$ include idle vehicles and those partially occupied (i.e., with one passenger on-board) inside zone $i$; 
among them, the nearest one is assigned instantly.\footnote{We assume that the assignment will not be changed in the future. This simplification ignores the possibility that a future caller (or suitable vehicle) may show up at a location closer to the assigned vehicle (or passenger) when the assigned vehicle is en-route to pick up the assigned passenger. 
}

As in \citet{Daganzo2019Book} and \citet{Daganzo2020}, the new passenger is referred to as a ``caller" (i.e., one who calls for service), and the single onboard passenger of a vehicle is referred to as a ``seeker" (i.e., one who seeks another rider), respectively. If the caller is assigned to an idle vehicle, the idle vehicle immediately moves toward the caller's origin, and the caller becomes a seeker upon pickup. 
Otherwise, if the caller is assigned to a vehicle with a seeker, we say that the caller and the seeker are ``matched,'' and the vehicle with the seeker will detour from its original path to pick up the assigned caller before dropping off the seeker; upon pickup, the vehicle has two passengers onboard, and will next deliver the passenger with a closer destination (which could be within the current zone $i$ or in another zone $j\in\mathcal{K}\setminus\{i\}$). Note that with the above rules, if an idle vehicle is assigned to a caller (whose origin is within the same zone), the deadheading trip for pickup cannot cross zone borders.


\begin{figure}[t]
     \centering
     \begin{subfigure}[b]{0.47\textwidth}
         \centering
         \includegraphics[width=0.95\textwidth]{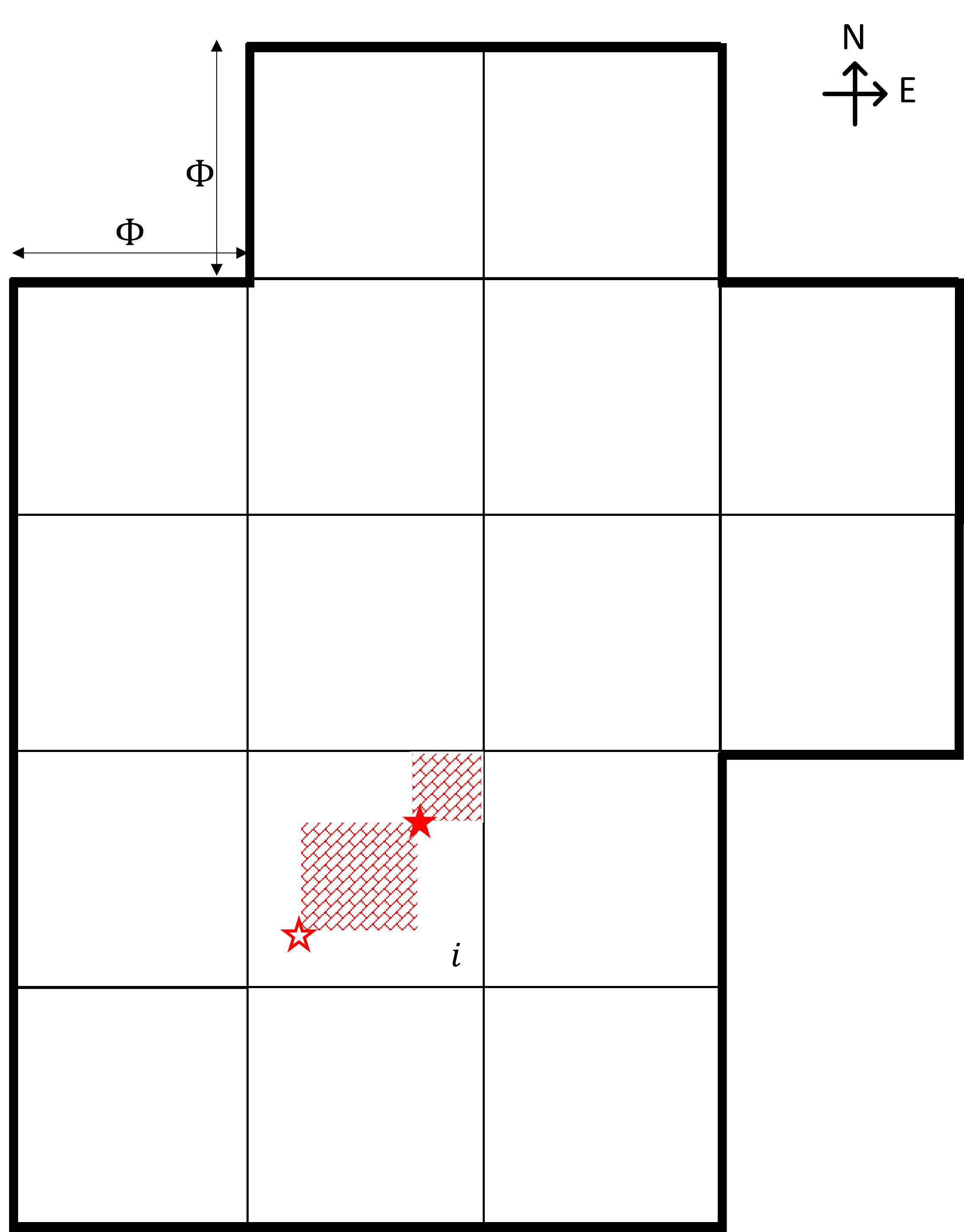}
         \caption{Seeker destination in zone $i$ vs. intra-zonal caller}
         \label{fig:matching example intra intra}
     \end{subfigure}
     \begin{subfigure}[b]{0.47\textwidth}
         \centering
         \includegraphics[width=0.95\textwidth]{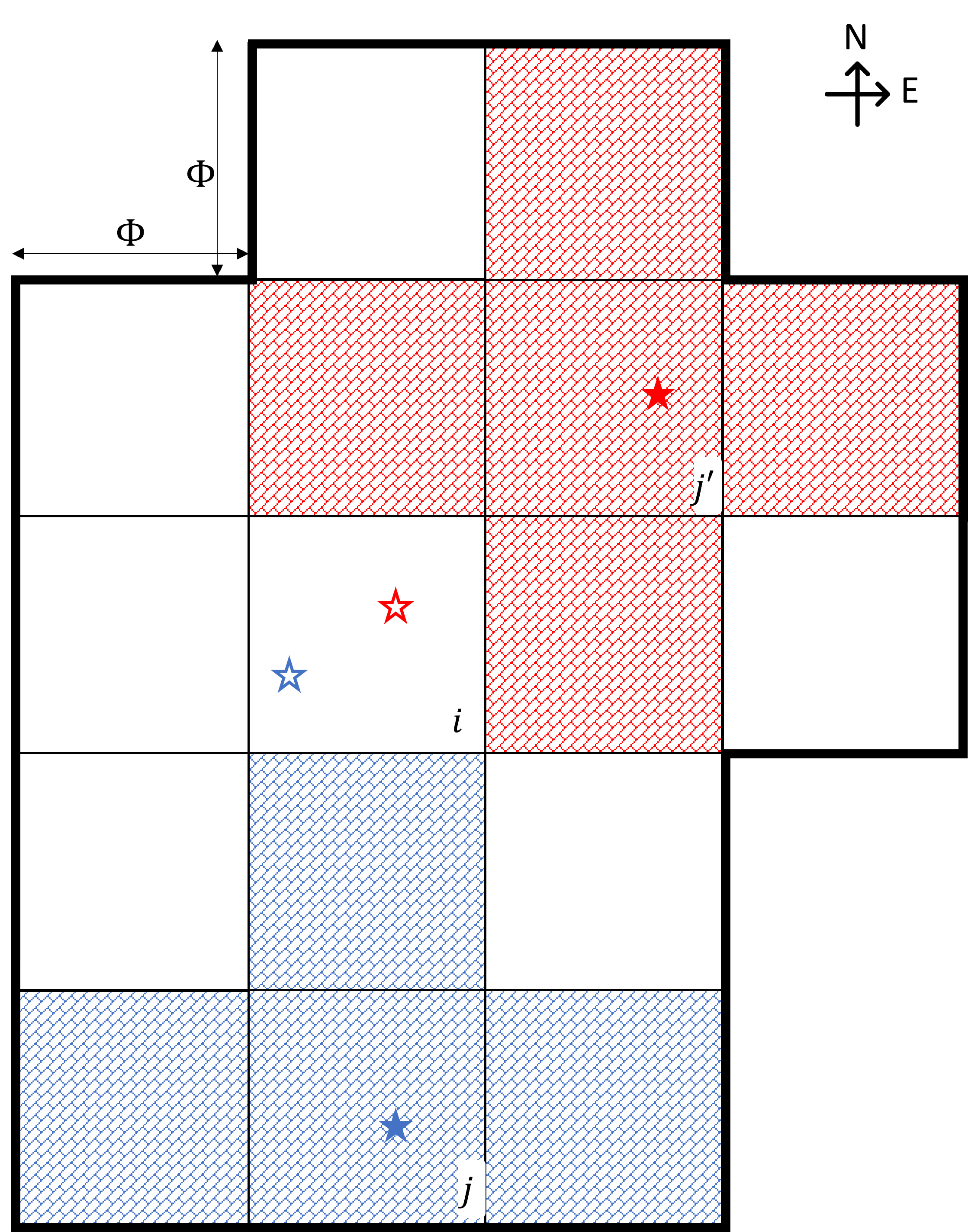}
         \caption{Seeker destination out of zone $i$ vs. inter-zonal caller}
         \label{fig:matching example inter inter}
     \end{subfigure}
     \begin{subfigure}[b]{\textwidth}
         \centering
         \includegraphics[width=\textwidth]{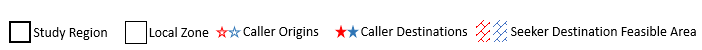}
     \end{subfigure}
    \caption{Feasible area for a seeker's destination under detour restriction rules given a caller's OD. (It should be noted that an intra-zonal caller's OD is continuously distributed in a local zone, and it could travel in NE, NW, SW, and SE directions; the destination zone of an inter-zonal caller could be in the E, NE, N, NW, W, SW, S, and SE directions. In the following, as an illustration, we present the detour requirements for one representative direction of an intra-zonal caller (NE) and two representative directions of an inter-zonal caller (S and NE). The matching criteria for other directions shall be obvious from symmetry.)}
    \label{fig:matching example 1}
\end{figure}

\begin{figure}[t]
     \centering
     \begin{subfigure}[b]{0.47\textwidth}
         \centering
         \includegraphics[width=0.9\textwidth]{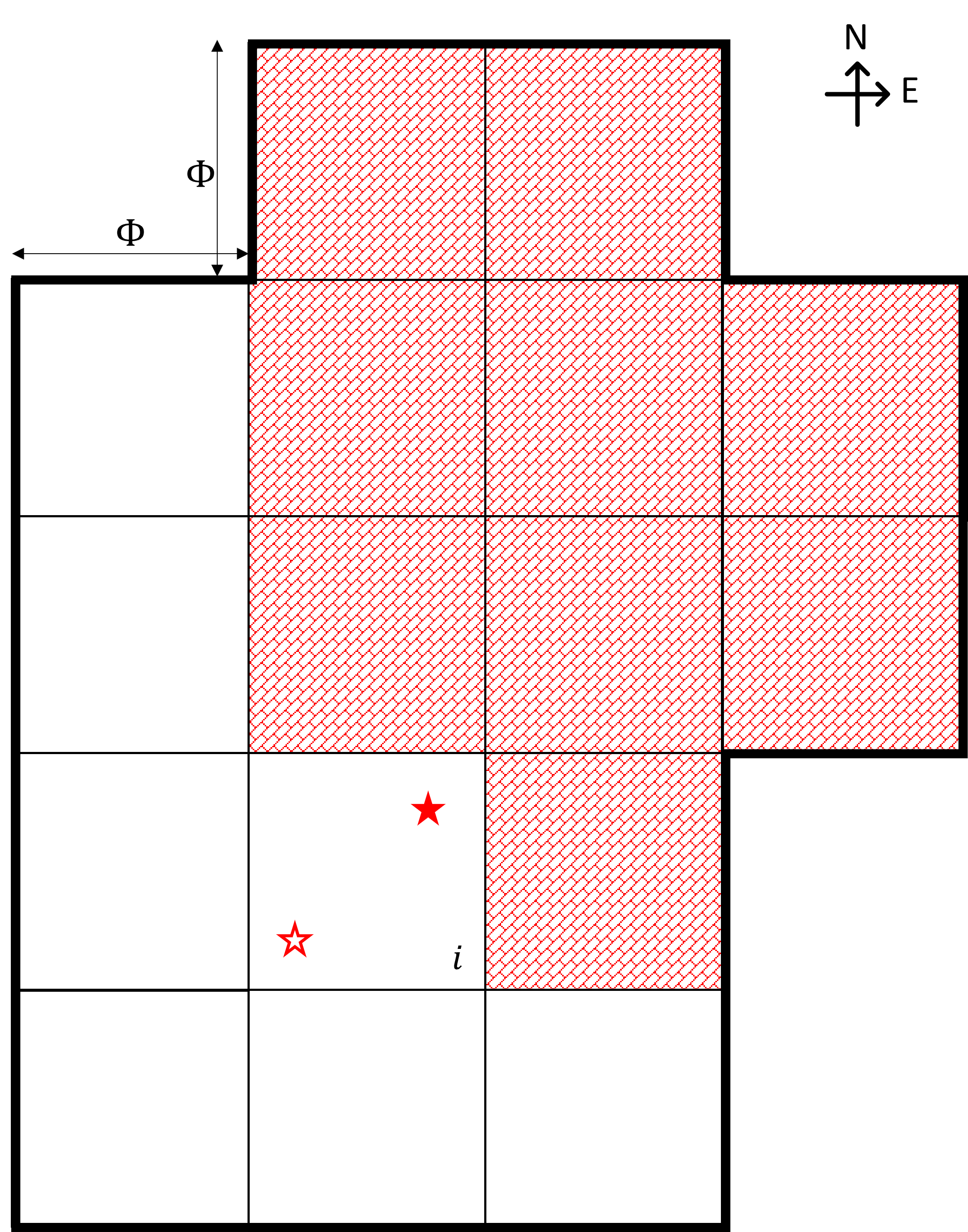}
         \caption{Seeker destination out of zone $i$ vs. intra-zonal caller}
         \label{fig:matching example inter intra}
     \end{subfigure}
     \begin{subfigure}[b]{0.47\textwidth}
         \centering
         \includegraphics[width=0.9\textwidth]{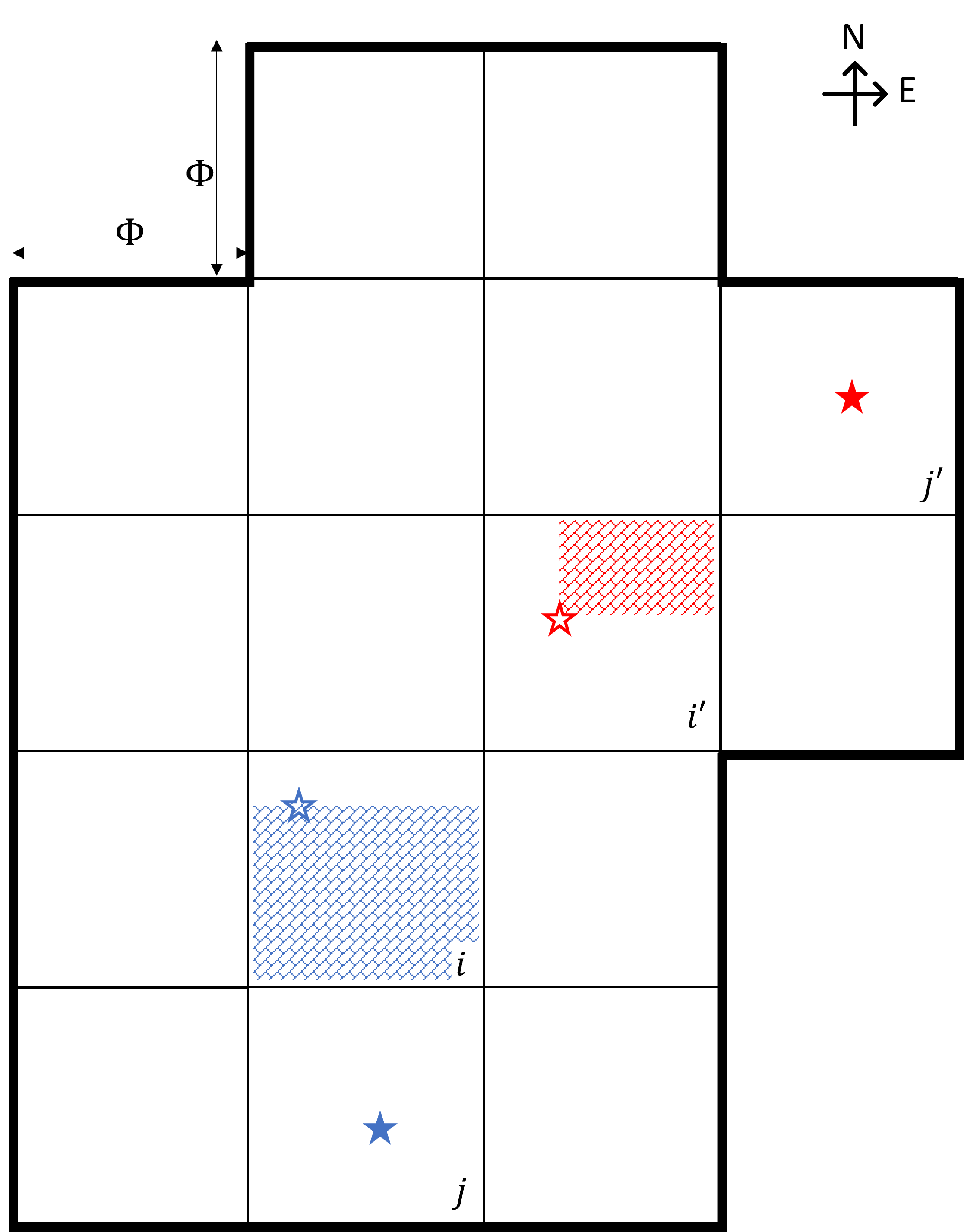}
         \caption{Seeker destination in zone $i$ vs. inter-zonal caller}
         \label{fig:matching example intra inter}
     \end{subfigure}
     \begin{subfigure}[b]{\textwidth}
         \centering
         \includegraphics[width=\textwidth]{images/Matching_legend.png}
     \end{subfigure}
    \caption{Feasible area for seeker's destination under detour restriction rules given a caller's OD (Part II).}
    \label{fig:matching example 2}
\end{figure}

However, not all nearby vehicles with a seeker are suitable to serve a caller, since long detours will occur if the destination of the caller is far away from that of the seeker. We assume that the TNC imposes case-dependent limits on possible detours from a caller's origin to the destinations of both this caller and a seeker.\footnote{No distance limit is imposed on the trip from a vehicle (either an idle vehicle or one with a seeker) to this caller, because we want to ensure that all callers receive service. This assumption can be relaxed, however, without affecting the overall modeling framework.} These rules, as summarized next, yield the seeker destination feasible area for any caller; i.e., if a seeker's destination is within this area, it is eligible to be matched with this caller.  
\begin{enumerate}[leftmargin=*, label={Case \arabic*:}, ref={Case \arabic*}]
    \item \label{case:1} When the seeker's destination is within the current zone of the vehicle, $i$, and the caller's trip is intra-zonal: we do not allow any detour because both passengers' remaining trips are short; see the shaded area in Figure \ref{fig:matching example intra intra} for an illustration of feasible area for the seeker's destination. 
    \item \label{case:2} When the seeker's destination is outside of zone $i$, and the caller's trip is inter-zonal toward a different zone $j$: we only allow detours inside local zones, but not at the zone level. For any given caller's OD zone pair $(i,j)$, this feasible area is fixed. Figure \ref{fig:matching example inter inter}, for example, illustrates the seeker destination feasible areas for a caller going in the S direction (shaded blue) and one going in the NE direction (shaded red), respectively.
    \item \label{case:3} When the seeker's destination is outside of zone $i$, but the caller's trip is intra-zonal: we allow local detours but require the travel directions of the seeker and the caller to be aligned. This rule avoids trapping the seeker in a local zone due to (possible) repeated pickups/deliveries of intra-zonal callers. As shown in Figure \ref{fig:matching example inter intra}, if the caller travels in the NE direction, it can only be matched with a seeker whose destination zone $j$ is to the N, E, or NE directions of the current zone $i$; i.e., $j\in G_i^{\text{N}}\cup G_i^{\text{NE}}\cup G_i^{\text{E}}$.
    \item \label{case:4} When the seeker's destination is within zone $i$, and the caller's trip is inter-zonal: the same rule as that for \ref{case:3} is applied to reduce the caller's detours after it is picked up by the vehicle. Under this rule, the vehicle will travel in a direction aligned with the caller's destination while delivering the seeker.
    The feasible area for the seeker's destination is the part of the current zone that is to the corresponding direction of the caller's destination; see the shaded areas in Figure \ref{fig:matching example intra inter} for two examples, respectively, when the caller goes south (blue) and northeast (red). 
\end{enumerate}

Before building the model, we further define a few sets that can help characterize the relative positions of the zones. First, we let ${\mathcal{V}}_{ij}\subseteq \{A^\text{E}_i, A^\text{N}_i, A^\text{W}_i, A^\text{S}_i\}$ as the set of zones that a vehicle in zone $i$ heading to zone $j$ may enter immediately after leaving zone $i$ without zone-level detours. For example, when $j\in G_i^{\text{NE}}$, we have $\mathcal{V}_{ij} = \{A_i^{\text{N}}, A_i^{\text{E}}\}$; when $j\in G_i^{\text{N}}$, we have $\mathcal{V}_{ij} = \{A_i^{\text{N}}\}$. 
To simplify the expression, we define $U^r$ as the set of two directions adjacent to direction $r$, i.e.,

\begin{equation}\label{eq:side directions}
\begin{aligned}
    & U^{\text{E}} = \{\text{SE}, \text{NE}\},\ U^{\text{NE}} = \{\text{E}, \text{N}\},\ U^{\text{N}} = \{\text{NE}, \text{NW}\},\ U^{\text{NW}} = \{\text{N}, \text{W}\},\\
    & U^{\text{W}} = \{\text{NW}, \text{SW}\},\ U^{\text{SW}} = \{\text{W}, \text{S}\},\ U^{\text{S}} = \{\text{SW}, \text{SE}\},\ U^{\text{SE}} = \{\text{S}, \text{E}\},
\end{aligned}
\end{equation}
and therefore, we can express the feasible travel directions as follows,
\begin{equation}\label{eq:feasible directions}
\begin{aligned}
    & \mathcal{V}_{ij} = \{A_i^r\},\quad &&\forall i\in \mathcal{K},\ j\in G_i^r,\ r\in \{\text{E}, \text{N}, \text{W}, \text{S}\}\\
    & \mathcal{V}_{ij} = \{A_i^{r'}|r'\in U^r\}.\quad &&\forall i\in \mathcal{K},\ j\in G_i^r,\ r\in \{\text{NE}, \text{NW}, \text{SW}, \text{SE}\}\\
\end{aligned}
\end{equation}

Next, we define $\Omega_{ij}$ to be the set of zones in the seeker's destination feasible area when the caller is inter-zonal traveling from zone $i$ to zone $j\in\mathcal{K}\setminus\{i\}$, and $\tilde{\Omega}_{ij} = \{k\in\Omega_{ij}|L_{ik}\geq L_{ij}\}$ the subset with zones no closer to zone $i$ than zone $j$. For example, if a caller travels from zone $5$ to zone $13$ in Figure \ref{fig:study region 1}, $\Omega_{5, 13} = \{6, 8, 9, 12, 13, 14, 16\}$, and $\tilde{\Omega}_{5, 13} = \{13, 14, 16\}$. Similarly, we define $\Omega^r_{ii}$ to be the set of zones in the seeker destination feasible area for an intra-zonal caller traveling in the $r$ direction, $r\in\{\text{NE}, \text{NW}, \text{SW}, \text{SE}\}$:

\begin{equation}\label{eq: destination feasible for intra caller}
\begin{aligned}
    & \Omega^{\text{NE}}_{ii} = G_i^{\text{E}}\cup G_i^{\text{NE}}\cup G_i^{\text{N}},\quad \Omega^{\text{NW}}_{ii} = G_i^{\text{N}}\cup G_i^{\text{NW}}\cup G_i^{\text{W}},\\ 
    & \Omega^{\text{SW}}_{ii} = G_i^{\text{W}}\cup G_i^{\text{SW}}\cup G_i^{\text{S}},\quad \Omega^{\text{SE}}_{ii} = G_i^{\text{S}}\cup G_i^{\text{SE}}\cup G_i^{\text{E}}.
    \end{aligned}
\end{equation}
Finally, we assume the following routing rules. A vehicle always takes a shortest path towards its next destination either when this destination is inside its current zone (e.g., to pick up or deliver a passenger), or when the vehicle is being rebalanced across zones. If a vehicle is delivering a single onboard inter-zonal seeker along the suggested zone-level path, it travels straight to leave the current zone, avoids any unnecessary turns, and takes a shortest path to drop off the passenger after entering the passenger's destination zone. If the zone-level path mandates that a vehicle makes a turn inside a zone, the turn occurs at a random location in anticipation of possible ride matches. 
If a vehicle has two onboard passengers destined outside of the current zone, 
it always travels straight to leave the current zone after picking up the second onboard passenger, follows a shortest path to enter the closer destination zone, and then takes a shortest path inside the zone toward the closer passenger destination. 
For example, if a vehicle in Figure \ref{fig:study region 1} delivers a single seeker from zone 5 to zone 13 via zone-level path 5-6-9-13, it first travels straight eastbound to zone 6, makes one turn northbound at a random location inside zone 6 (toward zone 9), travels straight through zone 9 to reach zone 13, and finally follows a shortest path to the seeker's destination.
If, instead, this vehicle currently carries two passengers from zone 5 to zones 13 and 14, it could travel straight eastbound to leave zone 5, turn north at the boundary between zones 5 and 6, go straight northbound all the way until entering zone 13 via the southwest corner. 

\begin{figure}[t]
    \centering
    \includegraphics[width=\textwidth]{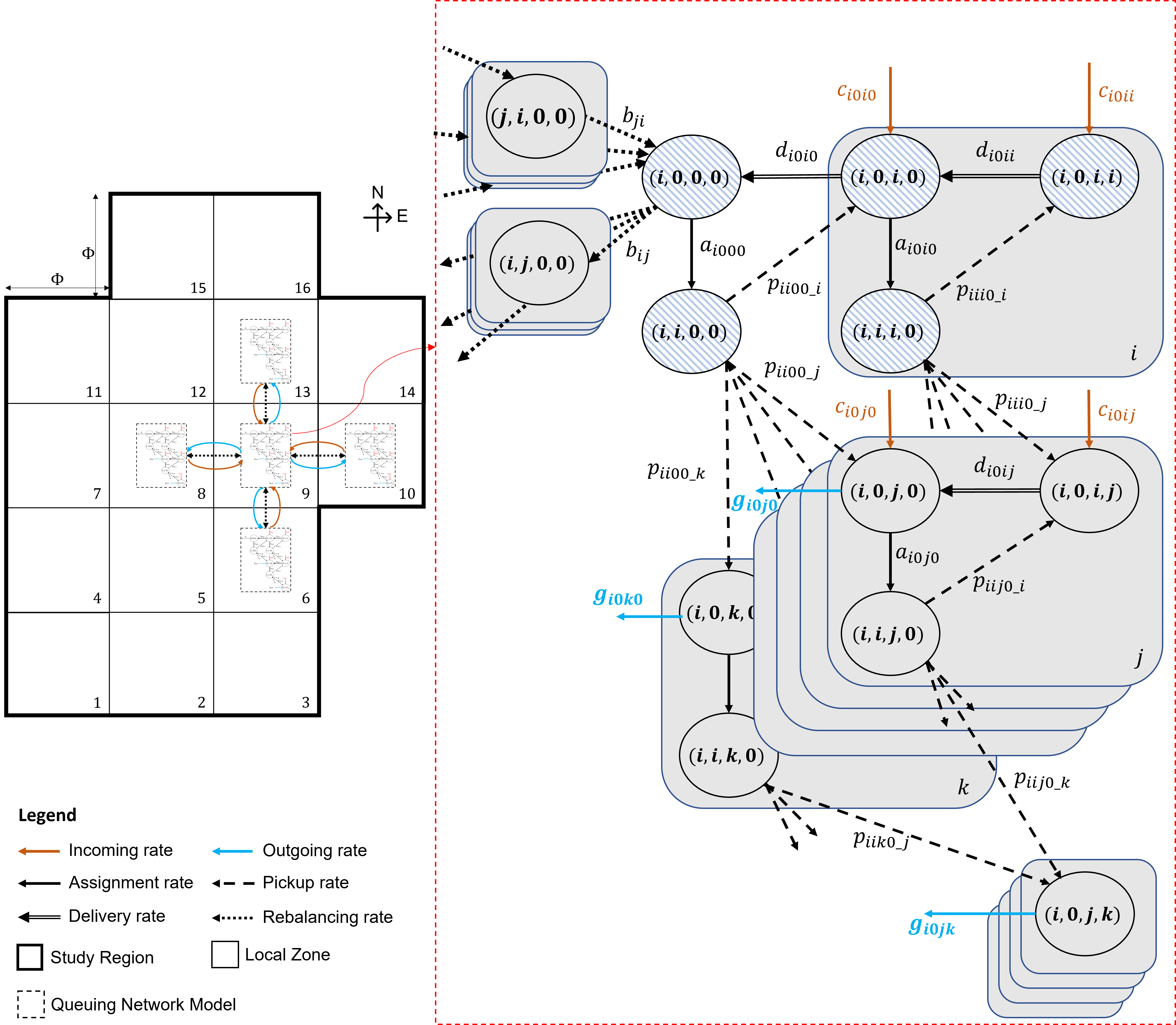}
    \caption[Caption for LOF]{State transition network of a multi-zonal ride-pooling service.}
    \label{fig:queuing model}
\end{figure} 

\subsection{Multi-zone Queuing Network Model}
We now extend the aspatial queuing network model in \citet{Daganzo2019Paper} into a multi-zone ride-pooling version that can describe ride matching, passenger assignment, and vehicle state evolution within and across multiple zones. 

\subsubsection{State and Transition}
The basic idea of the model is illustrated in Figure \ref{fig:queuing model}. We define the state of a vehicle by its current location (by zone) and its associated service task(s). Among the vehicles that are ``active'' in providing passenger service, some may be associated with no passenger at all (i.e., idling), some with one passenger (either assigned or onboard), or some with two passengers (both onboard, or one assigned and one onboard). The assigned passenger must have its origin inside the vehicle's current zone but its destination could be in any zone; similarly, each onboard passenger could have its destination in any zone. All possible states for an active vehicle in the system can be described by a vector $s = (s_0, s_1, s_2, s_3)$, where $s_0$ is the index of the vehicle's current zone; $s_1$ is the index of the origin zone of an assigned passenger, which shall take the value of $s_0$ if such a passenger exists (or 0 otherwise); and $s_2$ and $s_3$ respectively are the indices of the destination zones of up to two onboard passengers, if any, sorted in ascending order by the zone-level distance from the vehicle's current zone. If the vehicle has no assigned passenger, or fewer than two onboard passengers, the corresponding vector element takes the index of the dummy zone 0. Also, when the vehicle has only one onboard passenger, this passenger's destination zone is always indicated by index $s_2$ (and hence $s_3=0$). Therefore, for active vehicles in zone $s_0 = i$, we can know that $s_1 \in \{0, i\}$, $s_2 \text{ and } s_3 \in \mathcal{K}\cup \{0\}$, and $s_2$ and $s_3$ should satisfy the matching criteria from \ref{case:1} to \ref{case:4}. 

In the meantime, some vehicles cannot serve any passengers as they are being rebalanced; thus, their states can be described by their origin and destination zones only, i.e., $s = (s_0, s_1, 0, 0)$, with $s_0$ and $s_1$ now being the indices of the vehicle's origin and destination zones, respectively. 

The right-hand side of Figure \ref{fig:queuing model} illustrates all possible vehicle states in zone $i$, each of which is represented by a node. Vectors $(i, j, 0, 0)$ and $(j, i, 0, 0)$, $j\in\mathcal{K}\setminus\{i\}$ represent vehicles that are being rebalanced. All other vectors are for active vehicle states. For example, vector $(i, 0, 0, 0)$ represents the idle vehicle state; $(i, i, 0, 0)$ represents the vehicle state with one assigned caller and no onboard passenger; vectors $(i, 0, i, 0)$ and $(i, 0, j, 0)$, $j\in\mathcal{K}\setminus\{i\}$ represent the vehicle states with one onboard seeker whose destination is in zone $i$ and zone $j$, respectively; vectors $(i, i, i, 0)$ and $(i, i, j, 0)$, $j\in\mathcal{K}\setminus\{i\}$ represent the vehicle states with a matched caller (who is yet to board the vehicle) and a matched seeker destined in zone $i$ and zone $j$, respectively; and vectors $(i, 0, i, i)$, $(i, 0, i, j)$, and $(i, 0, j, k)$, $j\in\mathcal{K}\setminus\{i\},\ k\in\tilde{\Omega}_{ij}$ represent the vehicle states with two onboard passengers whose destinations are both in zone $i$, one in zone $i$ and one in zone $j$, and one in zone $j$ and one in zone $k$, respectively. 

A vehicle may change its state when it receives a new assignment, picks up an assigned caller, delivers an onboard passenger, crosses a zone boundary, or starts or ends a rebalancing trip. Each state change is represented by an arrow in the figure. Let $a_s$ [veh/hr] denote the rate to assign callers to the vehicles in state $s$, illustrated by the black solid arrows in Figure \ref{fig:queuing model}; $p_{s\_i}$ [veh/hr] the rate for vehicles in state $s$ to pick up an assigned caller with the destination in zone $i\in\mathcal{K}$, illustrated by the dashed arrows; $d_s$ [veh/hr] the rate for vehicles in state $s$ to deliver an onboard passenger, illustrated by the double arrows; $c_s$ (illustrated by the orange solid arrows) and $g_s$ (illustrated by the blue solid arrows) the rates for vehicles in state $s$ to enter and exit a zone, respectively; and $b_{ij}$ [veh/hr] the rate for idle vehicles to be rebalanced from zone $i$ to zone $j\in\mathcal{K}\setminus\{i\}$, as illustrated by the dotted arrows. All aforementioned rates should be non-negative.

Figure \ref{fig:queuing model} illustrates how a vehicle may change its states. An idle vehicle in state $(i, 0, 0, 0)$ changes to state $(i, i, 0, 0)$ once it is assigned to a caller; such a transition occurs at rate $a_{i000}$. Then, the vehicle travels toward the origin of the assigned caller to pick up this caller (who could be intra-zonal or inter-zonal). 

If the caller is intra-zonal, the state-$(i, i, 0, 0)$ vehicle, upon picking up the intra-zonal caller, transitions to state $(i, 0, i, 0)$ and the caller becomes a seeker, which occurs at rate $p_{ii00\_i}$. Then, this vehicle moves on to deliver the seeker toward its destination. During the delivery trip, this state-$(i, 0, i, 0)$ vehicle may receive a new assignment of a matched caller, occurring at rate $a_{i0i0}$, and become a state-$(i, i, i, 0)$ vehicle. If this caller is also intra-zonal, the vehicle travels to pick up the assigned caller and changes to state $(i, 0, i, i)$, which occurs at rate $p_{iii0\_i}$; then, as the vehicle occupancy reaches the capacity of 2, the vehicle has to deliver the passenger with a closer destination and changes to state $(i, 0, i, 0)$ upon delivery, occurring at rate $d_{i0ii}$. Otherwise, if a state-$(i, 0, i, 0)$ vehicle receives no new assignment during the delivery trip, it drops off the intra-zonal seeker at its destination, occurring at rate $d_{i0i0}$, after which the vehicle becomes an idle one in state $(i, 0, 0, 0)$. These transitions are highlighted by the top part of the network (with shaded circles) in Figure \ref{fig:queuing model}, which is very similar to the queuing network of single-zone ride-matching service in \citet{Daganzo2020}.

The caller picked up by a vehicle in state $(i, i, 0, 0)$ or $(i, i, i, 0)$ could also be inter-zonal with its destination in zone $j\in\mathcal{K}\setminus\{i\}$. In such a case, a state-$(i, i, 0, 0)$ vehicle would change to state $(i, 0, j, 0)$ upon pickup, which occurs at rate $p_{ii00\_j}$. This vehicle, while en-route toward a zone boundary to deliver the passenger, may (i) further receive a new assignment of an intra-zonal caller, and transition to states $(i, i, j, 0)$ and $(i, 0, i, j)$ and then state $(i, 0, j, 0)$ -- this transition path is similar to that of vehicle state $(i, 0, i, 0)$ (as contained in a gray layer with subscript $i$), and may occur for all $j \in\mathcal{K}\setminus\{i\}$ (as represented by the multiple grey layers in Figure \ref{fig:queuing model}); (ii) receive a new assignment of an inter-zonal caller whose destination is in zone $k$, and transition to state $(i, 0, j, k)$; and (iii) receive no further assignment, and exit the current zone $i$, occurring at rate $g_{i0j0}$ -- in this case, the vehicle transitions to state $(i', 0, j, 0)$ where $i'\in\{A^{\text{E}}_i, A^{\text{N}}_i, A^{\text{W}}_i, A^{\text{S}}_i\}$ is one of zone $i$'s adjacent zones. 

Recall that when there are two onboard passengers, we always use $s_2$ to denote the closer destination zone. Therefore, if an inter-zonal caller with destination in zone $j\in\mathcal{K}\setminus\{i\}$ is picked up by a vehicle in state $(i, i, i, 0)$, then the vehicle will also transition to state $(i, 0, i, j)$ upon pickup. Similarly, by symmetry, a vehicle may enter state $(i, 0, j, k)$ in two ways: a state-$(i, i, j, 0)$ vehicle picking up a caller with destination in zone $k$ (at rate $p_{iij0\_k}$), or a vehicle in state $(i, i, k, 0)$ picking up a caller with destination in zone $j$ (at rate $p_{iiko\_j}$). This vehicle will no longer accept new assignments but leave zone $i$ to deliver two onboard inter-zonal passengers, at rate $g_{i0jk}$. 
For simplicity, we consider that state-$(i, 0, j, k)$ vehicles will remain in this state even after leaving zone $i$, until reaching zone $j$; after that, the vehicles transition to state $(j, 0, j, k)$.

Finally, vehicles from adjacent zones may enter or pass zone $i$ (while delivering passengers) at certain rates; e.g., $c_{i0i0}$ and $c_{i0j0}$ for vehicles carrying a seeker destined in zone $i$ and zone $j$, respectively; and $c_{i0ii}$ and $c_{i0ij}$ for vehicles with two onboard passengers that both go to zone $i$, and go to zone $i$ and zone $j$, respectively. In addition, as mentioned before, idle vehicles in state $(i, 0, 0, 0)$ may be rebalanced to another zone $j\in\mathcal{K}\setminus\{i\}$ at rate $b_{ij}$, after which they transition to state $(i, j, 0, 0)$. Meanwhile, idle vehicles from zone $j\in\mathcal{K}\setminus\{i\}$ can be rebalanced to zone $i$ at rate $b_{ji}$, i.e., transitioning from $(j, i, 0, 0)$ state to $(i, 0, 0, 0)$ state. 

The above-mentioned state transitions can be illustrated with a simple example in Figure \ref{fig:study region}. Suppose an idle vehicle in zone $5$, currently in state $(5, 0, 0, 0)$, is assigned to a caller going from zone 5 to zone 13, and its delivery trip follows the zone path 5-6-9-13. This vehicle will immediately transition to state $(5, 5, 0, 0)$ and deadhead to pick up this caller. After pickup, it stays in state $(5, 0, 13, 0)$ while still in zone $5$. If it is matched with a new caller (e.g., going from zone 5 to zone $14$) before exiting zone $5$, it transitions to state $(5, 5, 13, 0)$ immediately after the assignment, 
state $(5, 0, 13, 14)$ upon picking up the matched caller (and stays in this state while heading to zone $13$), state $(13, 0, 13, 14)$ upon entering zones $13$, state $(13, 0, 14, 0)$ upon dropping off the first passenger in zone $13$, state $(14, 0, 14, 0)$ upon leaving zone $13$ to enter zone $14$ (if it is not matched with other callers), and state $(14, 0, 0, 0)$ upon delivering the second passenger. Otherwise, if no new assignment ever occurs along its trip, the state-$(5, 0, 13, 0)$ vehicle changes to states $(6, 0, 13, 0)$, $(9, 0, 13, 0)$, and $(13, 0, 13, 0)$ upon entering zones $6$, $9$, and $13$ sequentially. 

Now, in the queuing network model, we have two sets of variables, which are the transition flow rates, and the number of vehicles in each state $s$, denoted by $\{n_s\}$. We can build a system of equations to describe the connections among these variables in the long-term steady state. First, the vehicle state transition flow in Figure \ref{fig:queuing model} should be conserved at all state nodes, such that the number of vehicles in any state $\{n_s\}$ shall remain steady. This is discussed next in Section \ref{sec:flow conservation}. Then, we relate the flow rates and $\{n_s\}$, based on the matching criteria, assignment strategy, demand generation process, and the Little's formula; this will be discussed in Section \ref{sec:assignment and duration}.


\subsubsection{Flow Conservation within and across Zones}\label{sec:flow conservation}
The flow in the queuing network, as shown in Figure \ref{fig:queuing model}, should be conserved within each zone as well as across zones. Within a specific zone, e.g., $i\in\mathcal{K}$, the summations of the in- and out-flow at each state node must be equal. The reader can verify that the list of such conservation equations are given in Eqn. \eqref{eq:flow conserve within zone} below. 

\begin{subequations}\label{eq:flow conserve within zone}
\renewcommand{\theequation}{\theparentequation.\arabic{equation}}
\begin{align}
0 =&\ \sum_{j\in\mathcal{K}\setminus\{i\}}(b_{ij} - b_{ji}) + a_{i000} - d_{i0i0},\quad && \forall i\in\mathcal{K} \label{eq:conserve i000}\\
0 =&\ \sum_{j\in \mathcal{K}}p_{ii00\_j} - a_{i000},\quad && \forall i\in\mathcal{K} \label{eq:conserve ii00}\\
0 =&\ a_{i0i0} + d_{i0i0} - c_{i0i0} - d_{i0ii} - p_{ii00\_i},\quad && \forall i\in\mathcal{K} \label{eq:conserve i0i0}\\
0 =&\ \sum_{j\in \mathcal{K}}p_{iii0\_j} - a_{i0i0},\quad && \forall i\in\mathcal{K} \label{eq:conserve iii0}\\
0 =&\ d_{i0ii} - p_{iii0\_i} - c_{i0ii}, \quad && \forall i\in\mathcal{K}\label{eq:conserve i0ii}\\
0 =&\ g_{i0j0} + a_{i0j0} - p_{ii00\_j} - c_{i0j0} - d_{i0ij}, \quad && \forall i\in\mathcal{K},\ j\in \mathcal{K}\setminus\{i\}\label{eq:conserve i0j0}\\
0 =&\ \sum_{k\in \Omega_{ij}\cup\{i\}}p_{iij0\_k} - a_{i0j0},\quad &&\forall i\in\mathcal{K},\ j\in \mathcal{K}\setminus\{i\}\label{eq:conserve iij0}\\
0 =&\ d_{i0ij} - c_{i0ij} - p_{iii0\_j} - p_{iij0\_i},\quad &&\forall i\in\mathcal{K},\ j\in \mathcal{K}\setminus\{i\}\label{eq:conserve i0ij}\\
0 =&\ g_{i0jk} - p_{iij0\_k} - \mathbbm{I}(k \neq j)\cdot p_{iik0\_j},\quad &&\forall i\in\mathcal{K},\ j\in \mathcal{K}\setminus\{i\},\ k\in\tilde{\Omega}_{ij}\label{eq:conserve i0jk}
\end{align}
\end{subequations}
where $\mathbbm{I}(k \neq j)=1$ if $k \neq j$, or 0 otherwise; it is used to avoid counting the pickup rate $p_{iij0\_j}$ twice when the destination zones of two onboard inter-zonal passengers are the same.

Vehicle flow must also be conserved across zones. Note that the outgoing flow of vehicles (with an inter-zonal seeker) from a zone, $g$, becomes incoming flow of one of its adjacent zones, $c$. Since vehicles with an inter-zonal seeker may be matched with new callers, we must track their movements across zones, and construct a relationship between the outgoing and incoming flow across zone boundaries. Despite all feasible travel paths a vehicle with an inter-zonal seeker could follow to reach its seeker's destinations, this vehicle can only move into at most two adjacent zones without creating zone-level detours; this is specified in Eqn. \eqref{eq:feasible directions}. The platform shall choose zone-level paths for the vehicles with an inter-zonal seeker (i.e., currently in state $(i, 0, j, 0)$), to enhance system efficiency; e.g., by impacting the probability for the vehicle to find a matched caller en-route. These path decisions collectively specify the proportion of state-$(i, 0, j, 0)$ vehicles that enter a zone $i'\in \{A^\text{E}_{i}, A^\text{N}_{i}, A^\text{W}_{i}, A^\text{S}_{i}\}$ upon leaving zone $i$, which is denoted by $\delta_{ij\_i'}\in[0, 1]$. 
The vehicle path decisions should satisfy the following constraints:

\begin{align}\label{eq:constraint delta}
& \sum_{i'\in \{A^\text{E}_i, A^\text{N}_i, A^\text{W}_i, A^\text{S}_i\}} \delta_{ij\_i'} = 1, \quad \sum_{i'\in \mathcal{V}_{ij}}\delta_{ij\_i'} = 1, \quad \forall i\in \mathcal{K},\ j\in \mathcal{K}\setminus\{i\}
\end{align}
which postulate that a vehicle with an inter-zonal seeker can only travel in a feasible direction horizontally or vertically into an adjacent zone. Then, it is easy to see that the flow of such vehicles between zones must satisfy the following:

\begin{equation}\label{eq:flow with a seeker conserve across zone}
\begin{aligned}
c_{i0i0} =& \sum_{i'\in \{A^\text{E}_i, A^\text{N}_i, A^\text{W}_i, A^\text{S}_i\}}g_{i'0i0}\cdot \delta_{i'i\_i},\quad &&\forall i\in \mathcal{K}\\
c_{i0j0} =& \sum_{i'\in \{A^\text{E}_i, A^\text{N}_i, A^\text{W}_i, A^\text{S}_i\}}g_{i'0j0}\cdot \delta_{i'j\_i},\quad &&\forall i\in \mathcal{K},\ \ j\in\mathcal{K}\setminus\{i\}
\end{aligned}
\end{equation}

For vehicles with two onboard inter-zonal passengers, the flow going out of zone $i$, $g_{i0jk}$, becomes incoming flow into the closer passenger destination zone, $c_{j0jk}$; then, the cross-zone flow conservation for such vehicles are as follows,
\begin{equation}\label{eq:flow with two pax conserve across zone}
\begin{aligned}
c_{i0ii} = \sum_{i'\in \mathcal{K}\setminus\{i\}}g_{i'0ii},\quad \forall i\in \mathcal{K};\quad\quad c_{i0ij} = \sum_{i'\in \tilde{\Omega}_{ji}\setminus\{i\}}g_{i'0ij},\quad \forall i\in \mathcal{K},\ \ j\in\mathcal{K}\setminus\{i\}
\end{aligned}
\end{equation}

\subsubsection{Flow Rates, Matching, and Expected Vehicle State Duration}\label{sec:assignment and duration}

Recall that a new caller in zone $i$ will be assigned to a nearest suitable vehicle in states $(i, 0, 0, 0)$, $(i, 0, i, 0)$, and $(i, 0, j, 0),\ \forall j\in\mathcal{K}\setminus\{i\}$. The caller could be inter-zonal (traveling from zone $i$ to any zone $j$), or intra-zonal (traveling in any direction $r \in \{\text{NE}, \text{NW}, \text{SW}, \text{SE}\}$ within zone $i$). \ref{sec:appendix matching} shows that the number of vehicles suitable to serve an intra-zonal caller traveling in direction $r$ within zone $i$, denoted by $N_{ii}^{r}$, is given by 

\begin{equation}\label{eq:suitable veh}
\begin{aligned}
    &N_{ii}^r = n_{i000} + \frac{2n_{i0i0}}{9} + \sum_{j\in \Omega^r_{ii}}n_{i0j0},\quad &&\forall i \in \mathcal{K}, \ r\in\{\text{NE}, \text{NW}, \text{SW}, \text{SE}\}\\
\end{aligned}
\end{equation}
and the number of vehicles suitable to serve an inter-zonal caller going from zone $i$ to zone $j$, denoted by $N_{ij}$, is given by

\begin{equation}\label{eq:suitable veh_ij}
\begin{aligned}
    &N_{ij} = n_{i000} + \frac{n_{i0i0}}{4} + \sum_{k\in \Omega_{ij}}n_{i0k0},\quad &&\forall i \in \mathcal{K}, \ j \in G_i^{\times}\\
    &N_{ij} = n_{i000} + \frac{n_{i0i0}}{2} + \sum_{k\in \Omega_{ij}}n_{i0k0},\quad &&\forall i \in \mathcal{K}, \ j \in G_i^{+}\\
\end{aligned}
\end{equation}

Since vehicles in each state are approximately uniformly scattered in each local zone, under the nearest assignment strategy, the fraction of callers assigned to vehicles in any suitable state is proportional to the number of vehicles in that state. \ref{sec:appendix matching} further shows that $\{N_{ij}\}$ and $\{N_{ii}^r\}$ are related to the caller pickup rates as follows:

\begin{subequations}\label{eq:pickup}
\renewcommand{\theequation}{\theparentequation.\arabic{equation}}
\begin{align}
    & p_{ii00\_i} = \sum_{\substack{r\in \{\text{NE}, \text{NW}, \text{SW}, \text{SE}\}}}\left(\frac{\lambda_{ii}}{4}\cdot\frac{n_{i000}}{N^r_{ii}}\right),\  &&\forall i\in \mathcal{K}\label{eq:pickup 1}\\
    & p_{iii0\_i} = \sum_{\substack{r\in \{\text{NE}, \text{NW}, \text{SW}, \text{SE}\}}}\left(\frac{\lambda_{ii}}{4}\cdot\frac{2n_{i0i0}}{9N^r_{ii}}\right),\ &&\forall i\in \mathcal{K}\label{eq:pickup 2}\\
    & p_{iij0\_i} = \frac{\lambda_{ii}}{4}\cdot\frac{n_{i0j0}}{N^r_{ii}},\ && \forall i\in \mathcal{K}, \  j \in G^r_i,\ r \in \{\text{NE}, \text{NW}, \text{SW}, \text{SE}\} \label{eq:pickup 3}\\
    & p_{iij0\_i} = \sum_{r\in U^{r'}}\left(\frac{\lambda_{ii}}{4}\cdot\frac{n_{i0j0}}{N^{r}_{ii}}\right),\ && \forall i\in \mathcal{K}, \  j \in G^{r'}_i, \ r' \in \{\text{E}, \text{N}, \text{S}, \text{W}\}\label{eq:pickup 4}\\
    & p_{ii00\_j} = \lambda_{ij}\cdot \frac{n_{i000}}{N_{ij}},\quad &&\forall i \in \mathcal{K},\ j \in \mathcal{K}\setminus\{i\}\label{eq:pickup 5}\\
    & p_{iii0\_j} = \lambda_{ij}\cdot \frac{n_{i0i0}}{4N_{ij}},\quad &&\forall i \in \mathcal{K},\ j \in G^{\times}_i\label{eq:pickup 6}\\
    & p_{iii0\_j} = \lambda_{ij}\cdot \frac{n_{i0i0}}{2N_{ij}},\quad &&\forall i \in \mathcal{K},\ j \in G^{+}_i\label{eq:pickup 7}\\
    & p_{iik0\_j} = \lambda_{ij}\cdot \frac{n_{i0k0}}{N_{ij}}.\quad &&\forall i \in \mathcal{K},\ j \in \mathcal{K}\setminus\{i\}, \ k \in \Omega_{ij}\label{eq:pickup 8}
\end{align}
\end{subequations}

For a vehicle in zone $i$ with a seeker (i.e., in state $(i, 0, i, 0)$ or $(i, 0, j, 0)$), it may receive a new assignment during its delivery trip (which occurs at rate $a_{i0i0}$ or $a_{i0j0}$, respectively); otherwise, the vehicle will finish the delivery trip in zone $i$ (which occurs at rate $d_{i0i0}$ or $g_{i0j0}$, respectively). The probability for the latter case to occur depends on the new passenger assignment rates and the length of the delivery trip inside zone $i$. Since new passengers arrive according to Poisson processes, \ref{sec:appendix duration} shows that the rates $d_{i0i0}$ and $g_{i0j0}$ satisfy the following,

\begin{subequations}\label{eq:delivery and outgoing}
\renewcommand{\theequation}{\theparentequation.\arabic{equation}}
\begin{align}
    d_{i0i0} = &c_{i0i0}\cdot e^{-\frac{a_{i0i0}}{n_{i0i0}}\cdot\frac{5\Phi}{6v}} + (p_{ii00\_i} + c_{i0ii})\cdot e^{-\frac{a_{i0i0}}{n_{i0i0}}\cdot\frac{2\Phi}{3v}} + p_{iii0\_i} \cdot e^{-\frac{a_{i0i0}}{n_{i0i0}}\cdot\frac{\Phi}{2v}},\quad &&\forall i\in\mathcal{K} \label{eq:delivery only}\\
    g_{i0j0} = &c_{i0j0}\cdot e^{-\frac{a_{i0j0}}{n_{i0j0}}\cdot\frac{\Phi}{v}} + (p_{ii00\_j} + c_{i0ij})\cdot e^{-\frac{a_{i0j0}}{n_{i0j0}}\cdot\frac{\Phi}{2v}} + (p_{iii0\_j} + p_{iij0\_i})\cdot e^{-\frac{a_{i0j0}}{n_{i0j0}}\cdot\frac{\Phi}{3v}}.\quad && \forall i\in\mathcal{K},\ j\in\mathcal{K}\setminus\{i\}  \label{eq:outgoing only}
\end{align}
\end{subequations}

Based on Eqn. \eqref{eq:flow conserve within zone} and Eqn. \eqref{eq:flow with a seeker conserve across zone}-\eqref{eq:delivery and outgoing}, we can express the transition flow rates $\{a_s\}$, $\{p_{s\_i}\}$, $\{d_s\}$, $\{g_s\}$, and $\{c_s\}$ as functions of $\{n_{i000}\}$, $\{n_{i0i0}\}$, $\{n_{i0j0}\}$, $\boldsymbol \delta$, and $\{b_{ij}\}$. We can further use these flow rates to quantify the number of vehicles in any of the remaining states, by multiplying the expected duration of that vehicle state and the vehicle arrival rate into that state, per the Little's Law \citep{little1961proof}. \ref{sec:appendix duration} gives detailed derivations on the expected duration in each vehicle state, which leads to the following relationships:

\begin{subequations}\label{eq:fleet}
\renewcommand{\theequation}{\theparentequation.\arabic{equation}}
\small
\begin{align}
    & n_{ii00} = \frac{0.63\Phi n_{i000}}{v}\left(\sum_{\substack{r\in\{\text{NE}, \text{NW},\\ \ \ \ \text{SW}, \text{SE}\}}}\frac{\lambda_{ii}}{4(N_{ii}^r)^{\frac{3}{2}}} + \sum_{j\in\mathcal{K}\setminus\{i\}}\frac{\lambda_{ij}}{(N_{ij})^{\frac{3}{2}}}\right),\quad\quad\quad\quad\quad\quad\quad\quad\quad\quad\quad\quad\quad\ \ \forall i\in\mathcal{K}\label{eq:fleet 1}\\
    & n_{iii0} = \frac{0.63\Phi n_{i0i0}}{2v}\left(\sum_{\substack{r\in\{\text{NE}, \text{NW},\\ \ \ \ \text{SW}, \text{SE}\}}}\frac{\lambda_{ii}}{9(N_{ii}^r)^{\frac{3}{2}}} + \sum_{j\in G^{\times}_i}\frac{\lambda_{ij}}{2(N_{ij})^{\frac{3}{2}}} + \sum_{j\in G^+_i}\frac{\lambda_{ij}}{(N_{ij})^{\frac{3}{2}}}\right),\quad\quad\quad\quad\quad\quad\quad\ \  \forall i\in\mathcal{K}\label{eq:fleet 2}\\
    & n_{iij0} = \frac{0.63\Phi n_{i0j0}}{v}\left(\frac{\lambda_{ii}}{4(N_{ii}^r)^{\frac{3}{2}}} + \sum_{k\in\Omega_{ij}}\frac{\lambda_{ik}}{(N_{ik})^{\frac{3}{2}}}\right),\quad\quad\quad\ \ \ \ \ \forall i\in \mathcal{K}, \  j \in G^r_i,\ r \in \{\text{NE}, \text{NW}, \text{SW}, \text{SE}\}\label{eq:fleet 3}\\
    & n_{iij0} = \frac{0.63\Phi n_{i0j0}}{v}\left(\sum_{r\in U^{r'}}\frac{\lambda_{ii}}{4(N_{ii}^{r})^{\frac{3}{2}}} + \sum_{k\in\Omega_{ij}}\frac{\lambda_{ik}}{(N_{ik})^{\frac{3}{2}}}\right),\quad\quad\quad\quad\ \ \forall i\in\mathcal{K},\ j\in G^{r'}_i,\ r'\in\{\text{E}, \text{N}, \text{W}, \text{S}\}\label{eq:fleet 4}\\
    & n_{i0ii} = \sum_{i'\in G^{+}_i}g_{i'0ii}\cdot\frac{5\Phi}{8v} + \sum_{i'\in G^{\times}_i}g_{i'0ii}\cdot\frac{23\Phi}{30v} + p_{iii0\_i}\cdot\frac{\Phi}{2v},\quad\quad\quad\quad\quad\quad\quad\quad\quad\quad\quad\quad\quad\quad\ \ \ \ \ \forall i\in\mathcal{K} \label{eq:fleet 5}\\
    & n_{i0ij} =\sum_{i'\in G^{+}_i\cap\tilde{\Omega}_{ji}}g_{i'0ij}\cdot\frac{5\Phi}{6v} + \sum_{i'\in G^{\times}_i\cap\tilde{\Omega}_{ji}}g_{i'0ij}\cdot\frac{\Phi}{v} +  p_{iij0\_i}\cdot\frac{2\Phi}{3v} + p_{iii0\_j}\cdot\tilde{T}_{i0ij},\quad \forall i\in\mathcal{K},\ j \in \mathcal{K}\setminus\{i\}\label{eq:fleet 6}\\
    &n_{i0jk} = p_{iij0\_k}\cdot \left(\tilde{T}_{i0jk} + \frac{L_{ij} - \Phi}{v}\right) + \mathbbm{I}(k\neq j)\cdot p_{iik0\_j}\cdot \left(\frac{L_{ij}-\Phi/2}{v}\right),\quad\ \ \forall i\in\mathcal{K},\ j \in G^{+}_i,\ k\in\tilde{\Omega}_{ij}\label{eq:fleet 7.1}\\
    &n_{i0jk} = p_{iij0\_k}\cdot \left(\tilde{T}_{i0jk} + \frac{L_{ij} - 3\Phi/2}{v}\right) + \mathbbm{I}(k\neq j)\cdot p_{iik0\_j}\cdot \left(\frac{L_{ij}-\Phi}{v}\right),\quad \forall i\in\mathcal{K},\ j \in G^{\times}_i,\ k\in\tilde{\Omega}_{ij}\label{eq:fleet 7.2}\\
    & n_{ij00} = b_{ij}\cdot\frac{L_{ij}}{v},\quad\quad\quad\quad\quad\quad\quad\quad\quad\quad\quad\quad\quad\quad\quad\quad\quad\quad\quad\quad\quad\quad\quad\quad\quad\quad\quad\ \ \ \ \ \forall i\in\mathcal{K},\ j\in G^{\times}_i\label{eq:fleet 8}\\
    & n_{ij00} = b_{ij}\cdot\frac{L_{ij} + \Phi/3}{v},\quad\quad\quad\quad\quad\quad\quad\quad\quad\quad\quad\quad\quad\quad\quad\quad\quad\quad\quad\quad\quad\quad\quad\quad \ \ \ \ \forall i\in\mathcal{K},\ j\in G^{+}_i\label{eq:fleet 9}
    \end{align}
\end{subequations}
where, for all $i \in \mathcal{K}$, $j \in \mathcal{K}\setminus\{i\}$ and $k\in\tilde{\Omega}_{ij}$,
\begin{equation*}
\begin{aligned}
\footnotesize
    &\tilde{T}_{i0ij} = \frac{\Phi/v}{c_{i0i0} + p_{ii00\_i} + d_{i0ii}}\left[
    \frac{5 c_{i0i0}/6}{1 - e^{-\frac{a_{i0i0}}{n_{i0i0}}\frac{5\Phi}{6v}}} + \frac{2 (c_{i0ii} + p_{ii00\_i})/3}{1 - e^{-\frac{a_{i0i0}}{n_{i0i0}}\frac{2\Phi}{3v}}} + \frac{p_{iii0\_i}/2}{1 - e^{-\frac{a_{i0i0}}{n_{i0i0}}\frac{\Phi}{2v}}}\right] - \frac{n_{i0i0}}{a_{i0i0}},\\
    &\tilde{T}_{i0jk} = \frac{\Phi/v}{c_{i0j0} + p_{ii00\_j} + d_{i0ij}}\left[
    \frac{c_{i0j0}}{1 - e^{-\frac{a_{i0j0}}{n_{i0j0}}\frac{\Phi}{v}}} + \frac{(c_{i0ij} + p_{ii00\_j})/2}{1 - e^{-\frac{a_{i0j0}}{n_{i0j0}}\frac{\Phi}{2v}}} + \frac{(p_{iii0\_j} + p_{iij0\_i})/3}{1 - e^{-\frac{a_{i0j0}}{n_{i0j0}}\frac{\Phi}{3v}}}\right] - \frac{n_{i0j0}}{a_{i0j0}}.
\end{aligned}
\end{equation*}

\subsection{Performance Evaluation}\label{sec:performance}
With the number of vehicles in each state $\{n_s\}$, we can evaluate the system performance, including the agency cost and passenger travel experience. There are many ways to formulate the system costs; in this paper, as an illustration, we assume the agency cost is proportional to the total vehicle hours in operation, and the passenger cost is proportional to the summation of their total travel times (from calling for service to arriving at the destination). 

The total number of active vehicles in zone $i$, denoted by $M_i$ [veh], is
\begin{equation*}\label{eq:one zone fleet}
    M_i = n_{i000} + n_{ii00} + n_{i0i0} + n_{iii0} + n_{i0ii} +  \sum_{j\in\mathcal{K}\setminus\{i\}}\left(n_{i0j0} + n_{i0ij} + n_{iij0} + \sum\limits_{k\in\tilde{\Omega}_{ij}}n_{i0jk}\right),\quad \forall i\in\mathcal{K}
\end{equation*}
and that of rebalanced vehicles is $M_b = \sum_{i\in\mathcal{K}}\sum_{j\in\mathcal{K}\setminus\{i\}}n_{ij00}$ [veh]. Therefore, the total fleet needed across all zones can be written as 

\begin{equation}\label{eq:total fleet}
M = \sum_{i\in\mathcal{K}}M_i + M_b.
\end{equation}
This quantity is also the total vehicle hours in operation per hour, i.e., with unit veh-hr/hr. Therefore, the total agency cost per hour is simply $\gamma M$ \$/hr.

Also, we can calculate the total trip time spent by all passengers per hour, denoted by $P$ [pax-hr/hr]. Its formula is based on that for $M_i$, except that each term $n_s$ is weighted by the number of assigned and onboard passengers per vehicle in state $s$; i.e.,

\begin{equation}\label{eq:total pax}
    P = \sum_{i\in\mathcal{K}}\left[n_{ii00} + n_{i0i0} + 2(n_{iii0} + n_{i0ii}) +  \sum_{j\in\mathcal{K}\setminus\{i\}}\left(n_{i0j0} + 2n_{i0ij} + 2n_{iij0} + 2\sum\limits_{k\in\tilde{\Omega}_{ij}}n_{i0jk}\right) \right].
\end{equation}
Therefore, the total passenger cost per hour is $\beta P$ \$/hr. The system-wide cost per passenger, denoted by $Z$ 
[\$/pax], hence, is as follows,

\begin{equation} \label{obj}
    Z = \frac{\gamma M+\beta P}{\sum\limits_{i, j\in\mathcal{K}}\lambda_{ij}}.
\end{equation}

\subsection{Optimization Model and Solution Approach
}\label{sec:optimization}
Eqn. \eqref{eq:flow conserve within zone} and Eqn. \eqref{eq:flow with a seeker conserve across zone}-\eqref{eq:fleet} collectively describe how, in the steady state of the queuing network, the following variables are related: the vehicle transition flow rates $\{a_s\}$, $\{p_{s\_i}\}$, $\{d_s\}$, $\{c_s\}$, $\{g_s\}$, $\{b_{ij}\}$; the vehicle counts in all states $\{n_s\}$; and vehicle path fractions $\boldsymbol \delta$. It can be easily verified that, conditional on any choice of $\boldsymbol \delta, \{b_{ij}\}$, and $\{n_{i000}\}$ (idle vehicle counts), all other variables can be solved from the system of equations, and so can the system performance $Z$. As such, the TNC can optimize their decisions of $\boldsymbol \delta, \{b_{ij}\}$, and $\{n_{i000}\}$ to minimize the system-wide cost per passenger for providing a multi-zonal ride-pooling service; i.e., 
\begin{align}
\min_{\boldsymbol \delta, \{n_{i000}\}, \{b_{ij}\}} \quad & Z
 \nonumber\\
\textrm{s.t.} \quad & \text{Eqn. } \eqref{eq:side directions}-\eqref{obj}\nonumber\\
&a_s, p_{s\_i}, d_s, c_s, g_s,
n_s\geq 0, &&\forall s, \forall i\in\mathcal{K} \label{eq:constraint non-neg}\\
&N_{ij} > 1,\ N^r_{ii} > 1,\quad &&\forall i\in\mathcal{K},\ j\in\mathcal{K}\setminus\{i\},\ r\in\{\text{NE, NW, SW, SE}\}\label{eq:constraint N_ij}\\
& n_{i000} > 0, \quad &&\forall i\in\mathcal{K}\label{eq:constraint n_i000}\\
&\delta_{ij\_i'}\in [0, 1],\quad &&\forall i\in\mathcal{K},\ j\in\mathcal{K}\setminus\{i\},\ i'\in \{A^\text{E}_i, A^\text{N}_i, A^\text{W}_i, A^\text{S}_i\}\label{eq:constraint delta domain}\\
&b_{ij}\geq 0,\quad &&\forall i\in\mathcal{K},\ j\in\mathcal{K}\setminus\{i\}\label{eq:constraint b}
\end{align}
where Eqn. \eqref{eq:constraint non-neg} enforces all flow rates and fleet sizes to be non-negative real numbers;\footnote{All vehicle counts $\{n_s\}$ are interpreted as expectations, and hence we allow them to be non-negative real numbers.} Eqn. \eqref{eq:constraint N_ij} ensures at least one suitable vehicle on average is available to serve a trip between any zone pair and in any direction; 
and Eqn. \eqref{eq:constraint n_i000}-\eqref{eq:constraint b} 
specify the value domains for the decision variables.

This optimization problem involves only continuous variables, and hence it can be solved by gradient-based search methods. However, it is highly nonlinear and nonconvex, so it is challenging to directly search all involved variables simultaneously. 
We propose to decompose the problem by observing that the set of rebalancing flow $\{b_{ij}\}$ only appear as $\{\sum_{j\in\mathcal{K}\setminus\{i\}}(b_{ij}-b_{ji}), \forall i\}$ in the flow conservation constraint \eqref{eq:conserve i000}, and as a linear summation term in the objective function (see Eqn. \eqref{eq:fleet}-\eqref{eq:total fleet}). Given any combination of $\boldsymbol \delta$ and $\{n_{i000}\}$, we can treat $\{\sum_{j\in\mathcal{K}\setminus\{i\}}(b_{ij}-b_{ji}), \forall i\}$ as $|\mathcal{K}|$ variables, and solve them together with $\{a_s\}$, $\{p_{s\_i}\}$, $\{d_s\}$, $\{c_s\}$, $\{g_s\}$, $\{n_{i0i0}\}$, and $\{n_{i0j0}\}$, from Eqn. \eqref{eq:flow conserve within zone} and Eqn. \eqref{eq:flow with a seeker conserve across zone}-\eqref{eq:delivery and outgoing}. Then, we can solve $\{b_{ij}\}$ from a transportation problem, using $\sum_{j\in\mathcal{K}\setminus\{i\}}(b_{ij}-b_{ji})$ as the ``demand" or ``supply" of zone $i$, and $L_{ij}/v$ or $(L_{ij} + \Phi/3)/v$ from Eqn. \eqref{eq:fleet 8}-\eqref{eq:fleet 9} as the shipment costs. After obtaining all flow rates, we can quantify the vehicle counts in all states based on Eqn. \eqref{eq:fleet}. Although this paper could not provide a rigorous proof to guarantee the existence or uniqueness of the solution to the system of nonlinear equations Eqn. \eqref{eq:side directions}-\eqref{obj} (or lack thereof), the proposed decomposition approach always finds a reasonable solution in our numerical experiments.

The basic solution algorithm can be summarized as follows. 
\begin{enumerate}[
    leftmargin=2 cm,
    label={\textit{Step} \arabic*:},
    ref={Step \arabic*}]
    \item \label{step 1} Randomly initialize variables $\boldsymbol \delta$ and $\{n_{i000}\}$; set the iteration number to be 1. 
    \item \label{step 2} Given the current $\boldsymbol \delta$ and $\{n_{i000}\}$, evaluate the system performance as follows.
    \begin{enumerate}[
    leftmargin=*,
    label={\textit{Step} 2.\arabic*:},
    ref={Step 2.\arabic*}]
    \item \label{step 2a} Solve Eqn. \eqref{eq:flow conserve within zone} and Eqn. \eqref{eq:flow with a seeker conserve across zone}-\eqref{eq:delivery and outgoing} to obtain values of $\{a_s\}$, $\{p_{s\_i}\}$, $\{d_s\}$, $\{c_s\}$, $\{g_s\}$, $\{n_{i0i0}\}$, $\{n_{i0j0}\}$, and $\{\sum_{j\in\mathcal{K}\setminus\{i\}}(b_{ij}-b_{ji})\}$.
    \item \label{step 2b} Use $\{\sum_{j\in\mathcal{K}\setminus\{i\}}(b_{ij}-b_{ji})\}$ as input, and solve $\{b_{ij}\}$ from the transportation problem.
    \item \label{step 2c} Solve vehicle counts $\{n_s\}$ from Eqn. \eqref{eq:fleet}. If constraints \eqref{eq:constraint non-neg}-\eqref{eq:constraint b} are satisfied, evaluate the system performance from Eqn. \eqref{eq:total fleet}-\eqref{obj}; otherwise, set $Z=\infty$.
    \end{enumerate}
    \item \label{step perturb} Record (or update) the best solution so far, if possible. Perturb $\boldsymbol \delta$ and $\{n_{i000}\}$ by small increments, and repeat \ref{step 2} to numerically compute the gradient. Then, conduct a gradient-based search to find a new solution $\boldsymbol \delta$ and $\{n_{i000}\}$; increase the iteration number by 1 and go to \ref{step 2}.
    \item \label{last step} If a maximum iteration number is reached, or if  the best objective value does not improve in a number of iterations, terminate the algorithm. 
\end{enumerate}

Before moving on to the numerical results, we suggest two optional treatments that may facilitate the solution process. First, note that the solution procedure for the system of nonlinear equations in \ref{step 2a} can be very sensitive to the initial guess. Therefore, we can try multiple guesses, and/or keep a memory of past solutions (i.e., from different $\boldsymbol \delta$ and $\{n_{i000}\}$ values) as candidate initial guesses. If no solution can be found after a number of such initial guesses, then we consider this combination of $\boldsymbol \delta$ and $\{n_{i000}\}$ to be infeasible (and set $Z=\infty$). 
Second, note that the cross-zone vehicle flow rates $\boldsymbol g$ and $\{c_{i0j0}\}$, $\{c_{i0ij}\}$ are coupled through Eqn. \eqref{eq:flow with a seeker conserve across zone}-\eqref{eq:flow with two pax conserve across zone} and Eqn. \eqref{eq:delivery and outgoing}, and they are linearly related once assignment rates $\{a_s\}$ and pickup rates $\{p_{s\_i}\}$ are known. Hence, it is beneficial to first calculate $\{p_{s\_i}\}$ based on Eqn. \eqref{eq:pickup}, obtain $\{a_s\}$ with Eqn. \eqref{eq:flow conserve within zone}, and use matrix elimination methods to solve $\boldsymbol g$ and $\{c_{i0j0}\}$, $\{c_{i0ij}\}$ from a linear system of equation. After that, $\{c_{i, 0, i, 0}\}$, $\{c_{i, 0, i, i}\}$, $\{d_{i, 0, i, 0}\}$ can also be computed. 

\section{Numerical Analysis}\label{sec:numerical}

\begin{figure}[t]
     \centering
     \begin{subfigure}[b]{0.4\textwidth}
         \centering 
         \includegraphics[width = 0.7\textwidth]{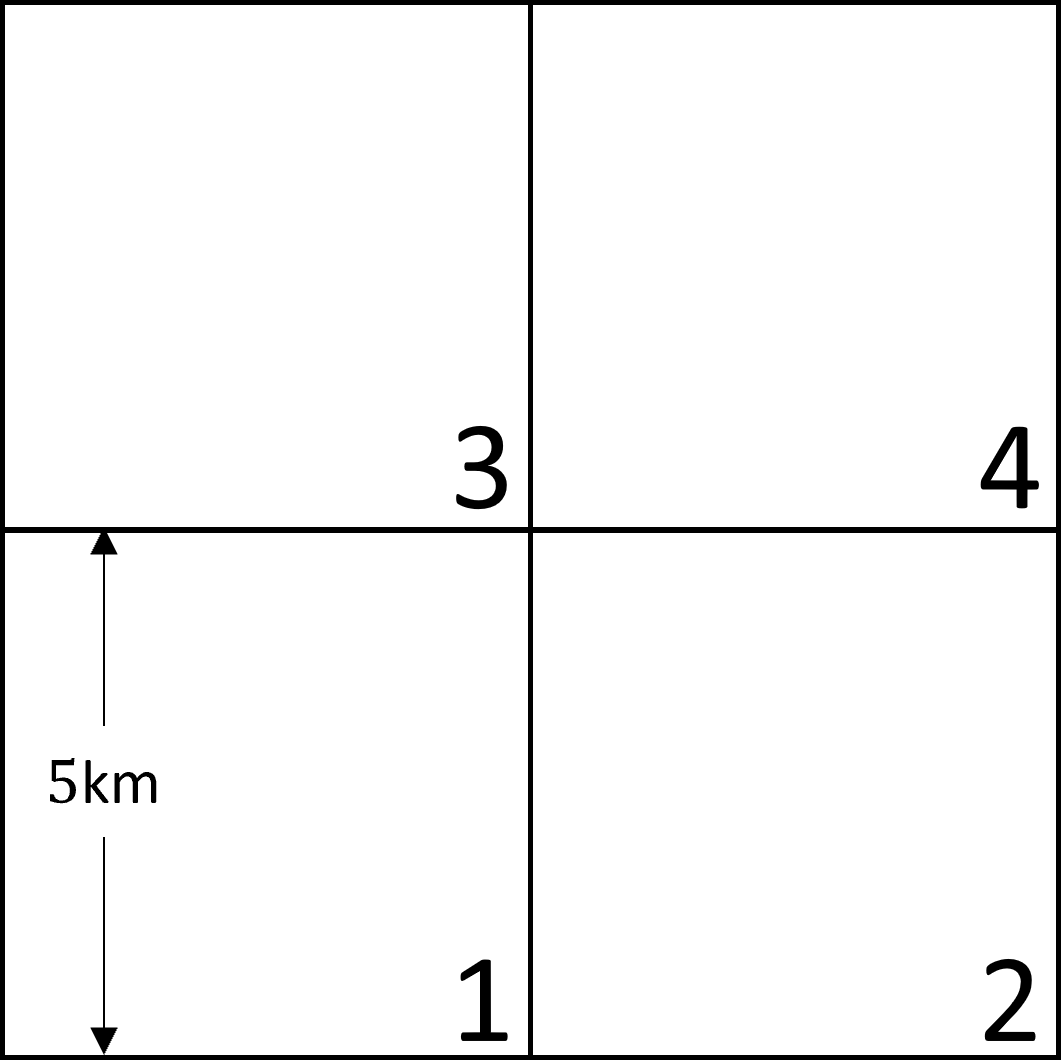} 
         \caption{Study region with 4 zones.}
         \label{fig:hypo region 4}
     \end{subfigure}
     \begin{subfigure}[b]{0.4\textwidth}
         \centering
         \includegraphics[width=0.7\textwidth]{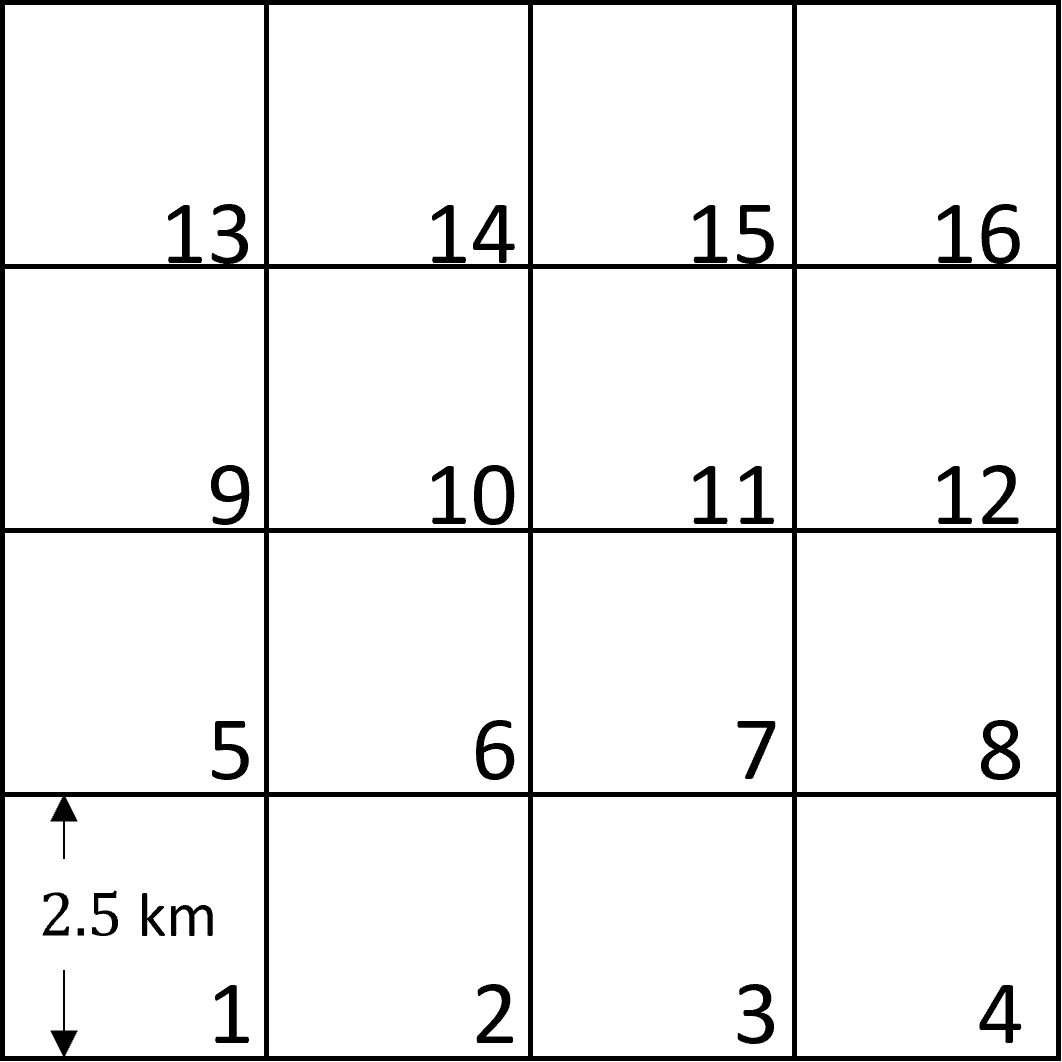}
         \caption{Study region with 16 zones.}
         \label{fig:hypo region 16}
     \end{subfigure}
    \caption{The hypothetical study region partitioned into $2\times 2$ zones and $4\times 4$ zones, respectively.}
    \label{fig:hypo region}
\end{figure}

In this section, we first consider a hypothetical study region to see how the optimal service decisions and system performance vary in different demand scenarios for different zone partitions. The model predictions are also verified with agent-based simulations. Later, we apply the proposed model to a real-world case study for downtown and suburban Chicago to draw further insights.


The ``root" function in Scipy \citep{2020SciPy-NMeth} is used to find the solution to the system of equations in \ref{step 2a}. The GNU Linear Programming Kit solver \citep{makhorin2008glpk} embedded in Pyomo \citep{hart2011pyomo,bynum2021pyomo} is used to solve the transportation problem in \ref{step 2b}. The Sequential Least Squares Programming algorithm \citep{kraft1988software} implemented in Scipy \citep{2020SciPy-NMeth} is used to conduct the gradient-based search in Step 3. 

In all cases, the default vehicle cruising speed is $v = 25$ km/hr, and the passengers' value-of-time is $\beta = 20$ \$/hr, unless specified otherwise. Following \citet{Liu2021}, the vehicle operation cost is set as $\gamma = 40 + 0.48v = 52$ \$/hr, in which $40$ is the driver wages per hour, and $0.48$ is the prorated vehicle cost per km for vehicle acquisition, depreciation, and fuel \citep{zoepf2018economics}. 

\begin{table}[!t]
\centering
\begin{adjustbox}{width=1\textwidth}
\begin{tabular}{ccccccccccccccc}
\Xhline{4\arrayrulewidth}
\multirow{2}{*}{} & \multicolumn{2}{c}{Zone 1} & \multicolumn{2}{c}{Zone 2} & \multicolumn{2}{c}{Zone 3} & \multicolumn{2}{c}{Zone 4} & \multicolumn{5}{c}{System Characteristics}\\  \cmidrule(lr){2-3} \cmidrule(lr){4-5} \cmidrule(lr){6-7} \cmidrule(lr){8-9} \cmidrule(lr){10-14} 
& $n_{1000}$ & $\delta_{14\_2}$ ($\delta_{14\_3}$) & $n_{2000}$ & $\delta_{23\_1}$ ($\delta_{23\_4}$) & $n_{3000}$ & $\delta_{32\_1}$ ($\delta_{32\_4}$) & $n_{4000}$ & $\delta_{41\_2}$ ($\delta_{41\_3}$) & \multicolumn{1}{c}{$\sum_i M_i$} & \multicolumn{1}{c}{$M_b$} & \multicolumn{1}{c}{$M$} & \multicolumn{1}{c}{$P/(\sum\lambda_{ij})$} & \multicolumn{1}{c}{$Z$}\\
\Xhline{4\arrayrulewidth}
S1 & 11 & 1.00 (0.00) & 12 & 0.00 (1.00) & 10 & 0.49 (0.51) & 11 & 1.00 (0.00) & 1711 & 490 & 2200 & 0.35 & 21.25 \\\hline
S2 & 8 & 0.00 (1.00) & 15 & 0.00 (1.00) & 13 & 0.39 (0.61) & 7 & 0.00 (1.00) & 1718 & 244 & 1962 & 0.35 & 19.79 \\\hline
S3 & 8 & 0.14 (0.86) & 14 & 0.00 (1.00) & 14 & 0.00 (1.00) & 7 & 0.54 (0.46) & 1719 & 241 & 1959 & 0.35 & 19.78 \\
\Xhline{4\arrayrulewidth}
\end{tabular}
\end{adjustbox}
\caption[Caption for LOF]{Optimal service design for the hypothetical example with $2\times 2$ zone partition.
}
\label{tab: numerical hypo 4}
\end{table}

\subsection{Hypothetical Case}
Here, we consider a small square region with side length $10$ km. It is partitioned into $2\times2$ zones, each with a side length $\Phi = 5$ km, as shown in Figure \ref{fig:hypo region 4}. 
We apply the proposed model to a set of heterogeneous demand scenarios: (i) scenario S1: $\lambda_{i1} = \lambda_{i2}=\lambda_{i4} = 200$ trip/hr, and $\lambda_{i3} = 1400$ trip/hr, $\forall i\in\mathcal{K}$; (ii) scenario S2: $\lambda_{i1} =\lambda_{i4} = 200$ trip/hr, $\lambda_{i2} = 600$ trip/hr, and $\lambda_{i3} = 1000$ trip/hr, $\forall i\in\mathcal{K}$; and (iii) scenario S3: $\lambda_{i1} =\lambda_{i4} = 200$ trip/hr, $\lambda_{i2} = \lambda_{i3} = 800$ trip/hr, $\forall i\in\mathcal{K}$. 
In each of these scenarios, the total trip rate is $\sum_{i,j} \lambda_{ij} = 8,000$ trip/hr, and every zone has the same outgoing trip rate but different incoming trip rates.

\subsubsection{Optimal Design and System Performance}
\begin{figure}[t]
     \centering
     \begin{subfigure}[b]{0.32\textwidth}
     \raggedleft
     \includegraphics[width=0.85\textwidth]{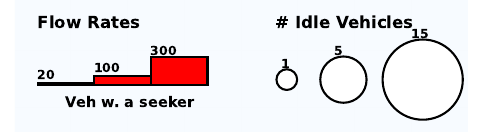}
     \end{subfigure}
     \begin{subfigure}[b]{0.32\textwidth}
     \raggedright
     \includegraphics[width=0.85\textwidth]{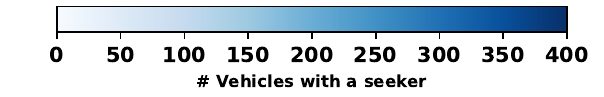}
     \end{subfigure}\\
     \begin{subfigure}[b]{0.32\textwidth}
         \centering
         \includegraphics[width=0.85\textwidth]{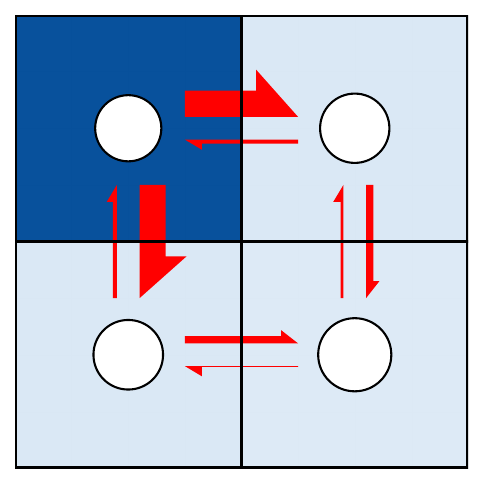}
         \caption{Scenario S1.}
         \label{fig:hypo S1}
     \end{subfigure}
     \begin{subfigure}[b]{0.32\textwidth}
         \centering
         \includegraphics[width=0.85\textwidth]{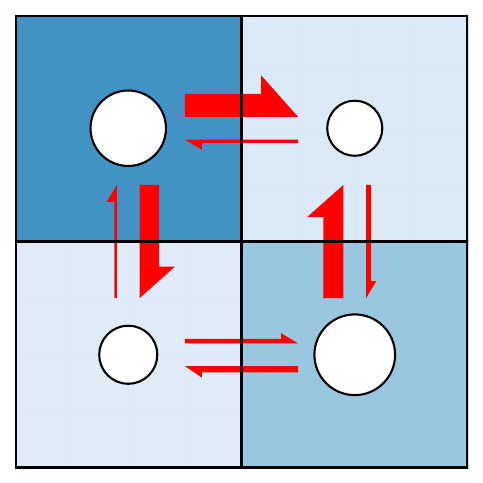}
         \caption{Scenario S2.}
         \label{fig:hypo S2}
     \end{subfigure}
     \begin{subfigure}[b]{0.32\textwidth}
         \centering
         \includegraphics[width=0.85\textwidth]{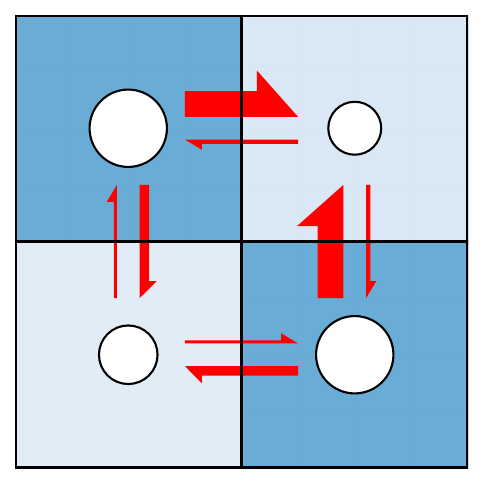}
         \caption{Scenario S3.}
         \label{fig:hypo S3}
     \end{subfigure}
    \caption{The hypothetical study region as $2\times 2$ zones, and illustration of idle vehicle count, seeker vehicle count, and cross-zone seeker flow in scenarios S1-S3.}
    \label{fig:hypo}
\end{figure}

The optimal deployment of idle vehicles and vehicle travel path decisions are summarized in Table \ref{tab: numerical hypo 4}.\footnote{These are the best solutions found by the proposed algorithm, which may or may not be the global optima due to non-convexity.} Vehicle counts are rounded to be integers, and we only display the vehicle path decisions between the diagonal zone pairs (for which vehicles have two feasible zone-level paths) --- those for other zone pairs are set by Constraint \eqref{eq:constraint delta} to take value 0 or 1. Figures \ref{fig:hypo S1}-\ref{fig:hypo S3} illustrate the numbers of idle vehicles (represented by the size of the white dots) and vehicles with a seeker (represented by the background color) in each zone. It also shows the flow rates of vehicles with an inter-zonal seeker (represented by the thickness of red arrows) between neighboring zones.

Scenario S1 mimics the morning rush in a city with a distinct CBD (i.e., zone 3). Despite the highly unbalanced travel demand, the difference in the numbers of idle vehicles in the four zones is not significant, which highlights the impact of ride-pooling. The high travel demand in zone 3 implies that the probability to match a seeker and a caller there is also high. Thus, zone 3 does not require a large number of idle vehicles; instead, there are about 350 vehicles with a seeker in zone 3 (as shown in Figure \ref{fig:hypo S1}), which are also suitable to serve new callers. In addition, since zone 2 is diagonally positioned to zone 3, passengers going from zone 1 (or zone 4) to zone 3 can take advantage of the ``passing" vehicles going from zone 2 towards zone 3. Hence, zone 1 and zone 4 require fewer idle vehicles than zone 2. By the symmetry of zone 1 and zone 4, it is reasonable that these two zones have a similar number of idle vehicles. 

It is also interesting to observe that the travel directions of vehicles are optimized as well. In Figure \ref{fig:hypo S1}, the flow of vehicles that carry a seeker out of zones 1, 2, and 4 is very small, possibly because most of these zones' inter-zonal passengers are going to zone 3, and hence it is likely that they will match up with each other before they leave their origin zone (i.e., only few of them leave their origin zone as a seeker). On the contrary, since the trip rate from zone 3 to zones 1, 2, and 4 is relatively low, the probability to match two inter-zonal passengers inside zone 3 is relatively low, and in turn, many of them leave zone 3 as a seeker. For the vehicle travel path decisions, $\delta_{32\_1} =0.49$ and $\delta_{32\_4}=0.51$, which makes sense as zone 1 and zone 4 are symmetric, while $\delta_{23\_1} = 0.00$ and $\delta_{23\_4} = 1.00$, possibly because the flow of vehicles with a seeker leaving zone 2 is too small to affect the matching in zones 1 and 4. 
In addition, the vehicles with a seeker going between zone 1 and zone 4 are routed to visit zone 2, i.e., $\delta_{14\_2}= \delta_{41\_2}=1.00,\ \delta_{14\_3}=\delta_{41\_3} = 0.00$, because zone 3 is the major destination zone with many suitable vehicles inside, and it is beneficial to route these seekers through zone 2 so as to better serve the passengers there.

In Scenario S2, a second, yet smaller, CBD emerges in zone 2 and diverts some of the trips away from zone 3. It is expected that more idle vehicles are now observed in zone 2, as more delivery trips end there and become idle. Interestingly, the number of idle vehicles in zone 3 also increases (as compared to Scenario S1), even though fewer 
vehicles now enter zone 3 to deliver passengers. A possible explanation is that the number of vehicles with a seeker has also decreased in zone 3 (e.g., by comparing the background color of zone 3 in Figures \ref{fig:hypo S2} and \ref{fig:hypo S1}), which increases its needs for idle vehicles. In the meantime, the numbers of idle vehicles in zone 1 and zone 4 both decrease, possibly because a more even trip distribution between zone 2 and zone 3 increases the probability to match passengers in zone 1 and zone 4 with seekers going between zone 2 and zone 3. Correspondingly, there is a smaller flow of vehicles with a seeker leaving zone 1 and zone 4. In addition, zone 1 now has one more vehicle than zone 4, and it is interesting to observe that more seekers going between zone 2 and zone 3 are routed through zone 4 (i.e., $\delta_{23\_4} = 1.00$, $\delta_{32\_4} = 0.61$), possibly to compensate for the fewer idle vehicles in zone 4. All seekers going between zone 1 and zone 4 are now routed through zone 3 (i.e., $\delta_{14\_3} = \delta_{41\_3} = 1.00$), which makes sense since passengers in zone 2 are well served by a large number of idle vehicles, while zone 3 experienced a loss of seekers (as zone 2 attracts more demand). 

As the CBD in zone 2 continues to grow and eventually becomes symmetrical to that in zone 3, as in Scenario S3, the numbers of idle vehicles in zone 2 and zone 3 are now the same, as expected. The rest of the optimal design is similar to that under S2.

Table \ref{tab: numerical hypo 4} also shows that, as demand distribution goes from monocentric (S1) to bicentric (S3), the total number of active vehicles increases but the number of rebalanced vehicles decreases. This is reasonable too. The trip destinations are spatially highly unbalanced in the monocentric scenario, which imposes a huge burden on vehicle rebalancing activities. However, in the meantime, since the trip destinations are concentrated in zone 3, the passenger matching probability is high and the ride-pooling service can be operated efficiently. As the trip destinations become more ``scattered'' in the bicentric scenario, the matching probability decreases, while the rebalancing operation is reduced as well. Overall, we see that the optimal fleet needed in the entire region, $M$, changes from 2,200 in Scenario S1, to 1,962 in Scenario S2, and to 1,959 in Scenario S3. On the other hand, the average passenger door-to-door travel time $P/(\sum\lambda_{ij})$ does not change significantly. As a result, the average system cost decreases from S1 to S3.

\subsubsection{Sensitivity to Zone Size}
\begin{figure}[t]
     \centering
     \begin{subfigure}[c]{0.32\textwidth}
     \raggedleft\includegraphics[width=0.85\textwidth]{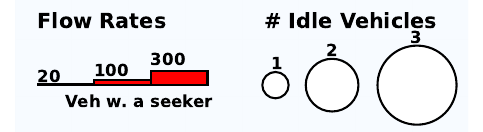}
     \end{subfigure}
     \begin{subfigure}[c]{0.32\textwidth}
     \raggedright\includegraphics[width=0.85\textwidth]{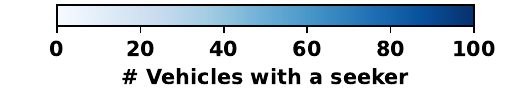}
     \end{subfigure}\\
     \begin{subfigure}[b]{0.32\textwidth}
         \centering\includegraphics[width=0.85\textwidth]{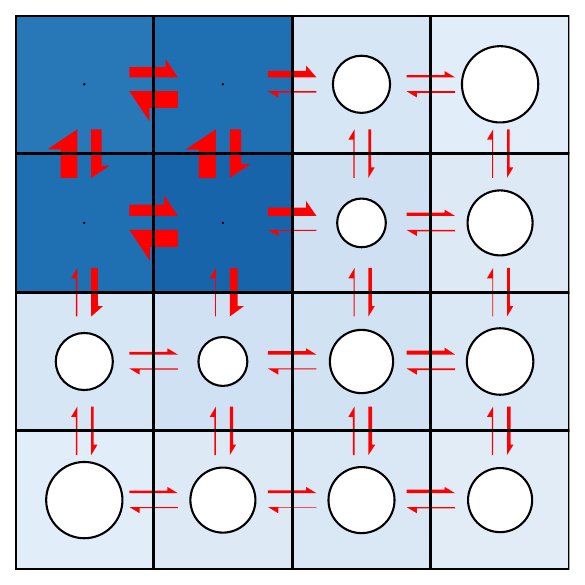}
         \caption{Scenario S1.}
         \label{fig:hypo S1 16}
     \end{subfigure}
     \begin{subfigure}[b]{0.32\textwidth}
         \centering
         \includegraphics[width=0.85\textwidth]{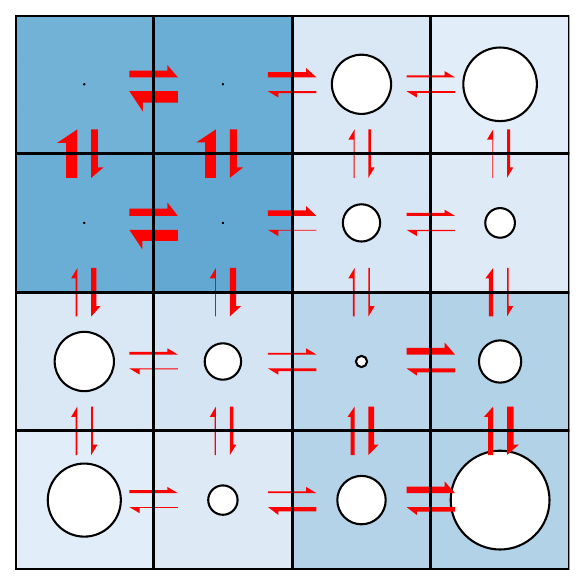}
         \caption{Scenario S2.}
         \label{fig:hypo S2 16}
     \end{subfigure}
     \begin{subfigure}[b]{0.32\textwidth}
         \centering
         \includegraphics[width=0.85\textwidth]{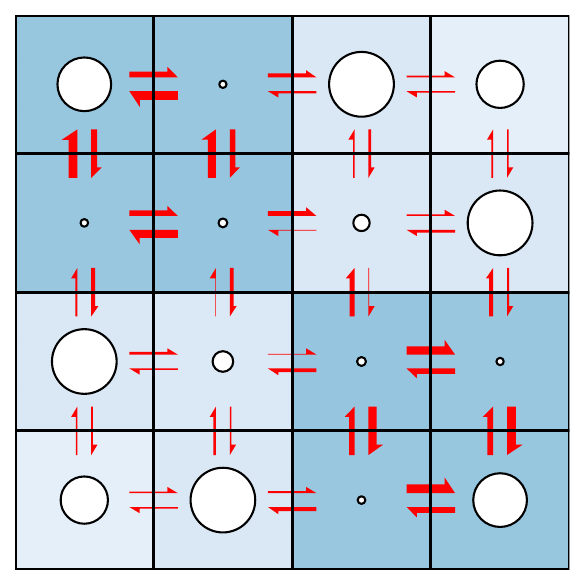}
         \caption{Scenario S3.}
         \label{fig:hypo S3 16}
     \end{subfigure}
    \caption{The hypothetical study region as $4\times 4$ zones, and illustration of idle vehicle count, seeker vehicle count, and cross-zone seeker flow in scenarios S1-S3.}
    \label{fig:hypo 16}
\end{figure}

In this subsection, we explore the impacts of zone partitions on system performance, by solving the same model when the region is partitioned into $4\times 4$ zones with side length $\Phi = 2.5$ km; see Figure \ref{fig:hypo region 16}. The same demand scenarios equivalent to S1-S3 are used; e.g., for S1, we use $\lambda_{ij} = 87.5$ trip/hr, $\forall i\in\mathcal{K},\ j\in\{9, 10, 13, 14\}$, and $\lambda_{ij} = 12.5$ trip/hr, $\forall i\in\mathcal{K},\ j\in \mathcal{K}\setminus\{9, 10, 13, 14\}$
. The optimal designs from the model, including the distribution of idle vehicles and vehicles with a seeker, and the flow of inter-zonal seekers between neighboring zones, are illustrated in Figure \ref{fig:hypo 16}. 
The results show similar patterns to those observed with $2\times 2$ zones; e.g., high-demand CBD zones (e.g., zones 9, 10, 13, and 14) require negligible idle vehicles; some zones see large volumes of passing seekers and hence can exploit them to serve local demand (e.g., $n_{6000}\leq n_{1000}$ and $n_{7000}\leq n_{4000}$ for S1-S3); and the optimal allocation of idle vehicles shows symmetry (e.g., $n_{1000} \approx n_{16,000}$, $n_{2000} \approx n_{12,000}$, and $n_{3000} \approx n_{8000}$ for S1-S3).

\begin{figure}
    \centering
    \includegraphics[width = \textwidth]{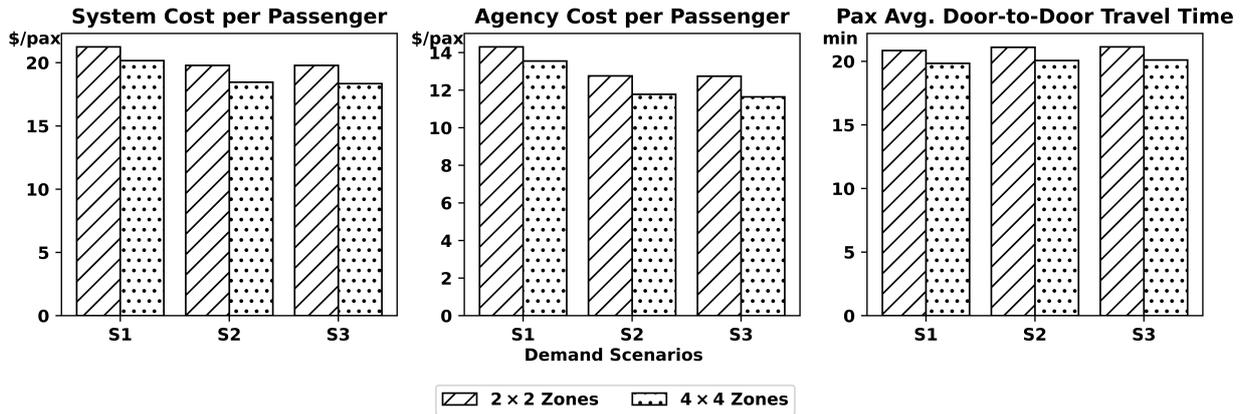}
    \caption{Comparison of the system performance with different zone partitions.}
    \label{fig:hypo partition comparison}
\end{figure}

The system's performance with $4\times 4$ zones is compared with that with $2\times 2$ zones in Figure \ref{fig:hypo partition comparison}. The overall trends across the three demand scenarios are similar; e.g., the system-wide cost and agency cost decrease as the demand distribution becomes more uniform, while the average passenger travel time remains relatively stable. However, as compared to the $2\times 2$ zone partition, $4\times 4$ zone partition reduces the system-wide cost by $5.12\%\sim7.28\%$, saves $5.25\%\sim8.61\%$ agency cost, and yields about $4.8\%$ lower average passenger travel time.
The results indicate that the ride-pooling service could achieve better service performance when the zones are smaller. This observation is reasonable, since the length of a detour from each seeker-caller match is bounded by the zone size --- as zones become smaller, the detour length per match also decreases, and in turn, the service can be operated more efficiently. However, it should be noted that the passenger travel experience could also deteriorate as zone size decreases, because some previous intra-zonal trips could become inter-zonal and then be feasible to possible matches that induce detours. Hence, the zone size may have complex effects on the performance of the ride-pooling service, and the impacts may depend heavily on the distribution of trip ODs and the matching strategies/criteria. 

The choice of zone size also has other implications. A smaller zone size could achieve better demand homogeneity within each zone and hence allow the model to capture the real-world demand distribution more accurately. Yet, we cannot partition the study region into arbitrarily small zones due to two reasons. First, if zones are too small, some zones may run out of suitable vehicles to serve passengers from time to time, which violates the model assumption that the TNC instantly assigns every passenger to a suitable vehicle. Second, smaller zones increase the computational burden for solving the system of nonlinear equations \eqref{eq:side directions}-\eqref{obj}. If we have $|\mathcal{K}|$ zones in a study region, the numbers of $n_{i0j0}$ and $n_{i0jk}$ variables are $O(|\mathcal{K}|^2)$ and $O(|\mathcal{K}|^3)$, respectively. On a laptop with Intel i7-7500U CPU @ 2.70GHz and 16GB memory, the average time to solve Eqn. \eqref{eq:side directions}-\eqref{obj} grow rapidly from about 0.2 seconds for $2\times 2$ zones to about 80 seconds for $4\times 4$ zones. Therefore, the zone size should be chosen carefully to balance model validity, data accuracy, and computational tractability.

\subsubsection{Simulation Verification}\label{sec:simulation}
In this subsection, we verify the accuracy of the proposed analytical model with an agent-based simulation platform \citep{Zhai2023}. 
The simulation inputs include: (i) demand characteristics, e.g., study region’s spatial geometry, demand rates, zone layout and size; and (ii) service supply characteristics, e.g., the total number of vehicles in each zone, vehicle speed, zone-level OD movement paths, and zone-to-zone rebalancing rates. Vehicles can travel on a dense grid network and make turns at any location. The passenger trips between each zone pair are generated from a time-homogeneous Poisson process, and their origins and destinations are uniformly distributed respectively within the corresponding zones. At the beginning of each simulation, all vehicles are idle and randomly distributed within each zone. When a passenger submits a travel request, the platform scans through all vehicles within its origin zone to identify the nearest suitable one based on the matching criteria. Due to randomness, it is possible that when a passenger calls for service, there are no suitable vehicles; in this case, that passenger is removed from the system and considered lost. A delivery trip follows the suggested zone-level path and the routing rules, and the simulation platform rebalances a random idle vehicle between two zones according to a Poisson process of the corresponding rate. 
In every two seconds, the simulation checks passengers’ travel requests and rebalancing operations, and updates vehicle states and locations. Each simulation’s duration is set to be 10 hours of operations (serving approximately 80,000 trips), and the first 1/3 of the simulation is considered as the warm-up period. During the remainder 2/3 of the simulation, the platform keeps track of all performance metrics; e.g., zonal vehicle distribution by type, passenger door-to-door travel time, and the number of lost passengers. 

This simulation platform is then used to corroborate the analytical model from Section \ref{sec:methodology} under a variety of fleet sizes, demand scenarios (S1-S3), and zone partition granularities ($2\times2$ vs. $4\times 4$ zones). To vary the fleet size for each zone partition and demand scenario, we first use the analytical model to obtain the optimal allocation of idle vehicles $\{n_{i000}^*:\forall i\}$ and zone-level vehicle movements $\{\delta^*_{ij\_i'}:\forall i, i',j\}$. We then find from the same model a new fleet size $\{M_i\}$ and rebalancing rates $\{b_{ij}\}$  based on the same $\{\delta^*_{ij\_i'}\}$ but a different number of idle vehicles at equilibrium $\{n_{i000}\}=q\cdot\{n_{i000}^*\}$, where the multiplicative factor $q\in[0.5,100]$ indicates “redundancy” of idle vehicles in the system. This step also predicts system performance statistics (e.g., average travel time) via the analytical model. Now, for each zone partition, demand scenario, and fleet size, we feed the corresponding $\{M_i\}$, $\{\delta^*_{ij\_i'}\}$ and $\{b_{ij}\}$  into the simulation platform to evaluate the system performance statistics (e.g., average travel time and percentage of lost demand). 

\begin{figure}[t]
     \centering
     \begin{subfigure}[b]{0.49\textwidth}
         \centering
         \includegraphics[width=\textwidth]{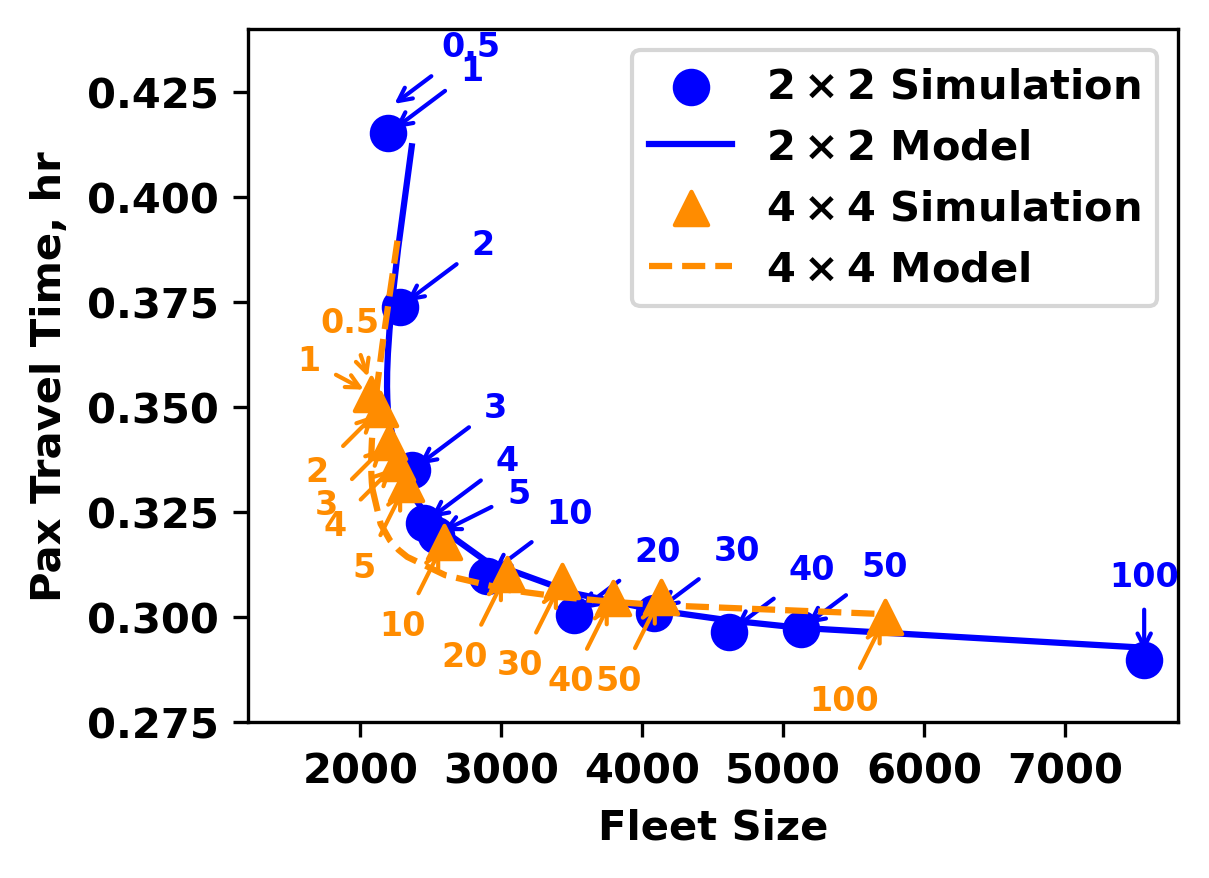}
         \caption{Travel time in scenario S1.}
         \label{fig:sim S1}
     \end{subfigure}
     \begin{subfigure}[b]{0.49\textwidth}
         \centering
         \includegraphics[width=\textwidth]{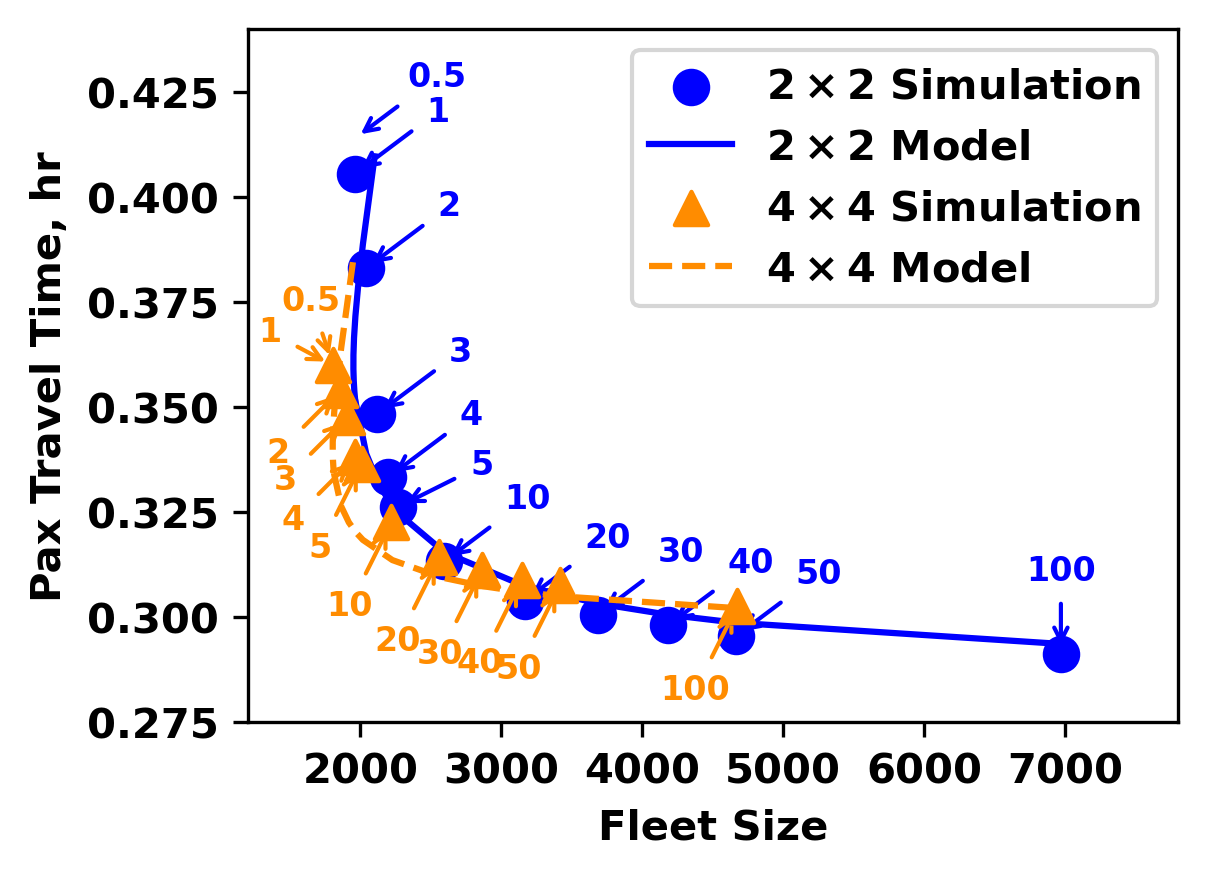}
         \caption{Travel time in scenario S2.}
         \label{fig:sim S2}
     \end{subfigure}
     \begin{subfigure}[b]{0.49\textwidth}
         \centering
         \includegraphics[width=\textwidth]{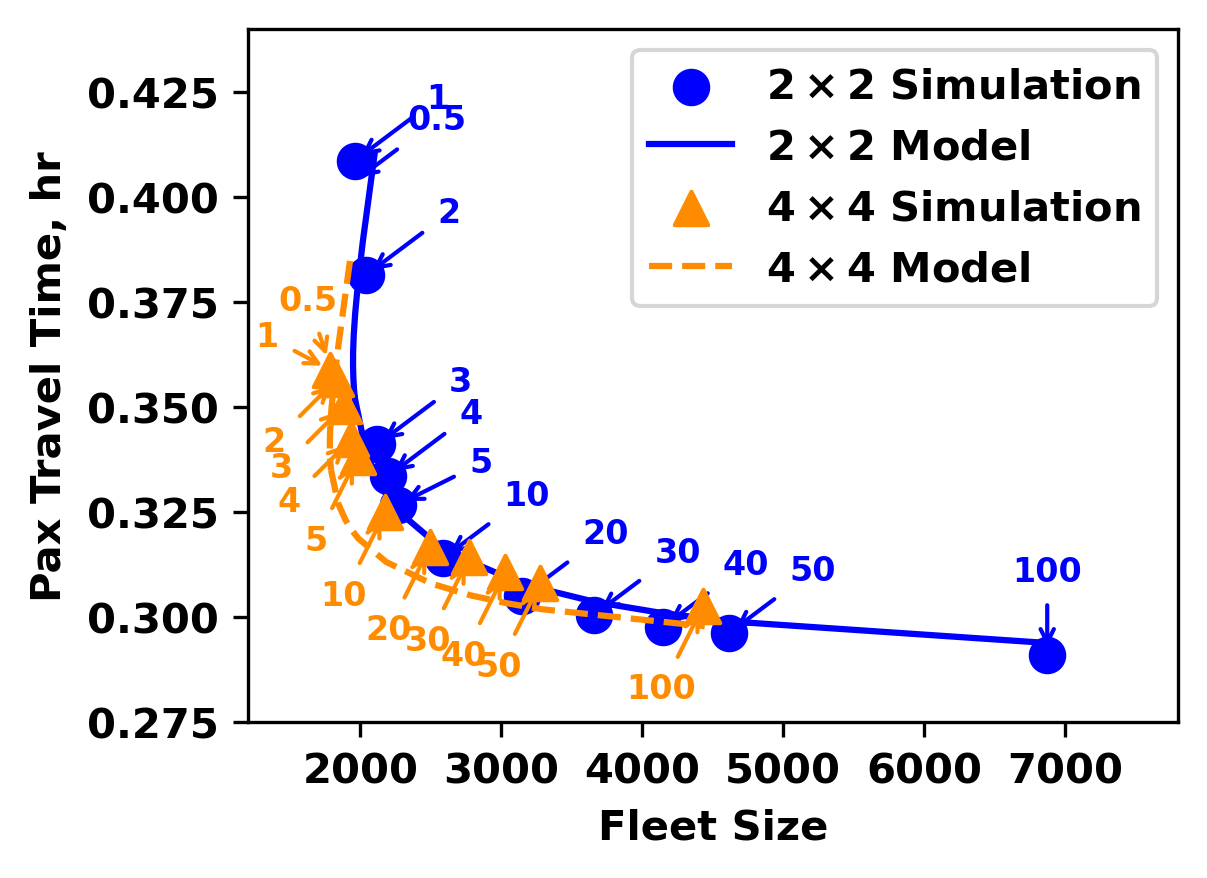}
         \caption{Travel time in scenario S3.}
         \label{fig:sim S3}
     \end{subfigure}
     \begin{subfigure}[b]{0.49\textwidth}
         \centering
         \includegraphics[width=\textwidth]{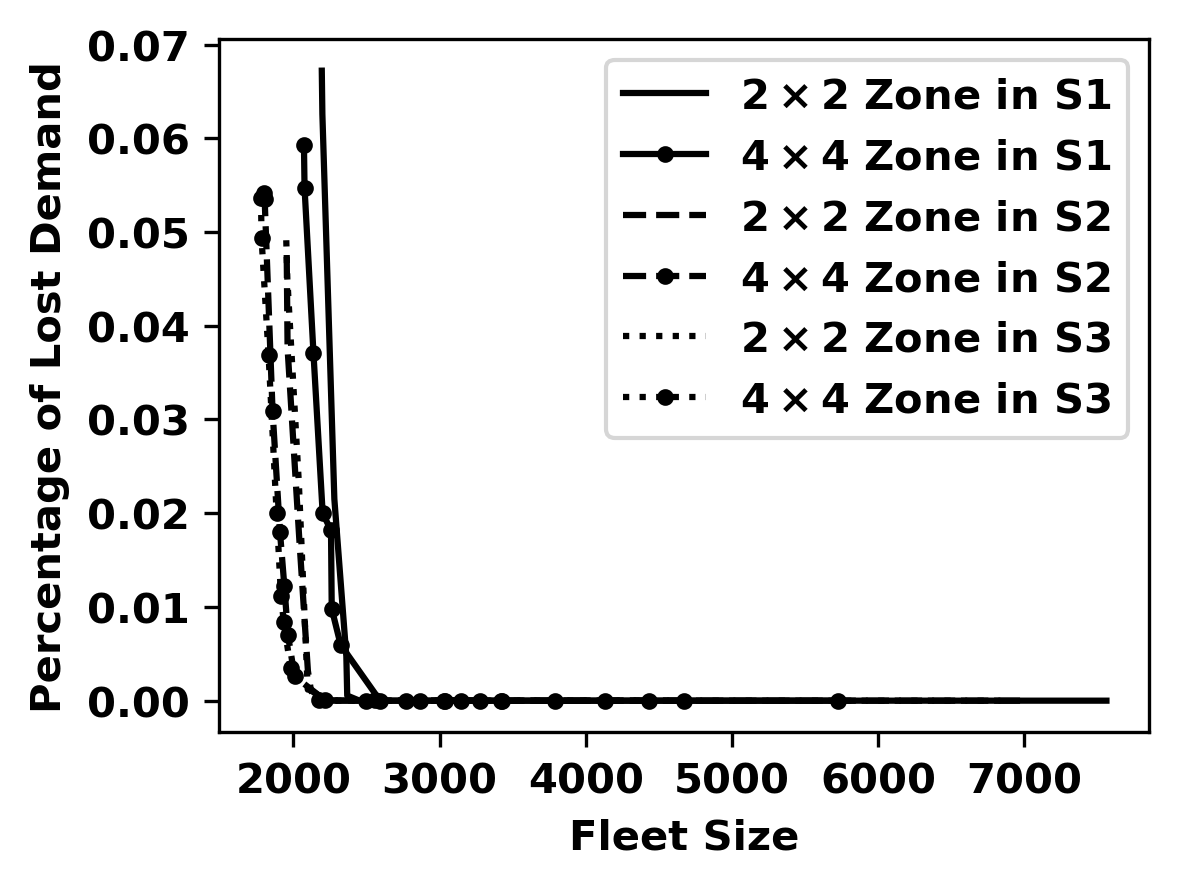}
         \caption{Percentage of lost demand from simulation.}
         \label{fig:sim demand}
     \end{subfigure}
    \caption{Comparisons between model results and simulation results.}
    \label{fig:sim}
\end{figure}

Figures \ref{fig:sim S1}-\ref{fig:sim S3} show the average passenger travel times from the analytical model (represented by curves), and those from the simulation (represented by markers), for all demand scenarios, fleet sizes, and zone partitions. The $q$ value is labeled near the corresponding marker. Overall, the simulation results and the model predictions fit very well for almost all cases, and the error gap tends to shrink as the fleet size increases. 

In addition, the backward bending parts of the model curves under small fleet sizes (e.g., around 2,000) illustrate possible presence of WGC, in which the passenger travel time increases as the fleet size increases. WGC is also observed in the simulations for small $q$ values, as the number of idle vehicles approaches zero in every zone due to randomness. In other words, the simulations under these $q$ values seem to converge to the inefficient equilibrium and inflate the measured passenger travel times
. Accordingly, as Figure \ref{fig:sim demand} shows, there is about 5\%-7\% of lost demand when the fleet size is less than 2,200, possibly because the system performance becomes very sensitive to the fleet size, and the naïve vehicle rebalancing strategy in the simulation platform (rebalancing at a constant rate) could not push the system to the efficient equilibrium.

However, as shown in Figure \ref{fig:sim}, as long as $q\geq3$, the simulated system performance always lies in an efficient equilibrium. 
Therefore, to enhance the system stability in the real-world implementation, we recommend that TNCs increase the optimal number of idle vehicles in each zone (normally a very small fraction of the fleet) slightly by a factor around 3 -- this option does not significantly increase the total fleet size but notably improve the system's performance against WGC. We may also adopt more advanced operation strategies to mitigate the WGC, such as the block matching strategy \citep{feng2022approximating}, adaptive matching radius \citep{xu2020supply}, or dynamic vehicle swaps \citep{shen4264716dynamic,ouyang2023measurement}. The integration of these advanced strategies with the proposed multi-zone modeling framework could be interesting topics for future research.

\subsection{Chicago Case Study}\label{sec:chicago}
\begin{figure}[t]
     \centering
     \begin{subfigure}[b]{0.4\textwidth}
         \centering
         \includegraphics[width=\textwidth]{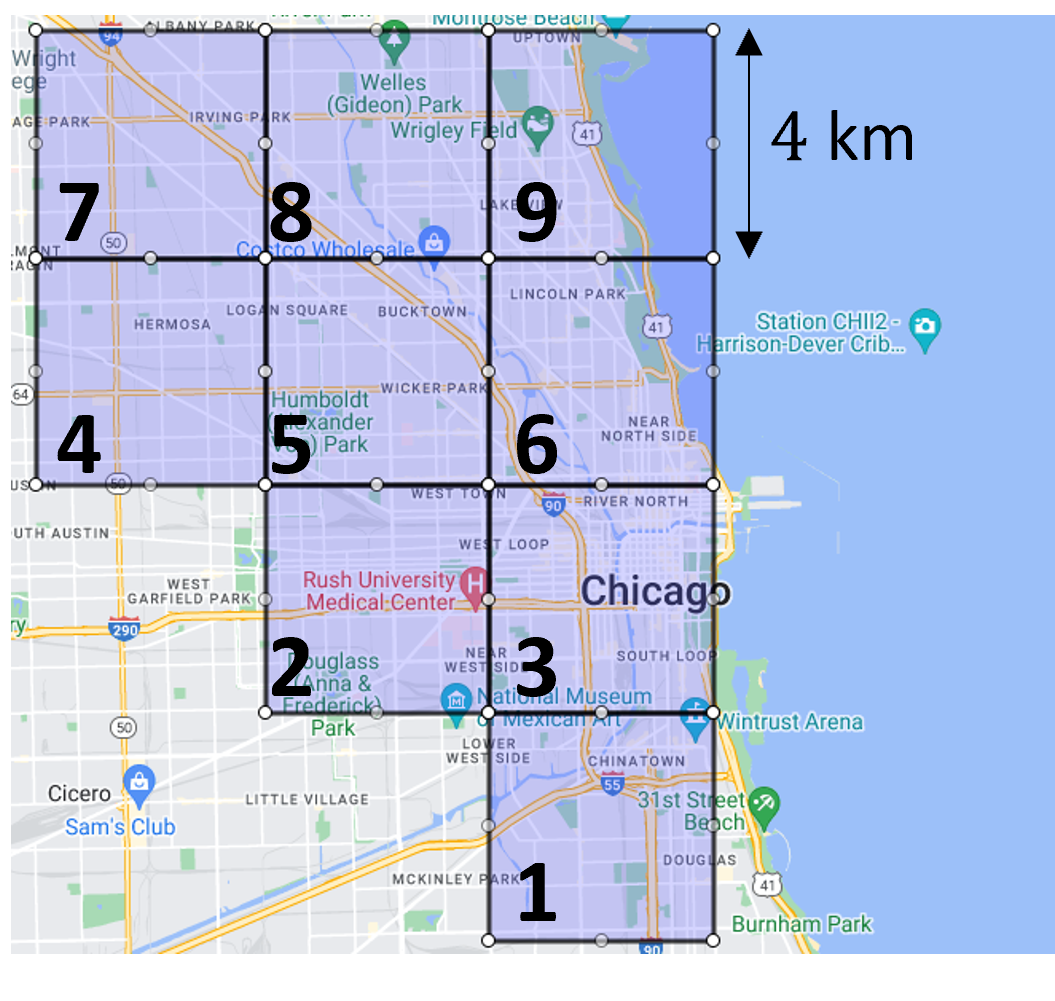}
         \caption{Partition of the study region.}
         \label{fig:chicago region}
     \end{subfigure}
     \begin{subfigure}[b]{0.505\textwidth}
         \centering
         \includegraphics[width=\textwidth]{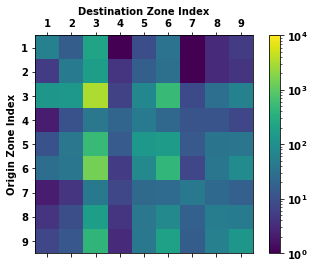}
         \caption{Trip rates of OD zone pair $\{\lambda_{ij}\}$.}
         \label{fig:chicago demand}
     \end{subfigure}
    \caption{Zone partition and taxi/TNC trip OD matrix of the Chicago case study.}
    \label{fig:demand}
\end{figure}

\begin{figure}[t]
     \centering
     \begin{subfigure}[b]{0.3\textwidth}
         \centering
         \includegraphics[width=\textwidth]{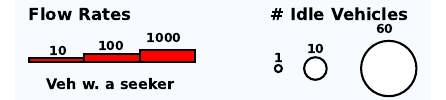}
     \end{subfigure}
     \begin{subfigure}[b]{0.3\textwidth}
         \centering
         \includegraphics[width=\textwidth]{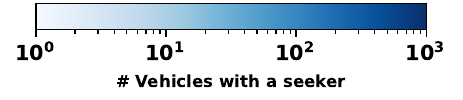}
     \end{subfigure}\\
     \begin{subfigure}[b]{0.3\textwidth}
         \centering
         \includegraphics[width=\textwidth]{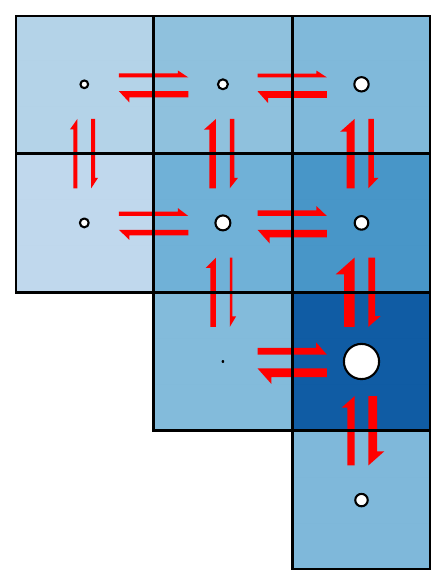}
         \caption{$q = 1,\ \beta = 10$ \$/hr.}
         \label{fig:chicago flow1}
     \end{subfigure}
     \begin{subfigure}[b]{0.3\textwidth}
         \centering
         \includegraphics[width=\textwidth]{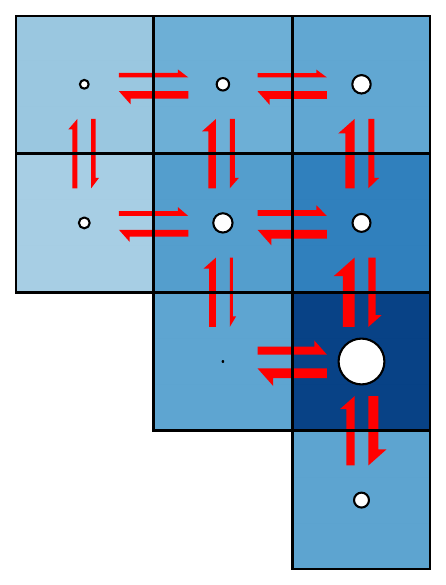}
         \caption{$q = 2,\ \beta = 10$ \$/hr.}
         \label{fig:chicago flow2}
     \end{subfigure}
     \begin{subfigure}[b]{0.3\textwidth}
         \centering
         \includegraphics[width=\textwidth]{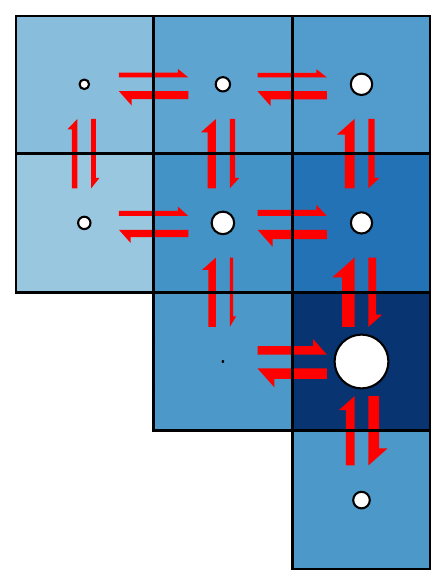}
         \caption{$q = 3,\ \beta = 10$ \$/hr.}
         \label{fig:chicago flow3}
     \end{subfigure}
    \caption{Illustration of idle vehicle count, seeker vehicle count, and cross-zone vehicle flow for the Chicago case study.}
    \label{fig:chicago flow}
\end{figure}

In this section, we apply the proposed model to a subarea of Chicago, which is partitioned into 9 zones with $\Phi = 4$ km, as shown in Figure \ref{fig:chicago region}. We extract TNC and taxi trip data for the morning rush (7:00 AM $\sim$ 9:00 AM) in the entire year of 2019 from the Chicago Data Portal \citep{cityofchicagoTaxi,cityofchicagoTNC}, and aggregate them into hourly zone-to-zone trip rates; see Figure \ref{fig:chicago demand}. The total trip rate for the entire study region is 9,520 trips per hour. The trips from and to those zones along the downtown and lake shore areas (e.g., zone 3, 6, 5, and 9) are orders of magnitude greater than those for others; e.g., the hourly incoming trip rate is 6,196 trip/hr for zone 3, but only 64 trip/hr for zone 4. Because the 2019 data may not be exhaustive, and because the TNC service is growing very rapidly, we follow \citet{lei2019path} and consider scaling all trip rates by a common factor $q\geq 1$. 

We solve the optimal service design and evaluate the system performance while varying the passenger value-of-time $\beta\in\{10, 50, 100\}$ \$/hr and the demand scaling factor $q\in\{1, 2, 3\}$. All other parameters are the same as those in the hypothetical example.

\begin{table}[!t]
\centering
\begin{adjustbox}{width=1\textwidth}
\begin{tabular}{cccccccccccccc}
\Xhline{4\arrayrulewidth}
$\beta$ & $q$  & $n_{1000}$ & $n_{2000}$ & $n_{3000}$ & $n_{4000}$ & $n_{5000}$ & $n_{6000}$ & $n_{7000}$ & $n_{8000}$ & $n_{9000}$ & $\sum_i M_i$ & $M_b$ &  $M$ \\
\Xhline{4\arrayrulewidth}
\multirow{3}{*}{10}  & 1 & 3  & 0 & 24  & 2  & 4  & 4  & 1  & 2  & 4  & 1652 & 256 & 1909 \\\cmidrule(lr){2-14} 
                     & 2 & 4  & 0 & 41  & 2  & 7  & 6  & 2  & 3  & 7  & 3191 & 511 & 3702 \\\cmidrule(lr){2-14} 
                     & 3 & 5  & 0 & 57  & 3  & 10 & 9  & 2  & 4  & 9  & 4707 & 767 & 5474 \\\hline
\multirow{3}{*}{50}  & 1 & 5  & 1 & 91  & 3  & 8  & 8  & 3  & 4  & 8  & 1741 & 264 & 2005 \\\cmidrule(lr){2-14} 
                     & 2 & 10 & 1 & 168 & 5  & 12 & 12 & 5  & 6  & 12 & 3359 & 513 & 3872 \\\cmidrule(lr){2-14} 
                     & 3 & 15 & 0 & 242 & 6  & 17 & 16 & 6  & 8  & 15 & 4953 & 760 & 5713 \\\hline
\multirow{3}{*}{100} & 1 & 8  & 5 & 163 & 6  & 13 & 20 & 6  & 8  & 14 & 1868 & 287 & 2155 \\\cmidrule(lr){2-14} 
                     & 2 & 11 & 7 & 307 & 8  & 20 & 30 & 9  & 12 & 21 & 3591 & 550 & 4141 \\\cmidrule(lr){2-14} 
                     & 3 & 17 & 7 & 446 & 10 & 26 & 38 & 11 & 15 & 27 & 5285 & 810 & 6095\\
\Xhline{4\arrayrulewidth}
\end{tabular}
\end{adjustbox}
\caption{Optimal design for the Chicago study region under varying $q$ and $\beta$.}
\label{tab: numerical chicago results}
\end{table}

\begin{figure}[t]
    \centering
    \includegraphics[width=\textwidth]{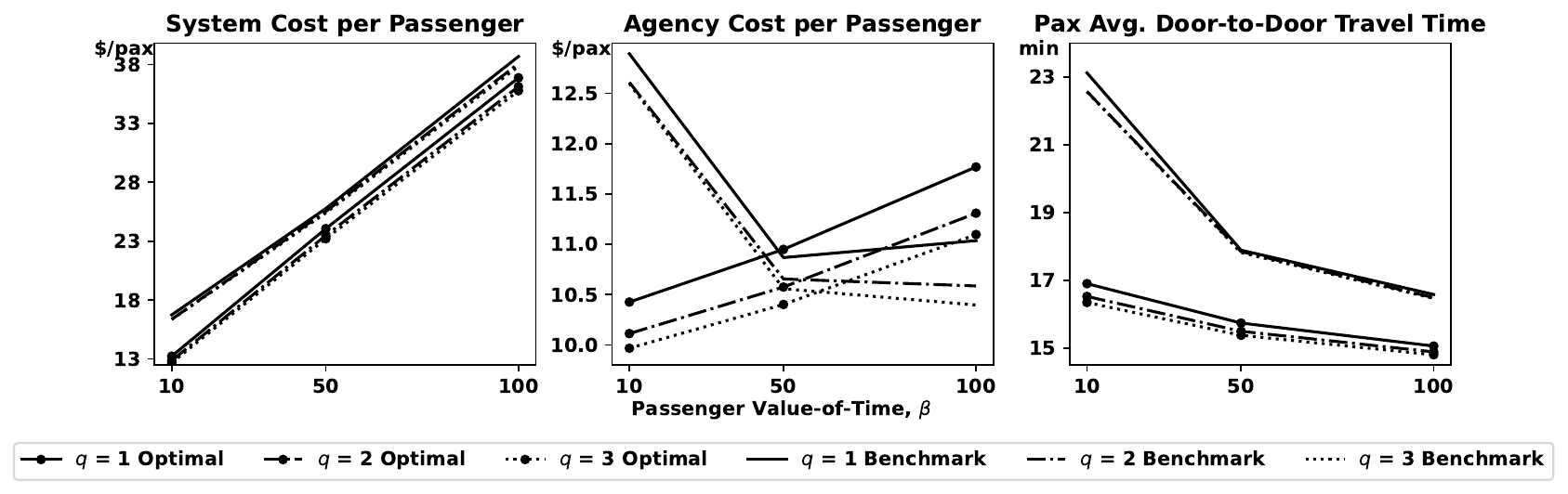}
    \caption{System performance under the optimal design and the benchmark design when varying $q$ and $\beta$.}
    \label{fig:chicago results}
\end{figure}

\begin{table}[!t]
\centering
\small
\begin{adjustbox}{width=1\textwidth}
\begin{tabular}{cccccccccccccc}
\Xhline{4\arrayrulewidth}
$\beta$ & $q$  & $n_{1000}$ & $n_{2000}$ & $n_{3000}$ & $n_{4000}$ & $n_{5000}$ & $n_{6000}$ & $n_{7000}$ & $n_{8000}$ & $n_{9000}$ & $\sum_i M_i$ & $M_b$ &  $M$ \\
\Xhline{4\arrayrulewidth}
\multirow{3}{*}{10}  & 1 & 8  & 0 & 3  & 1  & 0 & 0 & 7  & 0  & 0 & 2074 & 287 & 2361 \\\cmidrule(lr){2-14} 
                     & 2 & 13 & 0 & 5  & 2  & 0 & 0 & 10 & 0  & 0 & 4044 & 573 & 4618 \\\cmidrule(lr){2-14} 
                     & 3 & 18 & 0 & 5  & 2  & 0 & 0 & 14 & 0  & 0 & 6058 & 861 & 6919 \\\hline
\multirow{3}{*}{50}  & 1 & 20 & 2 & 8  & 5  & 1 & 1 & 16 & 3  & 1 & 1708 & 282 & 1990 \\\cmidrule(lr){2-14} 
                     & 2 & 33 & 2 & 12 & 6  & 1 & 1 & 24 & 3  & 1 & 3330 & 571 & 3901 \\\cmidrule(lr){2-14} 
                     & 3 & 45 & 2 & 14 & 7  & 1 & 1 & 31 & 3  & 1 & 4939 & 860 & 5799 \\\hline
\multirow{3}{*}{100} & 1 & 32 & 6 & 15 & 11 & 5 & 4 & 29 & 8  & 6 & 1732 & 288 & 2021 \\\cmidrule(lr){2-14} 
                     & 2 & 55 & 7 & 23 & 16 & 6 & 5 & 46 & 11 & 7 & 3313 & 563 & 3876 \\\cmidrule(lr){2-14} 
                     & 3 & 75 & 8 & 29 & 20 & 6 & 5 & 60 & 13 & 8 & 4868 & 842 & 5710\\
\Xhline{4\arrayrulewidth}
\end{tabular}
\end{adjustbox}
\caption{Benchmark design for the Chicago study region under varying $q$ and $\beta$.}
\label{tab: numerical chicago results benchmark}
\end{table}

Table \ref{tab: numerical chicago results} summarizes the optimal number of idle vehicles in each zone for all combinations of $q$ and $\beta$. When $q = 1$ and $\beta = 10$ \$/hr, the deployment of idle vehicles generally follows the demand distribution, as expected; e.g., zones 3, 5 and 9 have larger numbers, and zones 4 and 7 have smaller numbers. However, it is interesting to observe some exceptions. For example, zone 1 and zone 6 have a similar number of idle vehicles although the demand rate of zone 6 is much greater than that of zone 1, and zone 2 has moderate travel demand but idle vehicles are seldom allocated to zone 2. 
This discrepancy reflects the joint impacts of geometry and ride-pooling. For zone 1, since no vehicles from other zones would pass it on their delivery routes, the service in zone 1 heavily relies on idle vehicles inside. On the contrary, zones 2 and 6 
are between high-demand zone 3 and all northern and western zones (i.e., zones 4, 5, 7, 8, and 9), and hence the passengers in zones 2 and 6 can be served partially by passing vehicles. 

As $q$ increases, it is noted from Table \ref{tab: numerical chicago results} that the number of idle vehicles in high-demand zone 3 increases almost linearly with $q$. It is possible that even when $q=1$, passengers in zone 3 are well matched and ride-pooling opportunities have been fully utilized. Therefore, when demand increases, the number of idle vehicles in zone 3 must also increase. On the other hand, the idle vehicle counts in low-demand zones increase only mildly (or do not increase at all), probably because these zones can still benefit from a higher matching probability under higher demand. Overall, both the total number of active vehicles and that of rebalanced vehicles increase with $q$, almost proportionally, since there are more trips to serve and the absolute unevenness in demand distribution becomes more severe. More information about the vehicle counts and the flow between neighboring zones are visualized in Figure \ref{fig:chicago flow}, which shows that 
both the number of vehicles with a seeker in each zone and the flow of vehicles with a seeker leaving each zone increase with $q$. 

As $\beta$ increases, the system tends to emphasize more on passenger experience, and hence it is reasonable to observe (slightly) more idle vehicles for almost every zone. As a whole, the total number of active vehicles also increases, but only moderately, and almost all the increases come from idle vehicles -- as these vehicles most directly contribute to better passenger experience. Similarly, when the demand remains the same, the increasing value of $\beta$ does not affect the total number of rebalanced vehicles significantly, as expected. 

Figure \ref{fig:chicago results} illustrates how the optimal system performance varies with $q$ and $\beta$. A larger value of $q$ always yields a lower average system cost, average agency cost, and average passenger travel time; this exhibits economies of scale, as expected. Given $q$, the average passenger travel time decreases with $\beta$, while the agency cost increases --- overall, the average system-wide cost increases. 

Finally, we highlight the importance of explicitly addressing demand heterogeneity in mobility service design. We obtain a benchmark design for Chicago by optimizing the allocation of idle vehicles and vehicles' zone-level movements under the false assumption of homogeneous demand, i.e., the demand rate between every pair of zones is $\frac{\sum_{i,j}\lambda_{ij}}{|\mathcal{K}|^2}$. The benchmark designs under varying $q$ and $\beta$ values are summarized in Table \ref{tab: numerical chicago results benchmark}. The performance of these benchmark designs is then evaluated under the actual heterogeneous demand and compared with that of the corresponding optimal designs. Figure \ref{fig:chicago results} shows that the optimal design incurs a $20.93\%\sim22.42\%$ lower system-wide cost and a $26.91\%\sim27.57\%$ lower passenger travel time when $\beta = 10$\$/hr, and the improvement is more significant for larger $q$ values. 
The results indicate clearly that failure to address demand heterogeneity could lead to notably suboptimal service designs.

It is particularly interesting to see that the agency cost (or fleet size) of the benchmark design is $23.67\%\sim26.41\%$ higher than that of the optimal design when $\beta = 10$ \$/hr, but it decreases drastically as $\beta$ increases and even becomes smaller than that of the optimal design when $\beta =50$ and $\beta =100$ \$/hr. A possible reason is that the numbers of suitable vehicles in the benchmark design are very low in high-demand zones; as a result, 
the high-demand zones likely suffer heavily from WGC and a large number of vehicles must be deadheading to pick up passengers, which mandates a very large fleet size and incurs an extremely high agency cost. As $\beta$ increases, the number of idle vehicles in each zone increases, and therefore, WGC is less likely to appear, and the total fleet size of the benchmark design decrease. When $\beta$ is very large, e.g., $\beta = 100$\$/hr and $q = 1$, a large number of idle vehicles under the benchmark design leads to an increase in the total fleet size.

\section{Conclusion}\label{sec:conclusion}
This paper proposed a multi-zone queuing network model to describe how TNCs shall plan ride-pooling services for areas with spatially heterogeneous demand. 
The model incorporates TNC's decisions, including the deployment of vehicles in each zone, vehicle travel directions between zones, and the rate to rebalance vehicles between zones. A set of criteria are imposed to avoid excessive detours for matched passengers. The steady-state operation of the ride-pooling service is modeled by a multi-zone aspatial queuing network, composed of a set of vehicle states that describe vehicles' location and workload, and a set of flow variables that describe vehicles' transitions within and across zones. A large system of nonlinear equations is formulated based on the queuing network model, and its solution can be used to evaluate the expected system-wide cost in the steady state. The design of ride-pooling services is optimized as a constrained nonlinear program. A customized solution approach is proposed to decompose and solve the optimization problem. 

The proposed model and solution approach are verified with agent-based simulations, and tested through a series of hypothetical examples and a real-world case study. The results reveal interesting insights into the optimal service design characteristics; e.g., how zone partition size may affect system performance, how idle vehicles and partially-occupied vehicles complement each other in providing services, how the need for vehicles in a zone depends on not only its own travel demand but also those in other zones, and how important it is to properly address demand heterogeneity when designing ride-pooling services. 


The work in this paper can be extended in multiple directions. We limit the number of passengers allowed to share a vehicle, and choose specific operation strategies (e.g., picking up an assigned passenger holds a higher priority than delivering an on-board passenger), to reduce the difficulty of explaining the structure of the proposed multi-zone queuing network model. Yet, the multi-zone modeling framework in this paper can be directly built upon any variations of the original queuing network model (e.g., as in \citet{Daganzo2019Paper}), to address different service types (e.g., taxi and dial-a-ride), different vehicle capacities, and different operation strategies (e.g., whether or not a vehicle is allowed to pick up two passengers consecutively before making any delivery) --- these extensions only vary the “local” network structure within each zone and the transition links across zones.

In addition, although ride-pooling service is advantageous for high-demand areas (such that the matching probability is high), it may be beneficial to use non-shared taxi service where demand is very low. Therefore, an interesting extension to this paper is to allow the TNC platform to decide which zones to provide which type of service, or which combinations of services, and how the TNCs can optimally set the prices for different services (see similar efforts in \citet{zhang2021pool} and \citet{jacob2021ride}). 



In this paper, the system is spatially heterogeneous but temporally homogeneous. In the real world, service demand and supply could often be notably heterogeneous in the time dimension. It would be interesting to explore the system dynamics and possibly adaptive control laws (e.g., for rebalancing) based on either continuous (e.g., the fluid model in \citet{xu2021generalized}) or discrete models (e.g., cell transmission model in \citet{daganzo1994cell} and zone-based system in \citet{lei2019path}). 

Finally, this paper assumes constant vehicle speed in a Manhattan metric space, and the study region is partitioned into square zones for convenience of distance estimation. In future research, it is possible to (i) explicitly capture local traffic congestion by quantifying vehicle travel time as a function of the number of vehicles in each zone (e.g., using a macroscopic fundamental diagram as in \citet{daganzo2008analytical})
; and (ii) try other zone partition shapes (e.g., hexagons) under different distance metrics \citep{xie2015optimal} --- the multi-zone model framework still applies as long as the directional relationships for every pair of zones and their expected distances are computed accordingly. In so doing, we can also incorporate environmental impacts and externalities into the service design model. 

\section*{Acknowledgments}
The authors sincerely thank the editors and three anonymous reviewers for sharing very valuable insights and suggestions. We also acknowledge the help from Yuhui Zhai (UIUC) for developing the agent-based simulation platform. The first author was financially supported in part by the Illinois Department of Transportation via projects ICT-R27-SP45 and ICT-R27-SP53, and by a seed project from the University of Illinois Strategic Research Initiative Growth Fund. 

\newpage
\appendix
\renewcommand*\thetable{\arabic{table}}
\renewcommand*\thefigure{\arabic{figure}}

\section{List of Key Notation}\label{sec:appendix notation}
See Table \ref{tab:notation} for the list of key notation in this paper.
\begin{table}[h!]
\begin{tabular}{p{0.18\linewidth}p{0.8\linewidth}}
\Xhline{4\arrayrulewidth}
Notation & Description  \\ \hline
\\[-2ex]
\multicolumn{2}{l}{\textit{\textbf{Indices and Sets}}} \\
\\[-2ex]
$r$ & Index of travel direction. \\
$\mathcal{K}$ & Set of all zones in the study region, indexed by $i,\ j,\ k$.\\
$G^r_i$ & Set of zones to the $r\in\{\text{NE},\text{N},\text{NW},\text{W},\text{SW},\text{S},\text{SE},\text{E}, +, \times\}$ direction of zone $i$, where $G^+_i = G^{\text{E}}_i \cup G^{\text{N}}_i \cup G^{\text{W}}_i \cup G^{\text{S}}_i$ and $G^{\times}_i = G^{\text{NE}}_i \cup G^{\text{NW}}_i \cup G^{\text{SW}}_i \cup G^{\text{SE}}_i$. \\
$\mathcal{V}_{ij}$ & Set of zones that a vehicle going to zone $j$ can enter upon leaving zone $i$. \\
$U^r$ & Set of two directions adjacent to direction $r$.\\
$\Omega^r_{ii},\ \Omega_{ij}$ & Sets of zones in the destination feasible areas, for a caller going in the $r$ direction within zone $i$, and for a caller going from zone $i$ to zone $j$, respectively.\\
$\tilde{\Omega}_{ij}$ & Subset of zones in $\Omega_{ij}$ that are no closer to zone $i$ than zone $j$.\\
$A^r_i$ & Index of the adjacent zone on the $r$ side of zone $i$. \\
$s = (s_0, s_1, s_2, s_3)$ & Vehicle state vector, $s_0$ -- vehicle current zone, $s_1$ -- assigned caller's origin zone, $s_2$ and $s_3$ -- onboard passengers' destination zones.\\\hline
\\[-2ex]
\textit{\textbf{Parameters}} \\
\\[-2ex]
$\Phi$ [km] & Side length of the square local zone.\\
$v$ [km/hr] & Travel speed of ride-pooling vehicles.\\
$\gamma$ [\$/veh-hr] & TNC vehicle operation cost.\\
$\beta$ [\$/trip-hr] & Passenger value-of-time.\\
$\lambda_{ij}$ [trip/hr] & Generation rate of trips from zone $i$ to zone $j$.\\
$L_{ij}$ [km] & Shortest distance between zone centroids on road network.\\ 
$q$  & Scaling factor, for idle vehicle number in Section \ref{sec:simulation}, or demand in Section \ref{sec:chicago}.\\ \hline
\\[-2ex]
\textit{\textbf{Variables}}\\
\\[-2ex]
$a_s,\ p_{s\_i},\ d_s,\ c_s,\ g_s$ [veh/hr] & Rates for vehicles in state $s$ to receive assignments, pick up passengers destined in zone $i$, deliver passengers, and enter and leave its current zone, respectively.\\
$b_{ij00}$ [veh/hr] & Rate to rebalance idle vehicles from zone $i$ to zone $j$.\\
$n_s$ & Number of vehicles in state $s$.\\
$\tilde{T}_{i0ij},\ \tilde{T}_{i0jk}$ [hr] & Expected remaining trip durations of a state-$(i, 0, i, 0)$ vehicle and a state-$(i, 0, j, 0)$ vehicle after being matched with an inter-zonal caller, respectively.\\
$\delta_{ij\_i'}$ & \% of vehicles going to zone $j$ that enter zone $i'$ upon leaving zone $i$.\\
$N^r_{ii},\ N_{ij}$ & Numbers of vehicles suitable to serve a caller going in the $r$ direction within zone $i$, and to serve a caller going from zone $i$ to zone $j$, respectively.\\
$M_i$,\ $M_b$,\ $M$ [veh-hr/hr] & Vehicle hours in operation per hour, for active vehicles in zone $i$, all rebalanced vehicles in the system, and all vehicles in the system, respectively.\\
$P$ [trip-hr/hr] &  Total trip time spent by all passengers per hour. \\
$Z$ [\$/trip] & System-wide cost per passenger.\\
\Xhline{4\arrayrulewidth}
\end{tabular}
\caption{Key notation in this paper.}
\label{tab:notation}
\end{table}

\newpage
\section{Expected Passenger-Vehicle Matching, and Proofs of \texorpdfstring{Eqn. \eqref{eq:suitable veh}-\eqref{eq:pickup}}{Eqn. (8)-(10)}}\label{sec:appendix matching}

In this appendix, we derive formulas to show how the assignment rates and pickup rates are related to the number of vehicles in each state, as in Eqn. \eqref{eq:suitable veh}--\eqref{eq:pickup}, according to the assignment policy and matching criteria.

Recall that any caller in zone $i$ may be served by an idle vehicle in state $(i, 0, 0, 0)$, or a vehicle in state $(i, 0, i, 0)$ (i.e., with a seeker destined in zone $i$), or a vehicle in state $(i, 0, j, 0)$ (i.e., with a seeker destined in another zone $j$). Based on the matching criteria, all $n_{i000}$ idle vehicles are suitable to serve any caller, but only a portion of vehicles in states $(i, 0, i, 0)$ or $(i, 0, j, 0)$ can be suitable.\footnote{The matching can also be equivalently described from the perspective of a suitable vehicle, which may see multiple subgroups of feasible callers (i.e., callers with different destination zones), each arriving from an independent Poisson process. Each subgroup of feasible callers may be assigned to it according to the matching probability, and hence, the arrivals of all its feasible callers also follow a Poisson Process, whose rate (i.e., vehicle assignment rate) equals the summation of arrival rates across all feasible subgroups with the matching probability as the weight.}

If an intra-zonal caller is going, say, in the NE
direction, the destination feasible area for seekers in state-$(i, 0, i, 0)$ vehicles is indicated by Figure \ref{fig:matching example intra intra}, and the expected size of this feasible area is $2/9$ of the zone.\footnote{Recall that the callers' ODs are homogeneously distributed within the zone; readers can refer to Appendix 7B of \citet{Daganzo2019Book} for the derivation of this quantity.} Therefore, the expected number of suitable state-$(i, 0, i, 0)$ vehicles is $2n_{i0i0}/9$. A state-$(i, 0, j, 0)$ vehicle is suitable if its seeker's destination is in a zone of $\Omega^{\text{NE}}_{ii}$, as shown in Figure \ref{fig:matching example inter intra}; therefore, the number of suitable state-$(i, 0, j, 0)$ vehicles is $\sum_{j\in \Omega^{\text{NE}}_{ii}}n_{i0j0}$. Then, the total number of suitable vehicles for an intra-zonal caller going in the NE direction, $N^{\text{NE}}_{ii}$, is as follows,

\[N_{ii}^r = n_{i000} + \frac{2n_{i0i0}}{9} + \sum_{j\in \Omega^r_{ii}}n_{i0j0},\quad \forall i \in \mathcal{K},\ r\in\{\text{NE}\}.\]
By symmetry, the above formula holds for any travel direction $r\in \{\text{NE}, \text{NW}, \text{SW}, \text{SE}\}$, which yields Eqn. \eqref{eq:suitable veh}.

When the caller is inter-zonal, e.g., going from zone $i$ to zone $j$, we can obtain the number of suitable vehicles in a similar way. Among the seekers in state-$(i, 0, i, 0)$ vehicles, we expect $1/2$ of them to have their destinations aligned with that of the inter-zonal caller if $j\in G^{\text{S}}_i$ (e.g., see the blue shaded area in Figure \ref{fig:matching example intra inter}), or $1/4$ if $j\in G^{\text{NE}}_i$ (e.g., see the red shaded area in Figure \ref{fig:matching example intra inter}). Among the state-$(i, 0, k, 0)$ vehicles, they are suitable if their seekers' destination zones are in $\Omega_{ij}$ (e.g., see the shaded area in Figure \ref{fig:matching example inter inter}). Hence, the number of suitable vehicles for an inter-zonal caller going
from zone $i$ to zone $j$, $N_{ij}$, is as follows,

\begin{equation*}\label{eq:appendix suitable veh}
\begin{aligned}
    &N_{ij} = n_{i000} + \frac{n_{i0i0}}{4} + \sum_{k\in \Omega_{ij}}n_{i0k0},\quad &&\forall i \in \mathcal{K}, \ j \in G_i^r, \ r\in \{\text{NE}\}\\
    &N_{ij} = n_{i000} + \frac{n_{i0i0}}{2} + \sum_{k\in \Omega_{ij}}n_{i0k0},\quad &&\forall i \in \mathcal{K}, \ j \in G_i^r, \ r\in \{\text{S}\}.\\
\end{aligned}
\end{equation*}
Again, by directional symmetry, we can see that Eqn. \eqref{eq:suitable veh_ij} holds. 

Now we determine the caller pickup rates in Eqn. \eqref{eq:pickup}. Since the TNC always assigns the nearest suitable vehicle to a caller, and the vehicles in each state are approximately uniformly distributed in a zone, the probability for a caller to be assigned to a vehicle in state $s$ is proportional to the relative density of that type of vehicles.
Also, the rate for vehicles in state $s$ to be assigned to a type of callers should be equal to the rate for vehicles in state $s$ to pick up that type of callers. Therefore, we can determine the pickup rates based on all types of suitable vehicles, $n_{i000},\ n_{i0i0},\ \{n_{i0j0}\}_{j\neq i}$.

For example, in zone $i$, symmetry indicates that intra-zonal callers going in the $r\in$\{NE, NW, SW, SE\} direction are generated at rate $\lambda_{ii}/4$, and a fraction of $n_{i000}/N_{ii}^r$ of them are picked up by an idle vehicle. Therefore, the rate for idle vehicles to pick up intra-zonal callers, $p_{ii00\_i}$, can be obtained by summing the product of these two quantities, $(\lambda_{ii}/4)\cdot (n_{i000}/N_{ii}^r)$, across all four directions; i.e., 

\begin{equation*}
p_{ii00\_i} = \sum_{r\in \{\text{NE}, \text{NW}, \text{SW}, \text{SE}\}}\left(\frac{\lambda_{ii}}{4}\cdot\frac{n_{i000}}{N^r_{ii}}\right).
\end{equation*}
This gives Eqn. \eqref{eq:pickup 1}. The rate for state-$(i, 0, i, 0)$ vehicles to pick up intra-zonal callers, $p_{iii0\_i}$, can be derived similarly, except that the number of suitable state-$(i, 0, i, 0)$ vehicles is $2n_{i0i0}/9$ instead of $n_{i000}$; i.e., 

\begin{equation*}
p_{iii0\_i} =\sum_{r\in \{\text{NE}, \text{NW}, \text{SW}, \text{SE}\}}\left(\frac{\lambda_{ii}}{4}\cdot\frac{2n_{i0i0}}{9N^r_{ii}}\right). 
\end{equation*}
This gives Eqn. \eqref{eq:pickup 2}. Per the matching criterion in \ref{case:3}, an intra-zonal caller going in the $r\in$\{NE, NW, SW, SE\} direction can also be picked up by a state-$(i, 0, j, 0)$ vehicle, whose seeker is either (i) going to a zone $j$ in the exact same direction $r$ (i.e., $j \in G^r_i$); or (ii) going to a zone $j$ in an adjacent direction (i.e., $j \in G^{r'}_i$, $\forall r'\in U^r$). Conversely, from the vehicle's point of view, a state-$(i, 0, j, 0)$ vehicle with a seeker going to a zone on the $r\in$\{NE, NW, SW, SE\} direction can only serve an intra-zonal passenger of the exact same direction; the rate for pickup is

\begin{equation*}
p_{iij0\_i} = \frac{\lambda_{ii}}{4}\cdot\frac{n_{i0j0}}{N^r_{ii}},\quad \forall i\in \mathcal{K}, \ j \in G^r_i,\ r \in \{\text{NE}, \text{NW}, \text{SW}, \text{SE}\}.
\end{equation*}
This gives Eqn. \eqref{eq:pickup 3}. Yet, a state-$(i, 0, j, 0)$ vehicle with a seeker going in a direction $r'\in$\{E, N, W, S\} can serve intra-zonal passengers that are going in the two adjacent directions of $r'$, as given by $U^{r'}$; i.e.,

\begin{equation*}
p_{iij0\_i} =\sum_{r\in U^{r'}}\left(\frac{\lambda_{ii}}{4}\cdot\frac{n_{i0j0}}{N^{r}_{ii}}\right),\quad \forall i\in \mathcal{K}, \ j \in G^{r'}_i, \ r' \in \{\text{N, W, S, E}\}.
\end{equation*}
This gives Eqn. \eqref{eq:pickup 4}. 

When a caller is inter-zonal with its destination in zone $j$, the analysis is similar. We only need to note that the trip rate is now $\lambda_{ij}$, and the total number of suitable vehicles is now $N_{ij}$. The reader is encouraged to verify that the pickup rates for inter-zonal callers by different types of vehicles, as counterparts of those for intra-zonal callers, are given by Eqn. \eqref{eq:pickup 5}-\eqref{eq:pickup 8}.


\section{Expected State Duration, and Proofs of \texorpdfstring{Eqn. \eqref{eq:delivery and outgoing}-\eqref{eq:fleet}}{Eqn. (11)-(12)}}\label{sec:appendix duration}
This appendix derives the expected vehicle travel distance in each state, first among the ``assigned" and ``rebalanced" states, and then the ``delivery" states. These distances will be used to quantify the rates at which vehicles finish passenger delivery trip within zone $i$, as in Eqn. \eqref{eq:delivery and outgoing}, and to obtain the count of vehicles in each state, as in Eqn. \eqref{eq:fleet}. 

\subsection{Assigned and rebalanced state duration}
In zone $i$, a vehicle is en-route to pick up an assigned passenger if it is in state $(i, i, 0, 0)$, $(i, i, i, 0)$, and $(i, i, j, 0)$, $\forall j\in\mathcal{K}\setminus\{i\}$. The expected pickup travel distance under the nearest vehicle assignment rule, per discussion in \citet{Daganzo2019Book}, is simply inversely proportional to the number of suitable vehicles; i.e., $0.63\Phi/\sqrt{N_{ii}^r}$ if the passenger is intra-zonal and going in the $r\in\{\text{NE}, \text{NW}, \text{SW}, \text{SE}\}$ direction, or $0.63\Phi/\sqrt{N_{ij}}$ if the passenger is inter-zonal and going to zone $j\in\mathcal{K}\setminus\{i\}$. Dividing these distances by the vehicle speed $v$ yields the expected trip durations. Therefore, by Little's Law \citep{little1961proof}, the vehicle count in state $(i, i, 0, 0)$ can be obtained by summing the product of each expected trip duration, $\frac{0.63\Phi}{v\sqrt{N_{ii}^r}}$ or $\frac{0.63\Phi}{v\sqrt{N_{ij}}}$, and the corresponding pickup rate, $\frac{\lambda_{ii}}{4}\frac{n_{i000}}{N_{ii}^r}$ from Eqn. \eqref{eq:pickup 1} or $\lambda_{ij}\frac{n_{i000}}{N_{ij}}$ from Eqn.\eqref{eq:pickup 5}, across all intra-zonal travel directions $r\in\{\text{NE}, \text{NW}, \text{SW}, \text{SE}\}$ and all inter-zonal passenger destination zones $j\in\mathcal{K}\setminus\{i\}$. The result leads to Eqn.\eqref{eq:fleet 1}.

The vehicle counts in states $(i, i, i, 0)$ and $(i, i, j, 0)$ shall be similar, except that (i) the rates for state-$(i, i, i, 0)$ vehicles and state-$(i, i, j, 0)$ vehicles to pick up an intra-zonal passenger are given by Eqn. \eqref{eq:pickup 2} and Eqn. \eqref{eq:pickup 3}-\eqref{eq:pickup 4}, respectively; and (ii) the rates for state-$(i, i, i, 0)$ vehicles and state-$(i, i, j, 0)$ vehicles to pick up an inter-zonal passenger are given by Eqn. \eqref{eq:pickup 6}-\eqref{eq:pickup 7} and Eqn. \eqref{eq:pickup 8}, respectively. The readers can then verify from Little's formula that the vehicle counts in these two states are given by Eqn. \eqref{eq:fleet 2}-\eqref{eq:fleet 4}.

Finally, it shall be clear that the expected travel distance for a vehicle to be rebalanced from a random point in zone $i$ to a random point in zone $j$ shall be equal to the distance between the zone centroids, $L_{ij}$, if these two zones are not horizontally or vertically aligned (i.e., $j\in G^{\times}_i$). Otherwise, if $j\in G^+_i$, we need to add an additional ``lateral" travel distance that is equal to one third of the zone width \citep{Daganzo2019Book}; i.e., the expected distance becomes $L_{ij} + \frac{\Phi}{3}$. Little's formula then leads to Eqn. \eqref{eq:fleet 8}-\eqref{eq:fleet 9}.

\subsection{Delivery state duration}
In zone $i$, a vehicle is en-route to drop off an onboard passenger if it is in state $(i, 0, i, 0)$, $(i, 0, i, i)$, $(i, 0, j, 0)$, $(i, 0, i, j)$ or $(i, 0, j, k),\ \forall j\in\mathcal{K}\setminus\{i\},\ k\in\tilde{\Omega}_{ij}$. The expected travel distances for vehicles in these states are more complicated, because these vehicles may have one or two onboard passenger(s) who could be intra-zonal or inter-zonal, and these vehicles may transition into and out of the current state due to multiple types of events: picking up or delivering a passenger, receiving a new assignment, and crossing a zone boundary. 

\subsubsection{State \texorpdfstring{$(i, 0, i, i)$}{(\textit{i}, 0, \textit{i}, \textit{i})}}\label{sec:appendix i0ii}
We first look at the delivery trip of a vehicle with two onboard passengers destined zone $i$, i.e., in state $(i, 0, i, i)$. Let P1 and P2 represent the first onboard seeker and the later onboard caller, respectively, and their ODs are denoted by O1, D1, and O2, D2, accordingly, for the convenience of explanation. As shown in Figure \ref{fig:queuing model}, the vehicle could have transitioned from either state $(i, i, i, 0)$, at rate $p_{iii0\_i}$, upon picking up P2 at O2; or state $(i', 0, i, i)$ (where $i'\in \mathcal{K}\setminus\{i\}$), at rate $g_{i'0ii}$, upon entering zone $i$. 

\begin{figure}[t]
     \centering
     \begin{subfigure}[b]{0.6\textwidth}
         \centering
         \includegraphics[width=\textwidth]{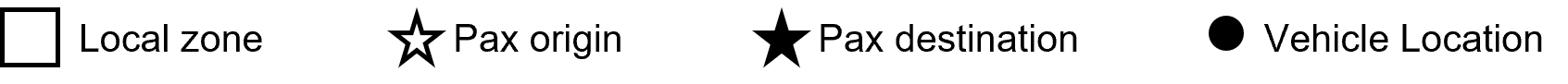}
     \end{subfigure}\\
     \begin{subfigure}[t]{0.4\textwidth}
         \centering
         \includegraphics[width=0.8\textwidth]{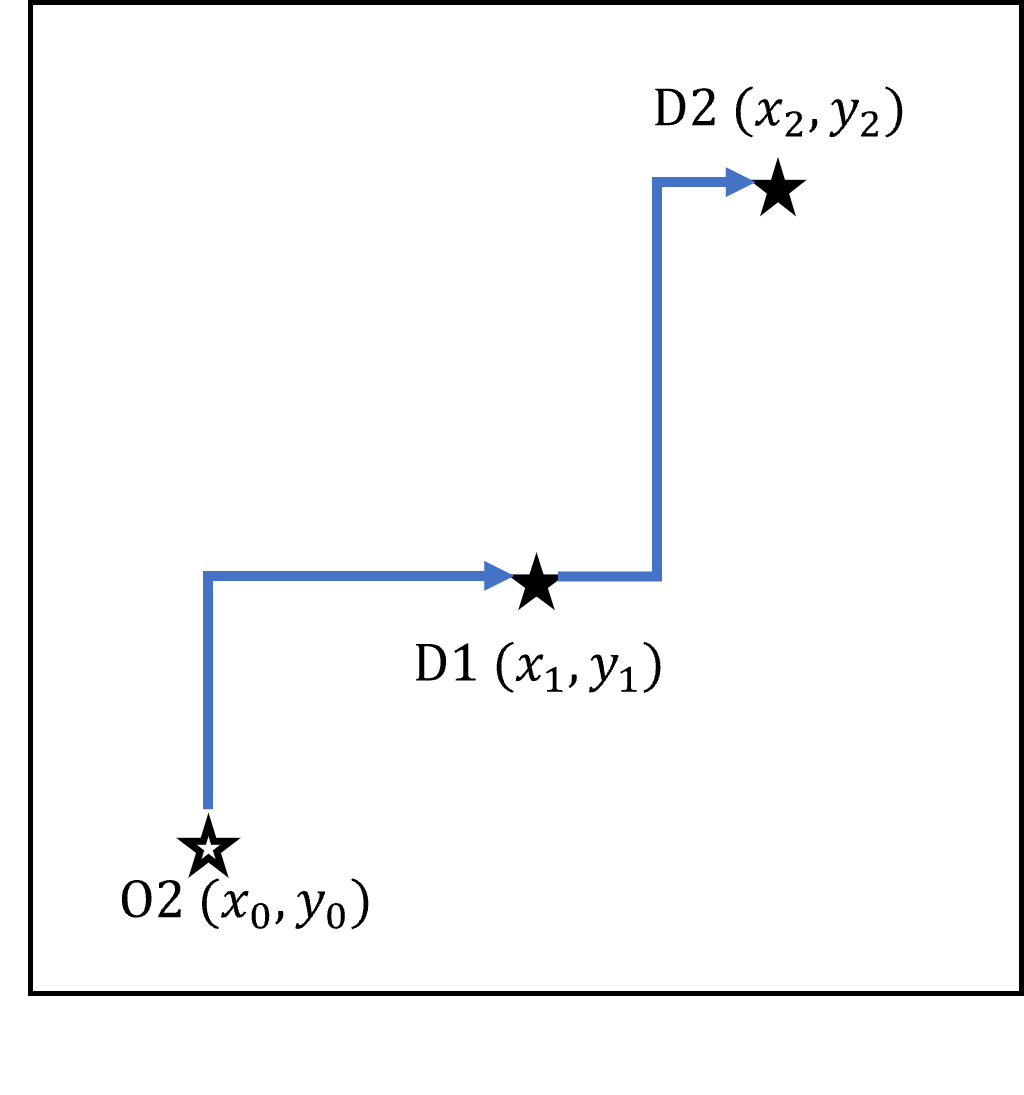}
         \caption{At least one trip originates in zone $i$.}
         \label{fig:distance 1}
     \end{subfigure}
     \begin{subfigure}[t]{0.4\textwidth}
         \centering
         \includegraphics[width=0.8\textwidth]{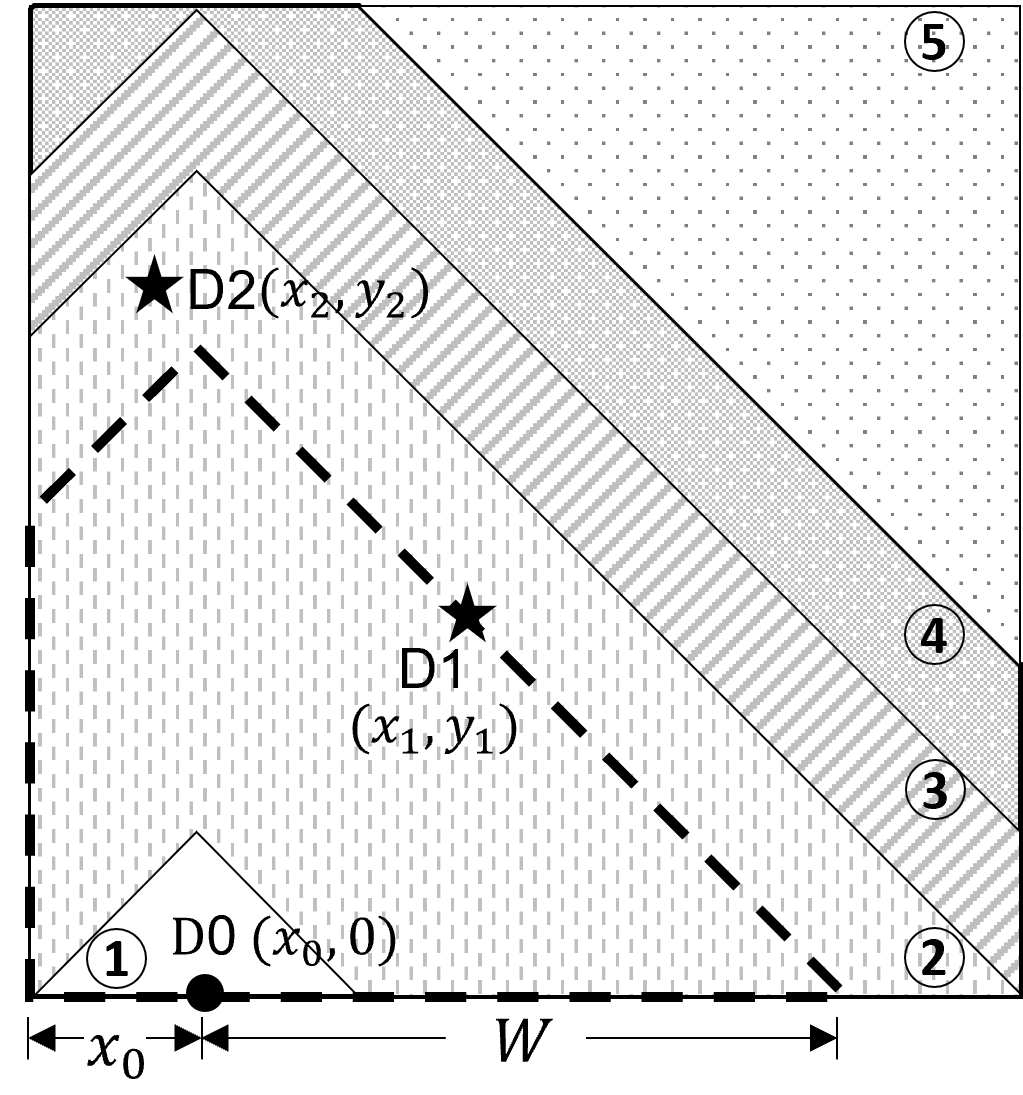}
         \caption{Both trips originate outside of zone $i$.}
         \label{fig:distance 2}
     \end{subfigure}
    \caption{Illustration of a state-$(i, 0, i, i)$ vehicle delivering the passenger with a closer destination, transitioning to state $(i, 0, i, 0)$ upon delivery, and delivering the other passenger.}
    \label{fig:distance}
\end{figure}

In the former case, the delivery trip starts at O2, and ends at the closer one of D1 and D2, where D1 and D2 should not cause detours based on matching criterion \ref{case:1}. Let $(x_0, y_0)$, $(x_1, y_1)$ and $(x_2, y_2)$ represent the coordinates of O2, D1, and D2, respectively. Note that $x_0$, $x_1$, $x_2$, $y_0$, $y_1$, and $y_2$ are uniformly distributed in ranges $[0, \Phi]$. Figure \ref{fig:distance 1} shows an example where D1 and D2 are aligned toward the NE direction of O2, and D1 is closer to O2 (i.e., $x_0\leq x_1\leq x_2$ and $y_0\leq y_1\leq y_2$). The expected distance between O2 and D1, conditioning on $x_0\leq x_1\leq x_2$ and $y_0\leq y_1\leq y_2$, is clearly,
\begin{equation}\label{eq:i0ii delivery trip}
\small\frac{\int_0^\Phi\int^\Phi_{x_0}\int^{\Phi}_{x_1}\int_0^\Phi\int^\Phi_{y_0}\int^{\Phi}_{y_1}\frac{x_1 - x_0 + y_1 - y_0}{\Phi^6}dy_2dy_1dy_0dx_2dx_1dx_0}{\Pr(x_0\leq x_1\leq x_2, y_0\leq y_1\leq y_2)} = \frac{\Phi}{2}.
\end{equation}
where $\Pr(x_0\leq x_1\leq x_2, y_0\leq y_1\leq y_2) = \frac{1}{36}$. Similarly, when D2 is closer to O2, or when the vehicle is traveling in other directions, the expected distance is also $\frac{\Phi}{2}$. 

In the latter case, the delivery trip starts when the vehicle enters zone $i$. The entry point into zone $i$, denoted by D0, could either be a random point on the zone border if O2 is in zone $i'\in G^{+}_i$, or at a vertex of zone $i$ if O2 is in zone $i'\in G^{\times}_i$. By matching criterion \ref{case:2}, detour within zone $i$ is allowed when delivering two such inter-zonal passengers, i.e., D1 and D2 are uniformly distributed within zone $i$. When D0 is a random location on the zone boundary, as shown in Figure \ref{fig:distance 2}, the length of the delivery trip is $\min\{|x_1-x_0|+|y_1 - 0|, |x_2-x_0|+|y_2 - 0|\}$. Let $W$ denote the delivery trip length. Since it depends on the $x$-coordinate of D0, $x_0$, we can obtain its complementary cumulative distribution function conditional on $x_0$, denoted by $\bar{F}_{W|x_0}(w|x)$, which equals the conditional probability that both D1 and D2 are outside of the dashed-line polygon in Figure \ref{fig:distance}\subref{fig:distance 2}). When $x_0\leq \frac{\Phi}{2}$, $\bar{F}_{W|x_0}(w|x)$ is given as follows,

\begin{equation}\label{eq:complementary CDF}
    \bar{F}_{W|x_0}(w|x)= 
    \begin{cases}
        \left(1 - \frac{w^2}{\Phi^2}\right)^2, & \text{if } w\leq x\\
        \left(1 - \frac{w^2+2xw-x^2}{2\Phi^2}\right)^2, & \text{if } x\leq w\leq \Phi-x\\
        \left(1 - \frac{2\Phi w+ 2\Phi x-2x^2-\Phi^2}{2\Phi^2}\right)^2, & \text{if } \Phi-x\leq w\leq \Phi\\
        \left(1 - \frac{6\Phi w+ 2\Phi x-2x^2-3\Phi^2-2w^2}{2\Phi^2}\right)^2, & \text{if } \Phi\leq w\leq \Phi+x\\
        \frac{(2\Phi-w-x)^4}{4\Phi^4}, & \text{if } \Phi+x\leq w\leq 2\Phi-x.
    \end{cases}
\end{equation}
These five cases correspond to areas 1-5 in Figure \ref{fig:distance}\subref{fig:distance 2}, respectively. We first calculate $\mathbb{E}[W|x_0\leq\frac{\Phi}{2}]=\frac{1}{\Phi/2}\int^{\Phi/2}_0\int^{2\Phi-x}_0\bar{F}_{W|x_0}(w|x)dwdx = \frac{5\Phi}{8}$. By symmetry, we know that $\mathbb{E}[W] = \frac{5\Phi}{8}$. It is noted that when D0 is at a zone's corner, the complementary CDF of the first delivery trip length can be obtained by setting $x=0$ in Eqn. \eqref{eq:complementary CDF}, and the simple integral shows that the expected distance of the first delivery trip in this case is $\frac{23\Phi}{30}$. These results are summarized in Table \ref{tab:travel dist}.

Upon completing the first delivery trip, the vehicle transitions into state $(i, 0, i, 0)$ and starts to deliver the remaining passenger (i.e., going between D1 to D2). The expected distance between D1 and D2 is also $\frac{\Phi}{2}$ (similar to Eqn. \eqref{eq:i0ii delivery trip}) in the former case, and simply $\frac{2\Phi}{3}$ in the latter case \citep{Daganzo2019Book}.

\subsubsection{State \texorpdfstring{$(i, 0, i, 0)$}{(i, 0, i, 0)}}
A state-$(i, 0, i, 0)$ vehicle 
is en-route to deliver its only onboard seeker to the destination within the current zone. As shown in Figure \ref{fig:queuing model}, this vehicle could have transitioned into the current state from (i) state $(i', 0, i, 0)$ for $i'\in \{A^{\text{E}}_i, A^{\text{N}}_i, A^{\text{W}}_i, A^{\text{S}}_i\}$, at inflow rate $c_{i0i0}$, upon entering zone $i$ at a random boundary point; (ii) state $(i, i, 0, 0)$, at inflow rate $p_{ii00\_i}$, upon picking up an assigned caller at a random interior point; or (iii) state $(i, 0, i, i)$, at inflow rate $d_{i0ii}$, upon delivering one of its onboard passengers to the destination. In case (i), the trip to deliver the seeker is between the random entry point on a zone boundary and the seeker's destination (i.e., a random interior point), with expected distance $\frac{5\Phi}{6}$ = $\frac{\Phi}{2}$ (longitudinal) + $\frac{\Phi}{3}$ (lateral). In case (ii), the delivery trip is between an intra-zonal passenger's OD, which are two random points in a zone, with expected distance $\frac{2\Phi}{3}$. In case (iii), there are two possibilities, as discussed at the end of \ref{sec:appendix i0ii}, with expected trip length: (iii-a) $\frac{\Phi}{2}$ if this vehicle enters state $(i, 0, i, i)$ upon a recent pickup (at inflow rate $p_{iii0\_i}$); or (iii-b) $\frac{2\Phi}{3}$ if this vehicle enters state $(i, 0, i, i)$ by crossing a zone boundary (at inflow rate $c_{i0ii}$). 
These expected trip lengths are summarized in Table \ref{tab:travel dist}.

A state-$(i, 0, i, 0)$ vehicle could end the current state when it completes the delivery of the seeker, represented by rate $d_{i0i0}$, or when a new caller is assigned to this vehicle before it finishes the delivery, which occurs according to a Poisson Process with rate $a_{i0i0}/n_{i0i0}$ (i.e., the total assignment rate $a_{i0i0}$is shared by $n_{i0i0}$ uniformly distributed vehicles in this state). Hence, the inter-arrival time of assignments follows an exponential distribution with rate  $\frac{a_{i0i0}}{n_{i0i0}}$ for every $(i, 0, i, 0)$ vehicle. As a result, the probability for the vehicle to successfully complete the delivery (i.e., the assignment arrival time is greater than the trip duration) is $\exp(-\frac{a_{i0i0}}{n_{i0i0}}\cdot\frac{5\Phi}{6v})$, $\exp(-\frac{a_{i0i0}}{n_{i0i0}}\cdot\frac{2\Phi}{3v})$,
$\exp(-\frac{a_{i0i0}}{n_{i0i0}}\cdot\frac{\Phi}{2v})$, and $\exp(-\frac{a_{i0i0}}{n_{i0i0}}\cdot\frac{2\Phi}{3v})$ for cases (i), (ii), (iii-a) and (iii-b), respectively.\footnote{Here, we approximate the probability for a $(i, 0, i, 0)$ vehicle to receive no assignments during the delivery trip by using the average trip distance, for simplicity. By Jensen's inequality, this approximation underestimates the actual probability. However, the error should be small, especially when the zone size is small.
} Therefore, the successful delivery rate $d_{i0i0}$ shall be the weighted sum of all inflow rates $c_{i0i0}$, $p_{ii00\_i}$, $c_{i0ii}$, and $p_{iii0\_i}$, using their corresponding success probabilities as weights. This gives $d_{i0i0}$ as in Eqn. \eqref{eq:delivery only}.

\subsubsection{State \texorpdfstring{$(i, 0, j, 0)$}{(\textit{i}, 0, \textit{j}, 0)}}
The analysis of a state-$(i, 0, j, 0)$ vehicle is almost the same as that for a state-$(i, 0, i, 0)$ vehicle, except that its onboard seeker now is destined in zone $j$. Figure \ref{fig:queuing model} shows that this vehicle could have transitioned into the current state from (i) state $(i', 0, j, 0)$ with $i'\in \{A^{\text{E}}_i, A^{\text{N}}_i, A^{\text{W}}_i, A^{\text{S}}_i\}$, at rate $c_{i0j0}$, upon crossing a zone boundary; or (ii) state $(i, i, 0, 0)$, at rate $p_{ii00\_j}$, upon picking up a recent inter-zonal caller from a random location; 
or (iii) state $(i, 0, i, j)$, at rate $d_{i0ij}$, upon delivering a passenger at a random location inside zone $i$. In case (iii), the two onboard passengers of a state-$(i, 0, i, j)$ vehicle could: (iii-a) both board the vehicle outside of zone $i$ and enter zone $i$ at rate $c_{i0ij}$, or (iii-b) the last onboard one is picked up within zone $i$, 
at rate $p_{iii0\_j}$ or $p_{iij0\_i}$. This vehicle state could end when this vehicle crosses the zone border, at rate $g_{i0j0}$, or when a new assignment arrives before it crosses the zone border, at rate $a_{i0j0}$. 
If no new assignment is ever received before the vehicle exit the zone, the corresponding expected travel distances for cases (i) - (iii-a) are approximately $\Phi$, $\frac{\Phi}{2}$ and $\frac{\Phi}{2}$ (based on the routing rules), respectively, from simple geometry. 
In case (iii-b), the distance 
approximately follows a triangular distribution with mean $\frac{\Phi}{3}$, if we simply assume that 
the recent delivery trip goes between two uniformly distributed points in zone $i$.\footnote{This approximation ignores a few other possibilities, as discussed in state-$(i, 0, i, 0)$ vehicles' delivery trips.}
These results are summarized in Table \ref{tab:travel dist}. Note that new caller assignments arrive with rate $a_{i0j0}/n_{i0j0}$ to this particular vehicle, then using a similar logic to that for calculating $d_{i0i0}$ above, we can obtain the formula for $g_{i0j0}$ as in Eqn. \eqref{eq:outgoing only}.

\subsubsection{State \texorpdfstring{$(i, 0, i, j)$}{(\textit{i}, 0, \textit{i}, \textit{j})}}

\begin{figure}[t]
     \centering
     \begin{subfigure}[b]{\textwidth}
         \centering
         \includegraphics[width=\textwidth]{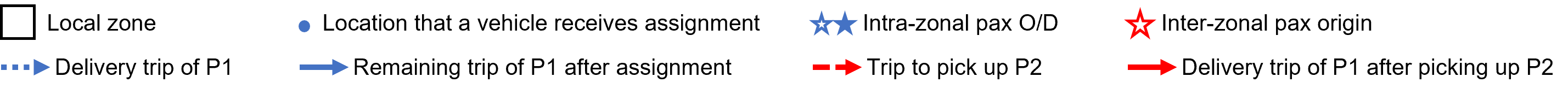}
     \end{subfigure}
     \begin{subfigure}[b]{0.242\textwidth}
         \centering
         \includegraphics[width=\textwidth]{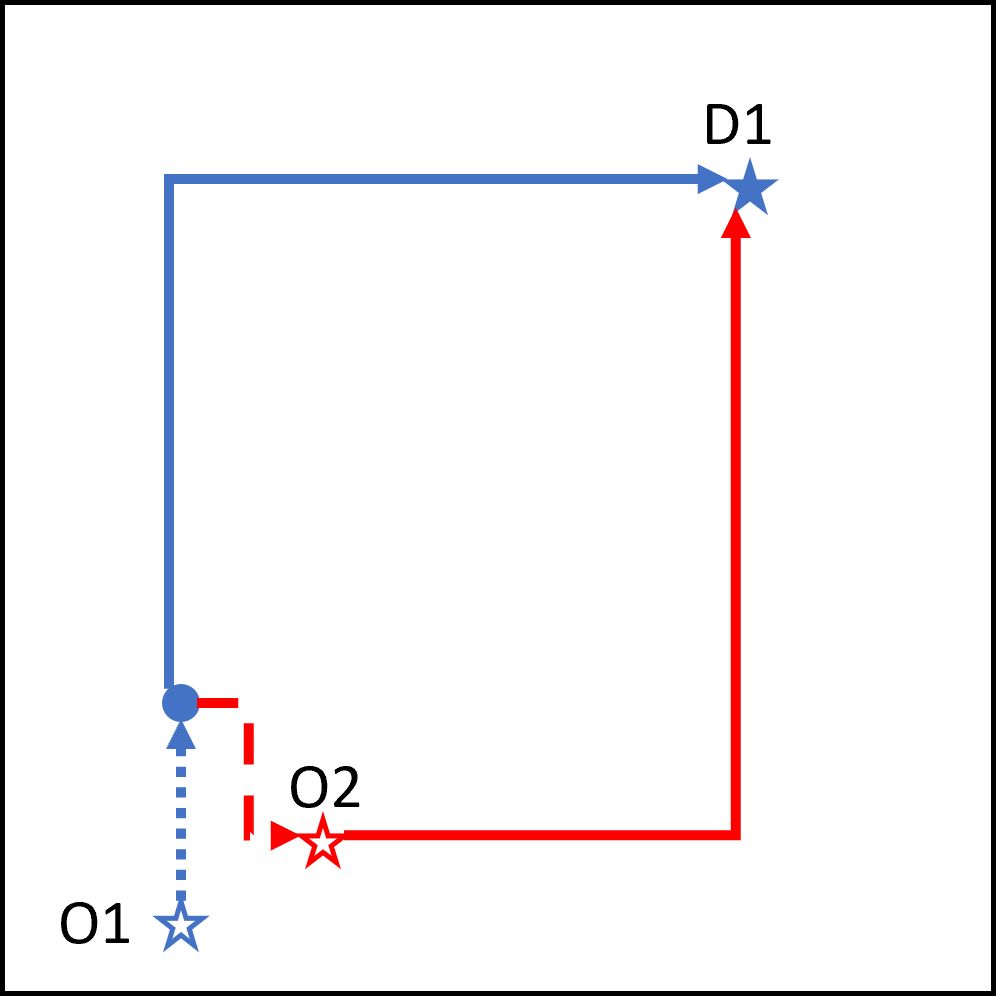}
         \caption{\tiny An $(i, 0, i, 0)$ veh matched w. an inter-zonal caller before turning.}
         \label{fig:i0ij 1}
     \end{subfigure}
     \begin{subfigure}[b]{0.242\textwidth}
         \centering
         \includegraphics[width=\textwidth]{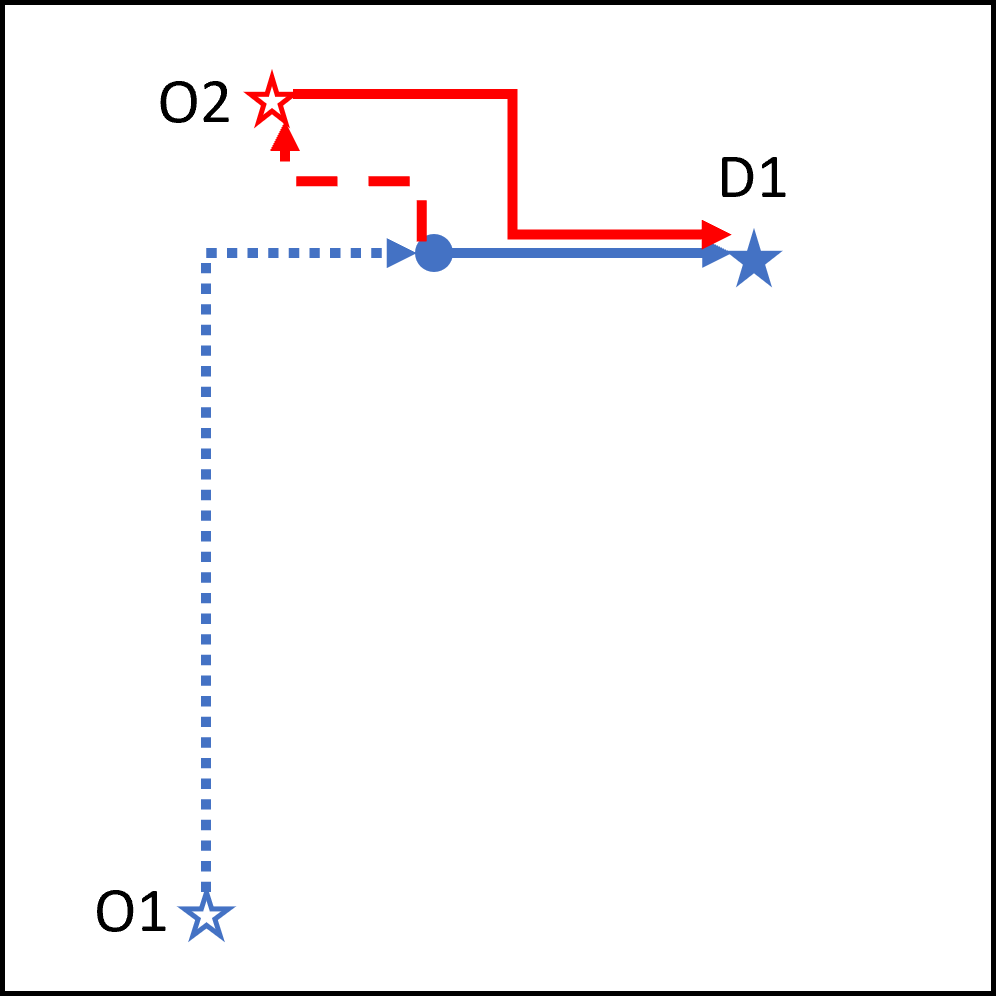}
         \caption{\tiny An $(i, 0, i, 0)$ veh matched w. an inter-zonal caller after turning.}
         \label{fig:i0ij 2}
     \end{subfigure}
     \begin{subfigure}[b]{0.242\textwidth}
         \centering
         \includegraphics[width=\textwidth]{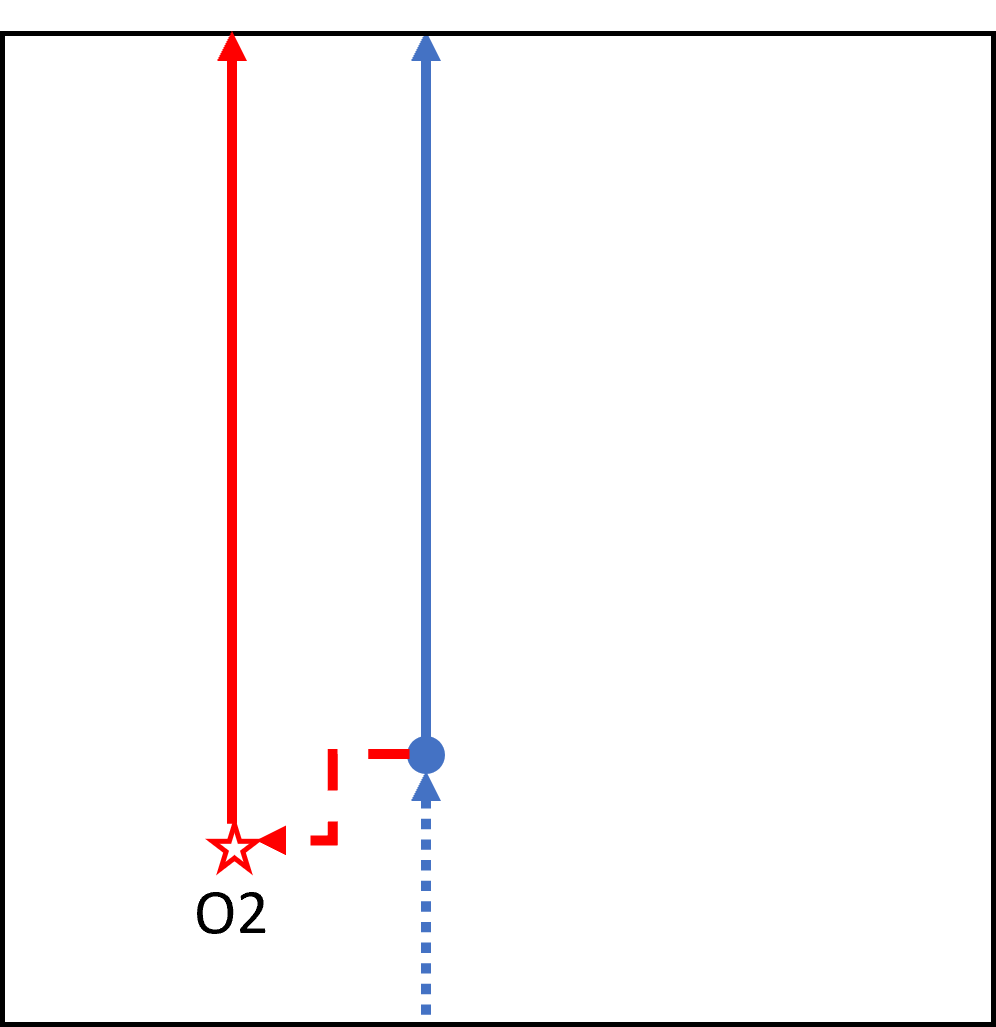}
         \caption{\tiny An $(i, 0, j, 0)$ veh matched w. an inter-zonal caller w/o turning.}
         \label{fig:i0jk 1}
     \end{subfigure}
     \begin{subfigure}[b]{0.242\textwidth}
         \centering
         \includegraphics[width=\textwidth]{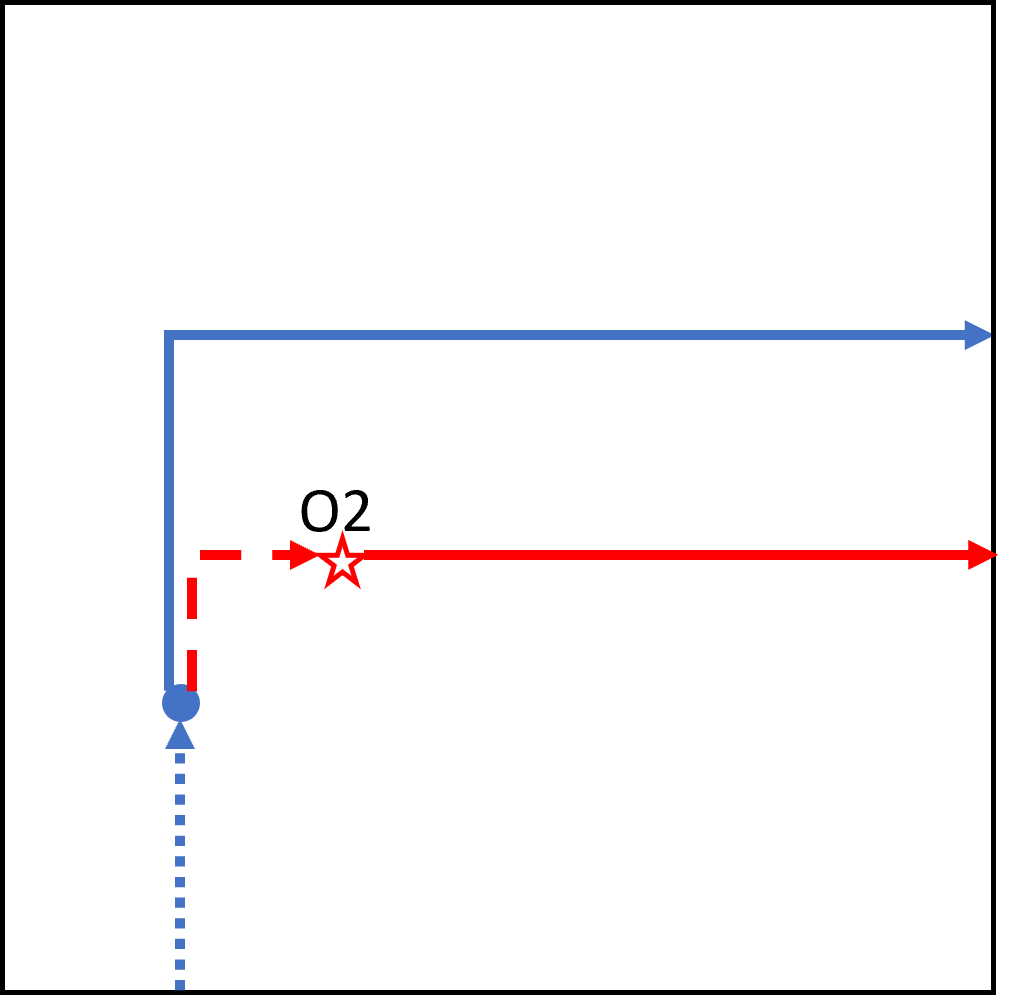}
         \caption{\tiny An $(i, 0, j, 0)$ veh matched w. an inter-zonal caller before turning.}
         \label{fig:i0jk 2}
     \end{subfigure}
    \caption{Example of trips for a vehicle with a seeker (P1), i.e., originally in state $(i, 0, i, 0)$ or $(i, 0, j, 0)$ to deliver P1, to be assigned an inter-zonal caller (P2), to pick up the caller, and to continue to deliver P1 after pickup.}
    \label{fig:i0ij}
\end{figure}

A vehicle in state $(i, 0, i, j)$ could have transitioned into the current state from (i) state $(i, i, i, 0)$ at rate $p_{iii0\_j}$, upon picking up an inter-zonal caller P2; or (ii) state $(i, i, j, 0)$ at rate $p_{iij0\_i}$, upon picking up an intra-zonal caller P2; or (iii) state $(i', 0, i, j)$ at rate $c_{i'0ij}$, upon entering zone $i$. 

In case (i), before being matching with P2 (destined in zone $j$) and transitioning to state $(i, i, i, 0)$, the 
vehicle is originally delivering P1 destined in zone $i$ in state $(i, 0, i, 0)$. After picking up P2, the next delivery leg goes from O2 to D1. 
Figure \ref{fig:i0ij 1} illustrates an example in which a state-$(i, 0, i, 0)$ vehicle originally delivers P1, receives an assignment of P2 while passing the solid blue dot and transitions into state $(i, i, i, 0)$, changes route to pick up P2 along the dashed red arrow, transitions into state $(i, 0, i, j)$ upon pickup, and then starts to deliver P1 from O2 along the solid red arrow. 
Due to the nearest assignment rule, O2 tends to be centered around the blue dot, instead of being uniformly distributed within zone $i$. Thus, the expected distance between O2 and D1 equals that between the blue dot and D1 (i.e., remaining delivery distance of the state-$(i, 0, i, 0)$ vehicle). 
Let $T_{i0i0}$ denote the arrival time of a new assignment to a state-$(i, 0, i, 0)$ vehicle, which follows an exponential distribution with rate $\frac{a_{i0i0}}{n_{i0i0}}$. Since the expected delivery trip length of this vehicle is $\frac{2\Phi}{3}$, the duration of the remaining delivery trip after receiving the assignment is $\frac{2\Phi}{3v} - T_{i0i0}$, and its expectation conditioning on the fact that the vehicle receives an assignment (i.e., this vehicle does transition to $(i, 0, i, j)$ state) can be calculated as $\mathbb{E}\left[\frac{2\Phi}{3v} - T_{i0i0}|T_{i0i0}\leq \frac{2\Phi}{3v}\right] = \left(1 - e^{-\frac{a_{i0i0}}{n_{i0i0}}\frac{2\Phi}{3v}}\right)^{-1}\int_0^{\frac{2\Phi}{3v}}\left(\frac{2\Phi}{3v} - t\right)\cdot\frac{a_{i0i0}}{n_{i0i0}}e^{-\frac{a_{i0i0}}{n_{i0i0}}t}dt = \frac{2\Phi}{3v}\left(1 - e^{-\frac{a_{i0i0}}{n_{i0i0}}\frac{2\Phi}{3v}}\right)^{-1} - \frac{n_{i0i0}}{a_{i0i0}}$. 
Similarly, we can also obtain the remaining delivery trip duration of a state-$(i, 0, i, 0)$ vehicle with expected trip length $\frac{5\Phi}{6}$ (if it transitions into state $(i, 0, i, 0)$ at rate $c_{i0i0}$) or $\frac{\Phi}{2}$ (if it comes from state $(i, 0, i, i)$ at rate $p_{iii0\_i}$), as summarized in Table \ref{tab:travel dist}. Then, by averaging these expectations using their occurrence rates $p_{ii00\_i} + c_{i0ii}$, $c_{i0i0}$, and $p_{iii0\_i}$ as weights, we can obtain the expected trip duration of a state-$(i, 0, i, j)$ vehicle for case (i), denoted by $\tilde{T}_{i0ij}$, as follows,
\begin{equation}\label{eq:i0ij remaining}
    \tilde{T}_{i0ij} = \frac{\Phi/v}{c_{i0i0} + p_{ii00\_i} + d_{i0ii}}\left[
    \frac{5 c_{i0i0}/6}{1 - e^{-\frac{a_{i0i0}}{n_{i0i0}}\frac{5\Phi}{6v}}} + \frac{2 (c_{i0ii} + p_{ii00\_i})/3}{1 - e^{-\frac{a_{i0i0}}{n_{i0i0}}\frac{2\Phi}{3v}}} + \frac{p_{iii0\_i}/2}{1 - e^{-\frac{a_{i0i0}}{n_{i0i0}}\frac{\Phi}{2v}}}\right] - \frac{n_{i0i0}}{a_{i0i0}}.
\end{equation}

In case (ii), this state-$(i, 0, i, j)$ vehicle's delivery trip goes between P2's OD, and the expected length is simply $\frac{2\Phi}{3}$. In case (iii), the delivery trip starts when the vehicle enters zone $i$ and ends at the passenger's destination within zone $i$ (i.e., a random interior point), and the expected trip length is $\frac{5\Phi}{6}$ if the vehicle comes from zone $i'\in G^{+}_i\cap\tilde{\Omega}_{ji}$ and enters zone $i$ at a random point on the border; or $\Phi$ if $i'\in G^{\times}_i\cap\tilde{\Omega}_{ji}$ and the vehicle enters zone $i$ at a zone corner.

As a final remark, we note that a state-$(i, 0, i, j)$ vehicle's delivery trip could be longer than a state-$(i, 0, i, 0)$ vehicle's remaining trip (i.e., Eqn. \eqref{eq:i0ij remaining} underestimates $\tilde{T}_{i0ij}$), if the assignment occurs after the vehicle has made a final turn towards D1; as illustrated in Figure \ref{fig:i0ij 2}. However, this error is found to be very insignificant in our numerical tests, and it is ignored for simplicity.

\subsubsection{State \texorpdfstring{$(i, 0, j, k)$}{(\textit{i}, 0, \textit{j}, \textit{k})}}
Similarly, a state-$(i, 0, j, k)$ vehicle could come from states $(i, 0, j, 0)$ or $(i, 0, k, 0)$, by being matched with an inter-zonal caller P2 on the way to deliver onboard inter-zonal seeker P1 (i.e., transitioning to state $(i, i, j, 0)$ or $(i, i, k, 0)$), and then picking up P2 at rate $p_{iij0\_k}$ or $p_{iik0\_j}$. If a state-$(i, 0, j, k)$ vehicle comes from state $(i, 0, j, 0)$, it always first delivers P1 (since zone $j$ is closer to zone $i$). Thus, it follows the original direction of P1's delivery trip, but starts from O2, to leave zone $i$; see Figure \ref{fig:i0jk 1}. The expected duration for this trip, denoted by $\tilde{T}_{i0jk}$, can be approximated by the expected duration from the blue dot to the zone boundary along the blue arrow. Similar analysis to that of state-$(i, 0, i, j)$ vehicles in Figure \ref{fig:i0ij 1} shows that:
\begin{equation}\label{eq:i0jk remaining}
\small
    \tilde{T}_{i0jk} = \frac{\Phi/v}{c_{i0j0} + p_{ii00\_j} + d_{i0ij}}\left[
    \frac{c_{i0j0}}{1 - e^{-\frac{a_{i0j0}}{n_{i0j0}}\frac{\Phi}{v}}} + \frac{(c_{i0ij} + p_{ii00\_j})/2}{1 - e^{-\frac{a_{i0j0}}{n_{i0j0}}\frac{\Phi}{2v}}} + \frac{(p_{iii0\_j} + p_{iij0\_i})/3}{1 - e^{-\frac{a_{i0j0}}{n_{i0j0}}\frac{\Phi}{3v}}}\right] - \frac{n_{i0j0}}{a_{i0j0}}.
\end{equation}
For a state-$(i, 0, j, k)$ vehicle coming from state $(i, 0, k, 0)$, the vehicle delivers P2 first. The feasible directions to leave zone $i$ for delivering P2, $\mathcal{V}_{ij}$, could be different from that of delivering P1, $\mathcal{V}_{ik}$; e.g., in Figure \ref{fig:i0jk 1}, the new zone-level movement to deliver P2 may require this vehicle to travel eastbound or westbound, instead of following P1's original delivery trip northbound.
In this case, the expected travel distance for this vehicle to leave zone $i$ is approximately $\frac{\Phi}{2}$. 

After leaving zone $i$, this vehicle stays in state $(i, 0, j, k)$ until reaching zone $j$ at the nearest boundary, with expected distance $L_{ij} - \Phi$ if $j\in G^{+}_i$, or $L_{ij}-\frac{3\Phi}{2}$ if $j\in G^{\times}_i$. Then, the expected travel distance in state $(i, 0, j, k)$ is the summation of the trip length within zone $i$ (i.e., $v\tilde{T}_{i0jk}$ or $\frac{\Phi}{2}$) and that between zone $i$ and zone $j$ (i.e., $L_{ij} - \Phi$ or $L_{ij}-\frac{3\Phi}{2}$), as summarized in Table \ref{tab:travel dist}.

Again, a state-$(i, 0, j, k)$ vehicle's trip to exit zone $i$ could be shorter than the remaining delivery trip of a state-$(i, 0, j, 0)$ vehicle (i.e., Eqn. \eqref{eq:i0jk remaining} overestimates $\tilde{T}_{i0jk}$), as the vehicle could go from O2 to leave zone $i$ directly, as illustrated in Figure \ref{fig:i0jk 2}.
This error is found to be very small and hence is ignored in the model.

\renewcommand{\arraystretch}{1.5}
\begin{table}[t]
\centering
\footnotesize
~\clap{
\begin{tabular}{>{\centering\arraybackslash}m{0.06\textwidth}>{\centering\arraybackslash}m{0.16\textwidth}>{\centering\arraybackslash}m{0.14\textwidth}m{0.34\textwidth}m{0.33\textwidth}}\hline
\multicolumn{1}{c}{\multirow{2}{*}{State}} & \multicolumn{1}{c}{\multirow{2}{0.13\textwidth}{\centering Vehicle Inflow}} & \multicolumn{1}{c}{\multirow{2}{0.14\textwidth}{\centering Expected Distance}} & \multicolumn{2}{c}{Delivery Trip}                   \\\cmidrule{4-5}
\multicolumn{1}{c}{}                       & \multicolumn{1}{c}{}                        & \multicolumn{1}{c}{} & \multicolumn{1}{c}{Start} & \multicolumn{1}{c}{End} \\\hline
\multirow{3}*{$(i, 0, i, i)$}  &  {$g_{i'0ii}, i'\in G^{+}_i$}  &  $\frac{5\Phi}{8}$    &   Random point on a boundary of zone $i$                    &     Closer of two random points in zone $i$                   \\
&  {$g_{i'0ii}, i'\in G^{\times}_i$}  &  $\frac{23\Phi}{30}$    &   A vertex of zone $i$                     &     Closer of two random points in zone $i$                   \\
& $p_{iii0\_i}$ & $\frac{\Phi}{2}$ & Random point in zone $i$ &  Closer of two random points that impose no detour  \\ \hline
\multirow{4}*{$(i, 0, i, j)$}  &  {$g_{i'0ij}, i'\in G^{+}_i\cap\tilde{\Omega}_{ji}$}  &  $\frac{5\Phi}{6}$    &  Random point on a boundary of zone $i$                    &     Random point in zone $i$                    \\
& {$g_{i'0ij}, i'\in G^{\times}_i\cap\tilde{\Omega}_{ji}$}  &  $\Phi$    &  A vertex of zone $i$                    &     Random point in zone $i$                    \\
& $p_{iii0\_j}$ & $v\tilde{T}_{i0ij}$ & Random point in zone $i$ &  Random point in zone $i$  \\
& $p_{iij0\_i}$ & $\frac{2\Phi}{3}$ & Random point in zone $i$ &  Random point in zone $i$  \\
 \hline
\multirow{4}*{$(i, 0, j, k)$}  &  {$p_{iij0\_k}, j\in G^{+}_i$}  &  $v\tilde{T}_{i0jk} + L_{ij} - \Phi$ & Random point in zone $i$                    &   Nearest boundary of zone $j$ \\
& {$p_{iij0\_k}, j\in G^{\times}_i$} & $v\tilde{T}_{i0jk} + L_{ij} - \frac{3\Phi}{2}$ & Random point in zone $i$ &  Nearest boundary of zone $j$\\
& {$p_{iik0\_j}, j\in G^{+}_i$} & $L_{ij} - \frac{\Phi}{2}$ & Random point in zone $i$ &  Nearest boundary of zone $j$\\
& {$p_{iik0\_j}, j\in G^{\times}_i$} & $L_{ij} - \Phi$ & Random point in zone $i$ &  Nearest boundary of zone $j$\\  \hline
\multirow{3}*{$(i, 0, i, 0)$}  &  $c_{i0i0}$  &  $\frac{5\Phi}{6}$    &  Random point on a boundary of zone $i$                     &     Random point in zone $i$                    \\ 
& $p_{ii00\_i},\ c_{i0ii}$ & $\frac{2\Phi}{3}$ & Random point in zone $i$ &  Random point in zone $i$  \\
 & $p_{iii0\_i}$ & $\frac{\Phi}{2}$ & Closer of two random interior points that impose no detour & Farther of two random interior points that impose no detour  \\  \hline
\multirow{3}*{$(i, 0, j, 0)$}  &  {$c_{i0j0}$}  &  $\Phi$ & Random point on a boundary of zone $i$                    &   Another boundary of zone $i$ \\
 & $p_{ii00\_j}$, $c_{i0ij}$ & $\frac{\Phi}{2}$ & Random point in zone $i$ &  A boundary of zone $i$\\
& $p_{iii0\_j}$, $p_{iij0\_i}$ & $\frac{\Phi}{3}$ & Random point in zone $i$ &  A boundary of zone $i$\\\hline
\end{tabular}
}
\caption{Expected vehicle travel distance to deliver onboard passengers.
}
\label{tab:travel dist}
\end{table}
\renewcommand{\arraystretch}{1}

\newpage
\nolinenumbers
\bibliographystyle{elsarticle-harv}\biboptions{authoryear}
\bibliography{references}

\end{document}